\documentclass[12pt]{article}
\renewcommand{\theequation}{\arabic{section}.\arabic{equation}}

\usepackage{graphicx}   %
\usepackage{enumerate}  %
\usepackage{latexsym}   %
\usepackage{amsmath}    %
\usepackage{amsbsy}     %
\usepackage{amsfonts}   %
\usepackage{mathrsfs}       
\usepackage{amssymb}    
\usepackage{array}      %
\usepackage{theorem}
\usepackage{comment}

\makeatletter
\@addtoreset{equation}{section}%
\renewcommand{\theequation}{\thesection.\@arabic\c@equation}
\makeatother

\usepackage[square,numbers,comma,sort&compress]{natbib}

\usepackage[colorlinks=true]{hyperref}
\hypersetup{urlcolor=blue, citecolor=red}

\newcommand{\qed}{\hbox{\rule{6pt}{6pt}}}

 \newtheorem{thm}{Theorem}
 \newtheorem{mainTh}{Main Theorem}
 
 \newtheorem{lem}{Lemma}[section]
 
 \theoremstyle{definition}
 \newtheorem{defn}{Definition}[section]
 \theoremstyle{remark}
 \newtheorem{rem}{Remark}[section]
 
 \numberwithin{equation}{section}
 \newtheorem{notn}{Notation}

\def\R{{\mathbb R}}  
\def\N{{\mathbb{N}}}
\def\Z{{\mathbb{Z}}}

\def\ds{\displaystyle}

\begin{document}
\vspace*{1.5cm}
\begin{center}
\textbf{\Large 
PHASE-FIELD SYSTEMS FOR GRAIN \\[0.1cm] 
BOUNDARY MOTIONS UNDER \\[0.1cm] 
ISOTHERMAL SOLIDIFICATIONS\\[1.5cm] 
}
\end{center}
\begin{center}
\textsc{Ken Shirakawa}
\\
Department of Mathematics, Faculty of Education, Chiba University
\\
1-33 Yayoi-cho, Inage-ku, Chiba 263-8522, Japan
\\
(sirakawa@faculty.chiba-u.jp)
\vspace{0.75cm}

\textsc{Hiroshi Watanabe}
\\
Department of General Education, Salesian Polytechnic 
\\
4-6-8 Oyamagaoka, Machida, Tokyo 194-0215, Japan
\\
(h-watanabe@salesio-sp.ac.jp)
\vspace{0.75cm}

\textsc{Noriaki Yamazaki}
\\
Department of Mathematics, Faculty of Engineering, Kanagawa University
\\
3-27-1 Rokkakubashi, Kanagawa-ku, Yokohama 221-8686, Japan
\\
(noriaki@kanagawa-u.ac.jp)
\end{center}
\vspace{1.5cm}

\noindent
\textbf{Abstract.}
Two main existence theorems are proved for two nonstandard systems of parabolic 
initial-boundary value problems.
The systems are based on the ``$ \phi $-$ \eta $-$ \theta $ model'' proposed by Kobayashi [RIMS K\^{o}ky\^{u}roku, \textbf{1210} (2001), 68--77] as a phase-field model of planar grain boundary motion under isothermal solidification. Although each of the systems has specific characteristics and mathematical difficulties, the proofs of the main theorems are based on the time discretization method by means of a common approximating problem. As a consequence, we provide a uniform solution method for a wide scope of parabolic systems associated with the $ \phi $-$ \eta $-$ \theta $ model.

\vfill
\noindent
------------------------------------------------------------
\\
{\footnotesize
This work is supported by Grants-in-Aid 24740099 (K.~S.), 25800086 and 26400138 (H.~W.), and 26400179 (N.~Y.) from the Japan Society for the Promotion of Science.
}
\newpage
\section*{\hypertarget{Intro}{Introduction}}
\ \vspace{-3ex}

Let $ 0 < T < \infty $ be a constant, $ 1 < N \in \N $ have a fixed value, and $ \Omega \subset \R^N $ be a bounded domain with a smooth boundary $ \partial \Omega $. We denote by $ \nu_{\partial \Omega} $ the unit outer normal vector on $ \partial \Omega $, and we set $ Q := (0, T) \times \Omega $ and $ \Sigma := (0, T) \times \partial \Omega $. 

Further, let $ \nu \geq 0 $ and $ u \in \R $ be constants. In this paper, two themes, concerning two nonstandard systems of parabolic variational inequalities, are addressed. In the first, we assume $ \nu > 0 $ and consider the following coupled system of parabolic type initial-boundary value problems, denoted by $ \mbox{(\hyperlink{(S;u)_nu}{S})}_\nu $. 
\begin{eqnarray}
&& \hspace{-5ex} \mbox{(\hypertarget{(S;u)_nu}{S})}_\nu:
\nonumber
\\
&& \hspace{-2ex} \left\{ ~ \parbox{14.5cm}{$ w_t -{\mit \Delta} w +\partial \gamma(w) +g_w(w, \eta) +\alpha_w(w, \eta)|\nabla \theta| +\nu \beta_w(w, \eta) |\nabla \theta|^2 \ni 0 $ \ in $ Q $,
\\[1ex]
$ \nabla w \cdot \nu_{\partial \Omega} = 0 $ \ on $ \Sigma $,
\\[1ex]
$ w(0, x) = w_0(x) $, \ $ x \in \Omega $;
} \right.
\label{1st.eq}
\\[1ex]
&& \hspace{-2ex} \left\{ ~ \parbox{14.5cm}{$ \eta_t -{\mit \Delta} \eta +g_\eta(w, \eta) +\alpha_\eta(w, \eta) |\nabla \theta| +\nu \beta_\eta(w, \eta) |\nabla \theta|^2 = 0 $ \ in $ Q $,
\\[1ex]
$ \nabla \eta \cdot \nu_{\partial \Omega} = 0 $ \ on $ \Sigma $,
\\[1ex]
$ \eta(0, x) = \eta_0(x) $, \ $ x \in \Omega $;
} \right.
\label{2nd.eq}
\\[1ex]
&& \hspace{-2ex} \left\{ ~ \parbox{14.5cm}{$ \displaystyle \alpha_0(w, \eta) \, \theta_t -{\rm div} \left( \alpha(w, \eta) \frac{\nabla \theta}{|\nabla \theta|} +2 \nu \beta(w, \eta) \nabla \theta \right) = 0 $ \ in $ Q $,
\\[1ex]
$ \displaystyle \left( \alpha(w, \eta) \frac{\nabla \theta}{|\nabla \theta|} +2 \nu \beta(w, \eta) \nabla \theta \right) \cdot \nu_{\partial \Omega} = 0 $ \ on $ \Sigma $,
\\[1ex]
$ \theta(0, x) = \theta_0(x) $, \ $ x \in \Omega $.
} \right.
\label{3rd.eq}
\end{eqnarray}
Here, $ w_0 = w_0(x) $, $ \eta_0 = \eta_0(x)$, and $ \theta_0 = \theta_0(x) $ are given initial data on $ \Omega $, $ \partial \gamma $ is the subdifferential of a proper lower semi-continuous (l.s.c.)\ and convex function $ \gamma = \gamma(w) $ on $ \R $, and $ g({}\cdot{}) = g(w, \eta) $, $ \alpha_0 = \alpha_0(w, \eta) $, $ \alpha = \alpha(w, \eta) $, and $ \beta = \beta(w, \eta) $ are given real-valued functions. The subscripts ``$ \empty_w $'' and ``$ \empty_\eta $'' denote differentials with respect to the corresponding variables.

The system $\mbox{(\hyperlink{(S;u)_nu}{S})}_\nu$ is based on the ``$ \phi $-$ \eta $-$ \theta $ model'' proposed by Kobayashi \citep{K} as a mathematical model of planar grain boundary motion under an isothermal solidification. Since this model was presented as an advanced version of the ``Kobayashi-Warren-Carter model'' of grain boundary motion, proposed by Kobayashi et al.\ \citep{KWC1,KWC2}, our themes are related to the previous work (e.g., \citep{GG,IKY08,IKY09,IKY12,KY,KG,KWC1,KWC2,MS,SW,SWY,WS13,WS14}) associated with the Kobayashi-Warren-Carter model. 
According to the modeling method of \citep{K}, $ \mbox{(\hyperlink{(S;u)_nu}{S})}_\nu $ is derived as a gradient system of a governing free energy, defined as follows:
\begin{align}
[w, \eta, \theta] \in & \; [H^1(\Omega) \cap L^\infty(\Omega)] \times [H^1(\Omega) \cap L^\infty(\Omega)] \times H^1(\Omega)
\nonumber
\\[1ex]
\mapsto & \; \mathscr{F}_\nu(w, \eta, \theta) := \frac{1}{2} \int_\Omega |\nabla w|^2 \, dx +\frac{1}{2} \int_\Omega |\nabla \eta|^2 \, dx +\int_\Omega \gamma(w) \, dx  
\label{freeEnergy}
\\
& +\int_\Omega g(w, \eta) \, dx +\int_\Omega \alpha(w, \eta) |\nabla \theta| \, dx +\nu \int_\Omega \beta(w, \eta)|\nabla \theta|^2 \, dx.
\nonumber
\end{align}

In this context, the constant $ u $ is the relative temperature with critical degree 0, and the unknown $ w = w(t, x) $ is an order parameter to indicate the solidification order of the polycrystal. The unknowns $ \eta = \eta(t, x) $ and $ \theta = \theta(t, x) $ are components of the vector field
$$
(t, x) \in Q \mapsto \eta(t, x) \left[ \rule{0pt}{10pt} \cos \theta(t, x), \sin \theta(t, x) \right] \in \mathbb{R}^2,
$$
which was adopted in \citep{KWC1,KWC2} as a vectorial phase field to reproduce the crystalline orientation in $ Q $. Here, the components $ \eta $ and $ \theta $ are order parameters to indicate, respectively, the orientation order and angle of the grain. In particular, $ w $ and $ \eta $ are taken to satisfy the constraints $ 0 \leq w, \eta \leq 1 $ in $ Q $, and the cases $ [w, \eta] \approx [1, 1] $ and $ [w, \eta] \approx [0, 0] $ are respectively assigned to ``the solidified-oriented phase'' and ``the liquefied-disoriented phase,'' which correspond to physically meaningful phases.
Hence, we suppose that
\begin{description}
\item[\textmd{(\hypertarget{g0}{g0})}](double-well graph) the function 
\begin{equation*}
[w, \eta] \in \R^2 \mapsto G(w, \eta) := \gamma(w) +g(w, \eta) \in (-\infty, \infty]
\end{equation*}
may have just two minimums, around $ [1, 1] $ and $ [0, 0] $, and moreover, if the temperature $ u $ is sufficiently lower than (resp. higher than) the critical degree, then this function has a unique minimizer around $ [1, 1] $ (resp. $ [0, 0] $) (cf.\ \citep{AC,Cag,Fix,Visintin}).
\end{description}
In light of (\hyperlink{g0}{g0}) and the existing phase transition models (e.g., \citep{AC,Cag,CL,CS,Fix,HSZ,Kenmochi01,PF,SIYK,SZ,Visintin}), we can consider the following settings as possible expressions of the double-well functions:
\begin{description}
\item[\textmd{(\hypertarget{g1}{g1})}](standard polynomial)
\begin{equation*}
\left\{ ~ 
\parbox{11cm}{$ \gamma(w) := 0 $
\\[1ex]
$ \displaystyle g(w, \eta) := c \left[ \frac{1}{4}w^2 (w -1)^2 -u w^2 \left( \frac{w}{3} -\frac{1}{2} \right) \right] +\frac{1}{2} (w -\eta)^2 $}
\right. \mbox{for all $ w, \eta \in \R $,}
\end{equation*}
and therefore $ G({}\cdot{}) = g({}\cdot{}) $ on $ \R^2 $ (cf.\ \citep{AC,Cag,Fix,PF,Visintin});
\vspace{-1ex}
\item[\textmd{(\hypertarget{g2}{g2})}](logarithmic constraint)
\begin{equation*}
\left\{ ~ 
\parbox{8cm}{$ \ds \gamma(w) := \displaystyle{\frac{1}{2}} \left( \rule{0pt}{10pt} w \log w +(1 -w) \log (1 -w) \right) $
\\[0.5ex]
\mbox{$\quad$} with $ \gamma(0) = \gamma(1) := 1 $
\\[1ex]
$ \ds g(w, \eta) := -\frac{c}{2} \left(w -u -\frac{1}{2} \right)^2 +\frac{1}{2} (w -\eta)^2 $}
\right. \mbox{ \ for all $ w, \eta \in \R $,}
\end{equation*}
and therefore the range of $ w $ is constrained to the open interval $ (0, 1) $ (cf.\ \citep{HSZ,SZ});
\item[\textmd{(\hypertarget{g3}{g3})}](non-smooth constraint)
\begin{equation*}
\left\{ ~ 
\parbox{8.25cm}{$ \ds \gamma(w) := I_{[0, 1]}(w) $ 
\\[1ex]
$ \ds g(w, \eta) := -\frac{c}{2} \left(w -u -\frac{1}{2} \right)^2 +\frac{1}{2} (w -\eta)^2 $}
\right. \mbox{for all $ w, \eta \in \R $,}
\end{equation*}
and therefore $ w $ is constrained to the compact interval $ [0, 1] $ (cf.\ \citep{CL,CS,Kenmochi01,SIYK,Visintin}).
\end{description}
Here, $ c > 0 $ is a constant, and for any $ K \subset \R $, $ I_{K} $ denotes the indicator function on $ K $, i.e.,
\begin{equation*}
\tau \in \R \mapsto I_{K}(\tau) := \left\{ \begin{array}{ll}
0, & \mbox{if $ \tau \in K $,}
\\[1ex]
\infty, & \mbox{otherwise.}
\end{array} \right.
\end{equation*}

Kobayashi \citep{K} adopted a setting such that
\begin{equation}\label{possSet}
\left\{ \hspace{-3ex}
\parbox{12.75cm}{
\vspace{-1.5ex}
\begin{itemize}
\item the functions $ \gamma $ and $ g({}\cdot{}) $ are given in accordance with (\hyperlink{g1}{g1}),
\vspace{-1ex}
\item $ \ds \alpha_0(w, \eta) = \alpha(w, \eta) := \eta^2/2 $ \ and \  $ \ds \beta(w, \eta) := w^2/2 $ \ \ for $ [w, \eta] \in \R^2 $.
\end{itemize}
\vspace{-1.5ex}
}
\right.
\end{equation}
Applying this, the original profile of the $ \phi $-$ \eta $-$ \theta $ model in \citep{K} is described by
\begin{eqnarray}
&& \hspace{-5ex} w_t -{\mit \Delta} w +c w(w -1) \left( w -u -{\frac{1}{2}} \right) +(w -\eta) +\nu w|\nabla \theta|^2 = 0 \mbox{ \ in $ Q $,}
\label{org.1st}
\\[1ex]
&& \hspace{-5ex} \eta_t -{\mit \Delta} \eta +(\eta -w) +\eta|\nabla \theta| = 0 \mbox{ \ in $ Q $,}
\label{org.2nd}
\\[1ex]
&& \hspace{-5ex} \eta^2 \theta_t -{\rm div}\left( \eta^2 \frac{\nabla \theta}{|\nabla \theta|} +2 \nu w^2 \nabla \theta \right) = 0 \mbox{ \ in $ Q $,}
\label{org.3rd}
\end{eqnarray}
with the initial-boundary conditions as in (\ref{1st.eq})--(\ref{3rd.eq}).
\medskip

From a mathematical point of view, there do not appear to be great differences between (\ref{1st.eq}) and (\ref{2nd.eq}). However, from the original profiles (\ref{org.1st})--(\ref{org.3rd}), it can be seen that (\ref{2nd.eq}) corresponds to the equation for the mobilities of grain boundaries (interfaces) as in Kobayashi-Warren-Carter \citep{KWC1}, while (\ref{1st.eq}) is an Allen-Cahn type equation to reproduce ``interfacial diffusions,'' as in the models of phase transitions. 
\medskip

Next, in our second theme, we consider a limiting system $ \mbox{(\hyperlink{(S;u)_nu}{S})}_\nu $ as $ \nu \searrow 0 $, similarly to the case with $ \nu = 0 $. This is denoted by $ \mbox{(\hyperlink{(S;u)_0}{S})}_0 $ and formally described as follows. 
\begin{eqnarray}
\hspace{-0.5cm}\mbox{(\hypertarget{(S;u)_0}{S})}_0: \hspace{1.5cm} &&
\nonumber
\\[2ex]
&& \hspace{-0.75cm}\left\{ ~ \parbox{11.5cm}{$ w_t -{\mit \Delta} w +\partial \gamma(w) +g_w(w, \eta) +\alpha_w(w, \eta)|D \theta| \ni 0 $ \ in $ Q $,
\\[1ex]
$ \nabla w \cdot \nu_{\partial \Omega} = 0 $ \ on $ \Sigma $,
\\[1ex]
$ w(0, x) = w_0(x) $, ~ $ x \in \Omega $;
} \right.
\label{1st.eq0}
\\[2ex]
&& \hspace{-0.75cm}\left\{ ~ \parbox{11.5cm}{$ \eta_t -{\mit \Delta} \eta +g_\eta(w, \eta) +\alpha_\eta(w, \eta) |D \theta| = 0 $ \ in $ Q $,
\\[1ex]
$ \nabla \eta \cdot \nu_{\partial \Omega} = 0 $ \ on $ \Sigma $,
\\[1ex]
$ \eta(0, x) = \eta_0(x) $, ~ $ x \in \Omega $;
} \right.
\label{2nd.eq0}
\\[2ex]
&& \hspace{-0.75cm}\left\{ ~ \parbox{11.5cm}{$ \displaystyle \alpha_0(w, \eta) \, \theta_t -{\rm div} \left( \alpha(w, \eta) \frac{D \theta}{|D \theta|} \right) = 0 $ \ in $ Q $,
\\[1ex]
$ \displaystyle \alpha(w, \eta) \frac{D \theta}{|D \theta|} \cdot \nu_{\partial \Omega} = 0 $ \ on $ \Sigma $,
\\[1ex]
$ \theta(0, x) = \theta_0(x) $, ~ $ x \in \Omega $.
} \right.
\label{3rd.eq0}
\end{eqnarray}
Taking these together, we set the functional
\begin{align}
[w, \eta, \theta] 
\in & \; [H^1(\Omega) \cap L^\infty(\Omega)] \times [H^1(\Omega) \cap L^\infty(\Omega)] \times BV(\Omega)
\nonumber
\\[1ex]
\mapsto & \;\mathscr{F}_0(w, \eta, \theta) := \frac{1}{2} \int_\Omega |\nabla w|^2 \, dx +\frac{1}{2} \int_\Omega |\nabla \eta|^2 \, dx +\int_\Omega \gamma(w) \, dx  
\label{freeEnergy0}
\\
& +\int_\Omega g(w, \eta) \, dx +\int_\Omega \alpha(w, \eta) |D \theta|
\nonumber
\end{align}
as the corresponding free energy in $ \mbox{(\hyperlink{(S;u)_0}{S})}_0 $. 
\medskip

Because of the absence of the term $ \nu \beta(w, \eta) $, the limiting system $ \mbox{(\hyperlink{(S;u)_0}{S})}_0 $ may appear to be a simplified version of $ \mbox{(\hyperlink{(S;u)_nu}{S})}_\nu $ when $ \nu > 0 $. However, it must be noted that the definition (\ref{freeEnergy0}) of the limiting free energy $ \mathscr{F}_0 $ in this system includes a nontrivial term $ \int_\Omega \alpha(w, \eta) |D \theta| $. The spatial gradients $ D \theta $ in (\ref{1st.eq0})--(\ref{3rd.eq0}) necessitate much more delicate mathematical treatment than $ \nabla \theta $ in (\ref{1st.eq})--(\ref{3rd.eq}).
\medskip

Now, the objective of this paper is to establish a uniform solution method for the systems $ \mbox{(\hyperlink{(S;u)_nu}{S})}_\nu $ for all $ \nu \geq 0 $, including their relaxed versions. We here adopt an analytical approach based on time discretization and set the goal to prove two main theorems:
\begin{description}
\item[Main Theorem \ref{MainTh01}.]There exists a solution to $ \mbox{(\hyperlink{(S;u)_nu}{S})}_\nu $ for each fixed $ \nu > 0 $.
\vspace{-1ex}
\item[Main Theorem \ref{MainTh02}.]There exists a solution to $ \mbox{(\hyperlink{(S;u)_0}{S})}_0 $.
\end{description}

The plan of this paper is as follows: In the next section, we set forth some specific notation. In Section~2, we state our two main theorems with proper definitions of the solutions to the respective systems. In Section~3, we present approximating problems for our systems and supply some auxiliary lemmas aimed at the method of obtaining the approximating solutions. The approximating problems are provided in the forms of time discretization of $ \mbox{(\hyperlink{(S;u)_nu}{S})}_\nu $ for $ \nu > 0 $, and the existence and uniqueness of the approximating solutions are proved in the following Section 4. Sections 5 and 6 are devoted to the proofs of Main Theorems 1 and 2, respectively. Finally, we add an Appendix to make supplementary statements for some preliminary facts and the solutions to our systems.

\section{Preliminaries}
\ \ \vspace{-3ex}

First we elaborate the notation used throughout. 
\begin{notn}[real analysis]\label{Note00}
\begin{em}
For arbitrary $ a_0, b_0 \in [-\infty, \infty] $, we define
\begin{equation*}
a_0 \vee b_0 := \max \{ a_0, b_0\} \mbox{ \ and \ } a_0 \wedge b_0 := \min \{ a_0, b_0\}, 
\end{equation*}
and for arbitrary $ -\infty \leq a \leq b \leq \infty $, we define the truncation function (operator) $ \mathscr{T}_a^b : \R \rightarrow [a, b] $ by letting
\begin{equation*}
r \in \R \mapsto \mathscr{T}_a^b r := a \vee (b \wedge r) \in [a, b].
\end{equation*}

Let $ d \in \N $ take any fixed value. We denote by $ |x| $ and $ x \cdot y $ the Euclidean norm of $ x \in \mathbb{R}^d $ and the standard scalar product of $ x, y \in \R^d $, respectively, as usual, i.e.,
\begin{equation*}
\begin{array}{c}
| x | := \sqrt{x_1^2 +\cdots +x_d^2} \mbox{ \ and \ } x \cdot y  := x_1 y_1 +\cdots +x_d y_d 
\\[1ex]
\mbox{ for all $ x = [x_1, \ldots, x_d], ~ y = [y_1, \ldots, y_d] \in \mathbb{R}^d $.}
\end{array}
\end{equation*}
For arbitrary $ x = [x_1, \ldots, x_d] \in \R^d $ and $ y = [y_1, \ldots, y_d] \in \mathbb{R}^d $, we say that $ x \leq y $ or $ y \geq x $ if $ x_k \leq y_k $, for $ k = 1, \ldots, d $.

The $ d $-dimensional Lebesgue measure is denoted by $ \mathscr{L}^d $. Also, unless otherwise specified, the measure-theoretic phrases such as ``a.e.,'' ``$ dt $,'' ``$ dx $'', and so on, are with respect to the Lebesgue measure in each corresponding dimension.
For a (Lebesgue) measurable function $ f : B \rightarrow [-\infty, \infty] $ on a Borel subset $ B \subset \R^d $, we denote by $ [f]^+ $ and $ [f]^- $, respectively, the positive and negative parts of $ f $, i.e.,
\begin{equation*}
[f]^+(x) := \mathscr{T}_0^\infty f(x) \mbox{ \ and \ } [f]^-(x) := - \mathscr{T}_{-\infty}^0 f(x), \mbox{ a.e.\ $ x \in B $.}
\end{equation*}
\end{em}
\end{notn}

\begin{notn}[abstract functional analysis]\label{Note01}
\begin{em}
For an abstract Banach space $ X $, we denote by $ |\mathrel{\cdot}|_X $ the norm of $ X $, and when $ X $ is a Hilbert space, we denote by $ (\,\cdot\,, \,\cdot\,)_X $ its inner product. For a subset $ A $ of a Banach space $ X $, we denote by $ {\rm int}(A) $ and $ \overline{A} $ the interior and the closure of $ A $, respectively. 

Fix $ 1 < d \in \N $. Then, for a Banach space $ X $ the topology of the product Banach space
\begin{center}
$ X^d := \overbrace{X \times \cdots \times X}^{\mbox{\scriptsize $ d $ times}} $
\end{center}
has the norm
\begin{equation*}
|z|_{X^d} := \sum_{k = 1}^d |z_k|_{X}, \mbox{ for $ z = [z_1, \ldots, z_d] \in X^d $.}
\end{equation*}
However, if $ X $ is a Hilbert space, then the topology of the product Hilbert space $ X^d $ has the inner product
\begin{equation*}
(z, \tilde{z})_{X^d} := \sum_{k = 1}^d (z_k, \tilde{z}_k)_X, \mbox{ for $ z = [z_1, \ldots, z_d] \in X^d $ and $ \tilde{z} = [\tilde{z}_1, \ldots, \tilde{z}_d] \in X^d $,}
\end{equation*}
and, hence the norm in this case is provided by
\begin{equation*}
|\zeta|_X := \sqrt{(z, z)_{X^d}} = \left( \rule{-2pt}{18pt} \right. \sum_{k = 1}^d |z_k|_X^2 \left. \rule{-2pt}{18pt} \right)^{\hspace{-0.5ex}1/2}, \mbox{ for $ z = [z_1, \ldots, z_d] \in X^d $.}
\end{equation*}
For a Banach space $ X $, we denote the dual space by $ X^* $. For a single-valued operator $ \mathscr{A} : X \rightarrow X^* $, we write
\begin{equation*}
\mathscr{A} z = [\mathscr{A} z_1, \ldots, \mathscr{A} z_d] \in [X^*]^d \mbox{ \ for any $ z = [z_1, \ldots, z_d] \in X^d $.}
\end{equation*}

For any proper lower semi-continuous (l.s.c.\ hereafter) and convex function $ \Psi $ defined on a Hilbert space $ X $, we denote by $ D(\Psi) $ its effective domain and by $ \partial \Psi $ its subdifferential. The subdifferential $ \partial \Psi $ is a set-valued map corresponding to a weak differential of $ \Psi $, and it has a maximal monotone graph in the product Hilbert space $ X^2 $. More precisely, for each $ z_0 \in X $, the value $ \partial \Psi(z_0) $ is defined as the set of all elements $ z_0^* \in X $ that satisfy the variational inequality
\begin{equation*}
(z_0^*, z -z_0)_X \leq \Psi(z) -\Psi(z_0) \mbox{ \ for any $ z \in D(\Psi) $,}
\end{equation*}
and the set $ D(\partial \Psi) := \{z \in X \mid \partial \Psi(z) \ne \emptyset\} $ \
is called the domain of $ \partial \Psi $. We often use the notation ``$ [z_0, z_0^*] \in \partial \Psi $ in $ X^2 $\,'' to mean ``$ z_0^* \in \partial \Psi(z_0) $ in $ X $ with $ z_0 \in D(\partial \Psi) $,'' by identifying the operator $ \partial \Psi $ with its graph in $ X^2 $.
\end{em}
\end{notn}
\begin{rem}\label{Rem.Time-dep.}
\begin{em}
It is often useful to consider the subdifferentials under time-dependent settings. In this regard, several general theories have been established by previous researchers (e.g., Kenmochi \citep{Kenmochi}, and \^{O}tani \citep{Otani}). From these (e.g., \citep[Chapter 2]{Kenmochi}), one can see the following fact: 
\begin{description}
\item[{(\hypertarget{Fact1}{Fact\,1})}]Let $ E_0 $ be a convex subset in a Hilbert space $ X $, let $ I \subset [0, \infty) $ be a time interval, and for any $ t \in I $, let $ \Psi^t : X \rightarrow (-\infty, \infty] $ be a proper l.s.c.\ and convex function such that $ D(\Psi^t) = E_0 $ for all $ t \in I $. Based on this, define a convex function $ {\Psi}^I : L^2(I; X) \rightarrow (-\infty, \infty] $, by setting
\begin{equation*}
\zeta \in L^2(I; X) \mapsto {\Psi}^I(\zeta) := \left\{ \begin{array}{ll}
\multicolumn{2}{l}{\ds \int_I \Psi^t(\zeta(t)) \, dt, \mbox{ if $ \Psi^{(\cdot)}(\zeta) \in L^1(I) $,}}
\\[1ex]
\infty, & \mbox{otherwise.}
\end{array} \right.
\end{equation*}
Here, if $ E_0 \subset D({\Psi}^I) $, and the function $ t \in I \mapsto \Psi^t(z) $ is integrable for any $ z \in E_0 $, then the following holds:
\begin{equation*}
\begin{array}{c}
[\zeta, \zeta^*] \in \partial {\Psi}^I \mbox{ \ in $ L^2(I; X)^2 $ if and only if}
\\[1ex]
\zeta \in D({\Psi}^I) \mbox{ and } [\zeta(t), \zeta^*(t)] \in \partial \Psi^t \mbox{ in $ X^2 $, a.e.\ $ t \in I $.}
\end{array}
\end{equation*}
\end{description}
\end{em}
\end{rem}
\begin{notn}[basic elliptic operators]\label{Note02}
\begin{em}
Let $ F : H^1(\Omega) \rightarrow H^1(\Omega)^* $ be the duality mapping, defined as
\begin{equation*}
\begin{array}{c}
\langle F \varphi, \psi \rangle_{*} := (\varphi, \psi)_{H^1(\Omega)} = (\varphi, \psi)_{L^2(\Omega)} +(\nabla \varphi, \nabla \psi)_{L^2(\Omega)^N}
\\[1ex]
\mbox{ \ for all $ \varphi, \psi \in H^1(\Omega) $,}
\end{array}
\end{equation*}
where $ \langle \, \cdot \,, \,\cdot\, \rangle_{*}  $ is the duality pairing between $ H^1(\Omega) $ and its dual $ H^1(\Omega)^* $. 

Let $ {\mit \Delta}_{\rm N} $ be the Laplacian operator subject to the zero Neumann boundary condition, i.e.,
\begin{equation*}
{\mit \Delta}_{\rm N} : z \in D_{\rm N}:= \left\{ \begin{array}{l|l}
z \in H^2(\Omega) & \nabla z \cdot \nu_{\partial \Omega} = 0 \mbox{ in $ L^2(\partial \Omega) $}
\end{array} \right\} \subset L^2(\Omega) \mapsto {\mit \Delta} z \in L^2(\Omega). 
\end{equation*}
As is well known, 
\begin{equation}\label{F}
F z = -{\mit \Delta}_{\rm N} z +z \mbox{ in $ L^2(\Omega) $ if $ z \in D_{\rm N} $.}
\end{equation}
\end{em}
\end{notn}

\begin{rem}\label{Note03}
\begin{em}
We here show a representative example of the subdifferential. Let $ d \in \N $ be fixed, and let $ V_{\rm D}^d : L^2(\Omega)^d \rightarrow [0, \infty] $ be a proper l.s.c.\ and convex function of the so-called Dirichlet type integral, i.e.,
\begin{equation}\label{Dirichlet^d}
z \in L^2(\Omega)^d \mapsto V_{\rm D}^d(z) := \left\{ \begin{array}{ll}
\multicolumn{2}{l}{\ds \frac{1}{2} \int_\Omega |\nabla z|_{\R^{d \times N}}^2 \, dx, \mbox{ \ if $ z \in H^1(\Omega)^d $,}}
\\[2ex]
\infty, & \mbox{otherwise.}
\end{array} \right.
\end{equation}
Then, with regard to the subdifferential $ \partial V_{\rm D}^d \subset [L^2(\Omega)^d]^2 $, it is known (see, e.g., \citep{Barbu} or \citep{Brezis}) that
\begin{equation}\label{elliptic01}
z \in L^2(\Omega)^d \mapsto \partial V_{\rm D}^d(z) = \left\{ \begin{array}{ll}
\multicolumn{2}{l}{\{ -{\mit \Delta}_{\rm N} z \}, \mbox{ if $ z \in D_{\rm N}^d $,}}
\\[1ex]
\emptyset, & \mbox{otherwise.}
\end{array} \right.
\end{equation}
In this light, $ \partial V_{\rm D}^d $ and $ -{\mit \Delta}_{\rm N} $ are identified as the maximal monotone graphs in $ [L^2(\Omega)^d]^2 $. 
\end{em}
\end{rem}

\begin{notn}[BV theory; cf.\ \citep{AFP,ABM,EG,G}]\label{Note04}
\begin{em}
Let $ d \in \N $, and let $ U \subset \mathbb{R}^d $ be an open set. We denote by $ \mathcal{M}(U) $ the space of all finite Radon measures on $ U $. The space $ \mathcal{M}(U) $ is known as the dual space of the Banach space $ C_0(U) $, i.e., $ \mathcal{M}(U) = C_0(U)^* $, where $ C_0(U) $ denotes the closure of the space $ C_{\rm c}(U) $ of all continuous functions having compact supports, in the topology of $ C(\overline{U}) $.

A function $ z \in L^1(U) $ is called a \textit{function of bounded variation} on $ U $ (or simply $ z \in BV(U)$)
if and only if its distributional gradient $ Dz $ is a finite Radon measure on $ U $, namely, $ Dz \in \mathcal{M}(U)^d $. Here, for any $ z \in BV(U) $ the Radon measure $ Dz $ is called the \textit{variation measure} of $ z $, and its total variation $ |Dz| $ is similarly the \textit{total variation measure} of $ z $. Additionally,
\begin{equation*}
|Dz|(U) = \sup \left\{ \begin{array}{l|l}
\displaystyle \int_U z \, {\rm div} \, \varphi \, dx & \varphi \in C_{\rm c}^{1}(U)^N \mbox{ \ and \ } | \varphi | \leq 1 \mbox{ on $ {U} $}
\end{array} \right\}.
\end{equation*}
The space $ BV(U) $ is a Banach space, with the norm
\begin{equation*}
|z|_{BV(U)} := |z|_{L^1(U)} +|Dz|(U) \mbox{ \ for any $ z \in BV(U) $.}
\end{equation*}
Additionally, we say that $ z_n \to z $ ``weakly-$ * $'' in $ BV(U) $ if $ z \in BV(U) $, $ \{ z_n \, | \, n \in \N \} \subset BV(U) $, $ z_n \to z $ in $ L^1(U) $, and $ Dz_n \to Dz $ weakly-$ * $ in $ \mathcal{M}(U) $ as $ n \to \infty $.

The space $ BV(U) $ has another topology, called ``strict topology,'' which has the following distance (cf.\ \citep[Definition 3.14]{AFP}):
\begin{equation*}
[\varphi, \psi] \in BV(U)^2 \mapsto |\varphi -\psi|_{L^1(U)} +\bigl| |D\varphi|(U) -|D\psi|(U) \bigr|.
\end{equation*}
In this regard, we say that $ z_n \to z $ \textit{strictly} in $ BV(U) $ if $ z \in BV(U) $, $ \{ z_n \, | \, n \in \N \} \subset BV(U) $, $ z_n \to z $ in $ L^1(U) $, and $ |Dz_n|(U) \to |Dz|(U) $ as $ n \to \infty $. 

Specifically, when the boundary $ \partial U $ is Lipschitz, the Banach space $ BV(U) $ is continuously embedded into $ L^{d/(d -1)}(U) $ and compactly embedded into $ L^p(U) $ for any $ 1 \leq p < d / (d -1) $ (cf.\ \citep[Corollary 3.49]{AFP} or \citep[Theorems 10.1.3--10.1.4]{ABM}). Additionally, if $ 1 \leq q < \infty $, then the space $ C^\infty(\overline{U}) $ is dense in $ BV(U) \cap L^q(U) $ for the intermediate convergence (cf.\ \citep[Definition 10.1.3 and Theorem 10.1.2]{ABM}), i.e., for any $ z \in BV(U) \cap L^q(U) $ there exists a sequence $ \{ z_n \, | n \in \N \} \subset C^\infty(\overline{U}) $ such that $ z_n \to z $ in $ L^q(U) $ and strictly in $ BV(U) $ as $ n \to \infty $.
\end{em}
\end{notn}
\begin{notn}[weighted total variation; cf.\ \citep{AB,AFP}]\label{Note05}
\begin{em}
In this paper, we define
\begin{equation*}
X_{\rm c}(\Omega) := \{
\varpi \in L^\infty(\Omega)^N  | \ \mbox{$ {\rm div} \, \varpi \in L^2(\Omega) $ \ and \ $ {\rm supp} \, \varpi $ is compact in $ \Omega $}
\},
\end{equation*}
\begin{equation*}
W_0(\Omega) := \left\{ \begin{array}{l|l}
\varrho \in H^1(\Omega) \cap L^\infty(\Omega) & \varrho \geq 0 \mbox{ a.e.\ in $ \Omega $}
\end{array} \right\},
\end{equation*}
\begin{equation}\label{W_c}
W_{\rm c}(\Omega) := \left\{ \begin{array}{l|l}
\varrho \in H^1(\Omega) \cap L^\infty(\Omega) & \parbox{5.5cm}{there exists $ c_\varrho > 0 $ such that $ \varrho \geq c_\varrho $ a.e.\ in $ \Omega $}
\end{array} \right\},
\end{equation}
and for any $ \varrho \in W_0(\Omega) $ and any $ z \in L^2(\Omega) $, we call the value $ {\rm Var}_\varrho(z) \in [0, \infty] $, defined as,
\begin{equation*}
{\rm Var}_\varrho(v) := \sup \left\{ \begin{array}{l|l}
\displaystyle \int_\Omega v \, {\rm div} \, \varpi \, dx & \parbox{3.5cm}{$ \varpi \in X_{\rm c} $ and \\$ | \varpi | \leq \varrho $ a.e.\ in $ \Omega $}
\end{array} \right\} \in [0, \infty],
\end{equation*}
``the total variation of $ v $ weighted by $ \varrho $,'' or the ``weighted total variation'' for short. 
\end{em}
\end{notn}
\begin{rem}\label{Rem.Note05}
\begin{em}
Referring to the general theories (e.g., \citep{AB,AFP,BBF}), we can confirm the following facts associated with the weighted total variations: 
\begin{description}
\item[{(\hypertarget{Fact2}{Fact\,2})}](cf.\ \citep[Theorem 5]{BBF}) For any $ \varrho \in W_0(\Omega) $, the functional $ z \in L^2(\Omega) \mapsto {\rm Var}_\varrho(z) \in [0, \infty] $ is a proper l.s.c.\ and convex function that coincides with the lower semi-continuous envelope of
\begin{equation*}
z \in W^{1, 1}(\Omega) \cap L^2(\Omega) \mapsto \int_\Omega \varrho |\nabla z| \, dx \in [0, \infty).
\end{equation*}
\item[{(\hypertarget{Fact3}{Fact\,3})}](cf.\ \citep[Theorem 4.3]{AB} and \citep[Proposition 5.48]{AFP}) If $ \varrho \in W_0(\Omega) $ and $ z \in BV(\Omega) \cap L^2(\Omega) $, then there exists a Radon measure $ |Dz|_\varrho \in \mathcal{M}(\Omega) $ such that
\begin{equation*}
|Dz|_\varrho(\Omega) = \int_\Omega  d|Dz|_\varrho = {\rm Var}_\varrho(z)
\end{equation*}
and 
\begin{equation}\label{|Dz|_beta(A)}
\begin{array}{c}
\left\{ ~ \parbox{12.25cm}{
$ |Dz|_\varrho(A) \leq |\varrho|_{L^\infty(\Omega)} |Dz|(A) $,
\\[1ex]
$ \displaystyle
|Dz|_\varrho(A) = \inf \left\{ \begin{array}{l|l}
\displaystyle \liminf_{n \to \infty} \int_A \varrho |\nabla \tilde{z}_n| \, dx & \parbox{3.4cm}{$ \{ \tilde{z}_n \, | \, n \in \N \} \subset W^{1, 1}(A) \cap L^2(A) $ such that \ $ \tilde{z}_n \to z $ in $ L^2(A) $  as $ n \to \infty $}
\end{array} \right\} $
} \right.
\end{array}
\end{equation}
for any open set $ A \subset \Omega $.
\item[{(\hypertarget{Fact4}{Fact\,4})}]If $ \varrho \in W_{\rm c}(\Omega) $ and $ z \in BV(\Omega) \cap L^2(\Omega) $, then for any open set $ A \subset \Omega $, it follows that
\begin{equation}\label{|Dz|_BV}
\left\{ ~ {\parbox{10cm}{
$ |Dv|_\varrho(A) \geq c_\varrho |Dz|(A) $ \ for any open set $ A \subset \Omega $,
\\[2ex]
$ D({\rm Var}_\varrho) = BV(\Omega) \cap L^2(\Omega) $, \ and
\\[1ex]
$ {\rm Var}_\varrho(z) = \sup \left\{ \begin{array}{l|l}
\displaystyle \int_\Omega z \, {\rm div} \, (\varrho \varphi) \, dx & \parbox{3cm}{$ \varphi \in X_{\rm c}(\Omega) $ and $ | \varphi | \leq 1 $ a.e.\ in $ \Omega $}
\end{array} \right\}, $}} \right.
\end{equation}
where $ c_\varrho $ is a constant as in (\ref{W_c}). 
\end{description}
Moreover, the following properties can be inferred from (\ref{|Dz|_beta(A)})--(\ref{|Dz|_BV}):
\begin{itemize}
\item[$\bullet$] $ |Dz|_c = c|Dz| $ in $ \mathcal{M}(\Omega) $ for any constant $c \geq 0 $ and $ z \in BV(\Omega) \cap L^2(\Omega) $;
\item[$\bullet$] $ |Dz|_\varrho = \varrho |\nabla z| \mathscr{L}^N $ in $ \mathcal{M}(\Omega) $, \ if $ \varrho \in W_0(\Omega) $ and $ z \in W^{1, 1}(\Omega) \cap L^2(\Omega) $.
\end{itemize}
\end{em}
\end{rem}

\begin{notn}[generalized weighted total variation; cf.\ {\citep[Section 2]{MS}}]\label{Note06}
\begin{em}
 \ For any \linebreak $ \varrho \in H^1(\Omega) \cap L^\infty(\Omega) $ and any $ z \in BV(\Omega) \cap L^2(\Omega) $ we define a real-valued Radon measure $ [\varrho |Dz|] \in \mathcal{M}(\Omega) $, as follows:
\begin{equation*}
[\varrho |Dz|](B) := |Dz|_{[\varrho]^+}(B) - |Dz|_{[\varrho]^-}(B) \mbox{ \ for any Borel set $ B \subset \Omega $.}
\end{equation*}
Note that $ [\varrho|\nabla z|](\Omega) $ can be thought of as a generalized version of the total variation of $ z \in BV(\Omega) \cap L^2(\Omega) $ weighted by the possibly sign-changing weight $ \varrho \in H^1(\Omega) \cap L^\infty(\Omega) $. So, hereafter, we simply refer to $ [\varrho|Dz|](\Omega) $ as the \textit{generalized weighted total variation}.
\end{em}
\end{notn}

\begin{rem}\label{Rem.Note06}
\begin{em}
With regard to the generalized weighted total variations, the following facts are verified in \citep[Section 2]{MS}:
\begin{description}
\item[{(\hypertarget{Fact5}{Fact\,5})}](strict approximation)
Let $ \varrho \in H^1(\Omega) \cap L^\infty(\Omega) $ and $ z \in BV(\Omega) \cap L^2(\Omega) $ be arbitrary fixed functions, and let $ \{ z_n \, | \, n \in \N \} \subset C^\infty(\overline{\Omega}) $ be any sequence such that
\begin{equation*}
z_n \to z \mbox{ in $ L^2(\Omega) $ \ and \ strictly \ in } BV(\Omega) \mbox{ as $ n \to \infty $.}
\end{equation*}
Then
\begin{equation*}
\int_\Omega \varrho |\nabla z_n| \, dx \to \int_\Omega d[\varrho |Dz|] \mbox{ \ as $ n \to \infty $.}
\end{equation*}
\item[{(\hypertarget{Fact6}{Fact\,6})}]For any $ z \in BV(\Omega) \cap L^2(\Omega) $, the mapping
\begin{equation*}
\ds \varrho \in H^1(\Omega) \cap L^\infty(\Omega) \mapsto \int_\Omega d[\varrho |Dz|] \in \R
\end{equation*}
is a linear functional, and moreover, if $ \varphi \in H^1(\Omega) \cap C(\overline{\Omega}) $ and $ \varrho \in H^1(\Omega) \cap L^\infty(\Omega) $, then
\begin{equation*}
\displaystyle \int_\Omega d[\varphi \varrho |Dz|] = \int_\Omega \varphi \, d[\varrho |Dz|].
\end{equation*}
\end{description}
\end{em}
\end{rem}

\begin{notn}[specific classes of functions]\label{Note07}
\begin{em}
Let $ X_0 $ be a Banach space, defined as
\begin{equation*}
X_0 := H^1(\Omega) \times H^1(\Omega) \times BV(\Omega).
\end{equation*}
For any $ 1 \leq p \leq \infty $ and any open interval $ I \subset \R $, we set
\begin{equation*}
L^p(I; BV(\Omega)) := \left\{ \begin{array}{l|l}
\tilde{\theta} & \parbox{8cm}{
$ \tilde{\theta} : I \rightarrow BV(\Omega) $ is measurable for the strict topology of $ BV(\Omega) $, and $ |\tilde{\theta}({}\cdot{})|_{BV(\Omega)} \in L^p(I) $
} 
\end{array} \right\}
\end{equation*}
and
\begin{equation*}
L^p(I; X_0) := \left\{ \begin{array}{l|l}
[\tilde{w}, \tilde{\eta}, \tilde{\theta}] & \parbox{7.7cm}{
$ \tilde{w}, \tilde{\eta} \in L^p(I; H^1(\Omega)) $ and $ \tilde{\theta} \in L^p(I; BV(\Omega)) $
} 
\end{array} \right\}.
\end{equation*}
For any $ 1 \leq p \leq \infty $ and any open interval $ I \subset \R $, we note that $ L^p(I; BV(\Omega)) $ and $ L^p(I; X_0) $ are normed spaces, with
\begin{equation*}
|\tilde{\theta}|_{L^p(I; BV(\Omega))} := \bigl| |\tilde{\theta}({}\cdot{})|_{BV(\Omega)} \bigr|_{L^p(I)}, \mbox{ for $ \tilde{\theta} \in L^p(I; BV(\Omega)) $}
\end{equation*}
and 
\begin{equation*}
\begin{array}{c}
|[\tilde{w}, \tilde{\eta}, \tilde{\theta}]|_{L^p(I; X_0)} := |\tilde{w}|_{L^p(I; H^1(\Omega))} +|\tilde{\eta}|_{L^p(I; H^1(\Omega))} +|\tilde{\theta}|_{L^p(I; BV(\Omega))} 
\\[1ex]
\mbox{for $ [\tilde{w}, \tilde{\eta}, \tilde{\theta}] \in L^p(I; X_0) $,}
\end{array}
\end{equation*}
respectively. 
\end{em}
\end{notn}

Finally, we mention the notion of functional convergences.
 
\begin{defn}[Mosco convergence; cf.\ \citep{Mosco}]\label{Def.Mosco}
\begin{em}
Let $ X $ be an abstract Hilbert space. Let $ \Psi : X \rightarrow (-\infty, \infty] $ be a proper l.s.c.\ and convex function, and let $ \{ \Psi_n \, | \, n \in \N \} $ be a sequence of proper l.s.c.\ and convex functions $ \Psi_n : X \rightarrow (-\infty, \infty] $, $ n \in \N $. We say that $ \Psi_n \to \Psi $ on $ X $, in the sense of Mosco \citep{Mosco}, as $ n \to \infty $, if and only if the following two conditions are fulfilled:
\begin{description}
\item[\textmd{(\hypertarget{m1}{m1})}](lower bound) $ \ds \liminf_{n \to \infty} \Psi_n(z_n^\dag) \geq \Psi(z^\dag) $ if $ z^\dag \in X $, $ \{ z_n^\dag \, | \, n \in \N \} \subset X $, and $ z_n^\dag \to z^\dag $ weakly in $ X $ as $ n \to \infty $;
\item[\textmd{(\hypertarget{m2}{m2})}](optimality) for any $ z^\ddag \in D(\Psi) $, there exists a sequence $ \{ z_n^\ddag  \, | \, n \in \N \} \subset X $ such that $ z_n^\ddag \to z^\ddag $ in $ X $ and $ \Psi_n(z_n^\ddag) \to \Psi(z^\ddag) $ as $ n \to \infty $.
\end{description} 
\end{em}
\end{defn}

\begin{defn}[$\Gamma$-convergence; cf.\ \citep{GammaConv}]\label{Def.Gamma}
\begin{em}
Let $ X $ be an abstract Hilbert space, $ \Psi : X \rightarrow (-\infty, \infty] $ be a proper functional, and $ \{ \Psi_n \, | \, n \in \N \} $ be a sequence of proper functionals $ \Psi_n : X \rightarrow (-\infty, \infty] $, $ n \in \N $. We say that $ \Psi_n \to \Psi $ on $ X $, in the sense of $ \Gamma $-convergence \citep{GammaConv}, as $ n \to \infty $ if and only if the following two conditions are fulfilled:
\begin{description}
\item[\textmd{(\hypertarget{gamma1}{$\gamma$1})}](lower bound) $ \ds \liminf_{n \to \infty} \Psi_n(z_n^\dag) \geq \Psi(z^\dag) $ if $ z^\dag \in X $, $ \{ z_n^\dag \, | \, n \in \N \} \subset X $, and $ z_n^\dag \to z^\dag $ (strongly) in $ X $ as $ n \to \infty $;
\item[\textmd{(\hypertarget{gamma2}{$\gamma$2})}](optimality) for any $ z^\ddag \in D(\Psi) $, there exists a sequence $ \{ z_n^\ddag  \, | \, n \in \N \} \subset X $ such that $ z_n^\ddag \to z^\ddag $ in $ X $ and $ \Psi_n(z_n^\ddag) \to \Psi(z^\ddag) $ as $ n \to \infty $.
\end{description} 
\end{em}
\end{defn}

\begin{rem}\label{Rem.MG}
\begin{em}
Note that if the functionals are convex, then Mosco convergence implies $ \Gamma $-convergence, i.e., the $ \Gamma $-convergence of convex functions can be regarded as a weak version of Mosco convergence. Additionally, in the $ \Gamma $-convergence of convex functions, we can see the following:
\vspace{1ex}
\begin{description}
\item[(\hypertarget{Fact7}{Fact\,7})]Let $ \Psi : X \rightarrow (-\infty, \infty] $ and $ \Psi_n : X \rightarrow (-\infty, \infty] $ be proper l.s.c.\ and convex functions on a Hilbert space $ X $ such that $ \Psi_n \to \Psi $ on $ X $, in the sense of $ \Gamma $-convergence, as $ n \to \infty $. Assume that
\begin{equation*}
\left\{ ~ \parbox{9cm}{
$ [z, z^*] \in X^2 $, ~ $ [z_n, z_n^*] \in \partial \Psi_n $ in $ X^2 $, $ n \in \N $,
\\[1ex]
$ z_n \to z $ in $ X $ and $ z_n^* \to z^* $ weakly in $ X $, as $ n \to \infty $.
} \right.
\end{equation*}
It then holds that $ [z, z^*] \in \partial \Psi \mbox{ in $ X^2 $ and } \Psi_n(z_n) \to \Psi(z) $ as $ n \to \infty $.
\end{description}
\end{em}
\end{rem}

\section{Statement of the main theorems}
\ \vspace{-3ex}

We begin by reiterating the assumptions in this paper. 
\begin{enumerate}
\item[(\hypertarget{A1}{A1})]$ \gamma : \R \rightarrow [0, \infty] $ is a proper l.s.c.\ and convex function such that
\begin{equation*}
\left\{ \hspace{-1ex} \parbox{11cm}{
\begin{itemize}
\item $ K_\gamma := D(\gamma) \subset \R $ is a closed interval, and $ K_\gamma \supset (0, 1) $; 
\item $ \gamma \in C^1({\rm int}(K_\gamma)) $, so that the subdifferential $ \partial \gamma $ coincides with the usual differential $ \gamma' $ on $ {\rm int}(K_\gamma) $.
\end{itemize}
} \right.
\end{equation*}
\item[(\hypertarget{A2}{A2})]$ g(\,\cdot\,,\cdot\,) \in C^2(\R^2)  $, and there exists a constant $ c_* \in \R $ such that
\begin{equation*}
\gamma(\tilde{w}) +g(\tilde{w}, \tilde{\eta}) \geq c_* \mbox{ \ for all $ \tilde{w}, \tilde{\eta} \in [0, 1]^2 $.}
\end{equation*}
\item[(\hypertarget{A3}{A3})]$ \alpha_0 \in W_{\rm loc}^{1, \infty}(\R^2) $ is a positive-valued function, and $ \alpha, \beta \in C^1(\R^2) \cap C^2([0, 1]^2) $ are nonnegative-valued convex functions.
\item[(\hypertarget{A4}{A4})]There exist two constants $ o_*, \iota_* \in \R $ such that
\begin{equation*}
\{ o_*, \iota_* \} \subset D(\partial \gamma) \mbox{ and } 0 \leq o_* < \iota_* \leq 1, 
\end{equation*}
and the subdifferential $ \partial \gamma $ and the partial differentials
$ g_w(\,\cdot\,,\cdot\,) = \frac{\partial}{\partial w} g(\,\cdot\,,\cdot\,) $, $ g_\eta(\,\cdot\,,\cdot\,) = \frac{\partial}{\partial \eta} g(\,\cdot\,,\cdot\,) $ fulfill that
\begin{equation}\label{a4-1}
\left\{ ~ 
\parbox{12cm}{
$ \ds \bigcap_{\tilde{\eta} \in [0, 1]} \bigl( \partial \gamma(o_*)  +g_w(o_*, \tilde{\eta}) \bigr) \cap (-\infty, 0] \ne \emptyset $ \ and \ $ \ds \min_{\tilde{w} \in [0, 1]} g_\eta(\tilde{w}, 0) \leq 0 $,
\\[1ex]
$ \ds \bigcap_{\tilde{w} \in [0, 1]} \bigl( \partial \gamma(\iota_*) +g_w(\iota_*, \tilde{w}) \bigr) \cap [0, \infty) \ne \emptyset $ \ and \ $ \ds \max_{\tilde{w} \in [0, 1]} g_\eta(\tilde{w}, 1) \geq 0 $.
} \right. 
\end{equation}
Furthermore, the partial differentials $ \alpha_w= \frac{\partial \alpha}{\partial w} $, $ \alpha_\eta = \frac{\partial \alpha}{\partial \eta} $, $ \beta_w= \frac{\partial \beta}{\partial w} $ and $ \beta_\eta = \frac{\partial \beta}{\partial \eta} $ fulfill that
\begin{equation}\label{a4-2}
\left\{ ~ \parbox{10cm}{
$ \ds \sup_{[\tilde{w}, \tilde{\eta}] \in \R^2} \bigl\{ \alpha_w(o_*, \tilde{\eta}), \alpha_\eta(\tilde{w}, 0), \beta_w(o_*, \tilde{\eta}), \beta_\eta(\tilde{w}, 0) \bigr\} \leq 0 $,
\\[1ex]
$ \ds \inf_{[\tilde{w}, \tilde{\eta}] \in \R^2} \bigl\{ \alpha_w(\iota_*, \tilde{\eta}), \alpha_\eta(\tilde{w}, 1), \beta_w(\iota_*, \tilde{\eta}), \beta_\eta(\tilde{w}, 1) \bigr\} \geq 0 $.
} \right.
\end{equation}
\end{enumerate}
\begin{rem}\label{Rem.Assumption}
\begin{em}
Assumptions (\hyperlink{A1}{A1})--(\hyperlink{A4}{A4}) cover all of the settings presented in (\hyperlink{g1}{g1})--(\hyperlink{g3}{g3}) and (\ref{possSet}). In particular, we note that assumptions (\hyperlink{A3}{A3})--(\hyperlink{A4}{A4}) encompass more interactional functions, such as
\begin{equation*}
[\tilde{w}, \tilde{\eta}] \in \R^2 \mapsto \bigl( |\tilde{w} -o_*|^p +|\tilde{\eta}|^q \bigr)^r \mbox{ with constants $ p, q, r > 1 $}
\end{equation*}
as possible expressions of the mobilities $ \alpha $ and $ \beta $. 
\end{em}
\end{rem}

Next, for descriptive convenience, we introduce the following abbreviations.

\begin{notn}\label{Note08}
\begin{em}
For any $ \tilde{v} = [\tilde{w}, \tilde{\eta}] \in \R^2 $, we abbreviate $ \alpha_0(\tilde{w}, \tilde{\eta}) $, $ \alpha(\tilde{w}, \tilde{\eta}) $, and $ \beta(\tilde{w}, \tilde{\eta}) $ by $ \alpha_0(\tilde{v}) $, $ \alpha(\tilde{v}) $, and $ \beta(\tilde{v}) $, respectively, and we set
\begin{equation*}
\left\{ ~ \parbox{9.8cm}{
$ [\nabla g](\tilde{v}) = [\nabla g](\tilde{w}, \tilde{\eta}) := [g_w(\tilde{w}, \tilde{\eta}), g_\eta(\tilde{w}, \tilde{\eta})] $
\\[1ex]
$ [\nabla \alpha](\tilde{v}) = [\nabla \alpha](\tilde{w}, \tilde{\eta}) := [\alpha_w(\tilde{w}, \tilde{\eta}), \alpha_\eta(\tilde{w}, \tilde{\eta})] $
\\[1ex]
$ [\nabla \beta](\tilde{v}) = [\nabla \beta](\tilde{w}, \tilde{\eta}) := [\beta_w(\tilde{w}, \tilde{\eta}), \beta_\eta(\tilde{w}, \tilde{\eta})] $
} \right. \mbox{for $ \tilde{v} = [\tilde{w}, \tilde{\eta}] \in \R^2 $.}
\end{equation*}
For any $ \tilde{v} = [\tilde{w}, \tilde{\eta}] \in [H^1(\Omega) \cap L^\infty(\Omega)]^2 $, let $ \Phi_\nu(\tilde{v};{}\cdot{}) = \Phi_\nu(\tilde{w}, \tilde{\eta}; {}\cdot{}) $ be a proper l.s.c.\ and convex function on $ L^2(\Omega) $ defined as
\begin{eqnarray*}
\vartheta \in L^2(\Omega) & \mapsto & \Phi_\nu(\tilde{v};\vartheta) = \Phi_\nu(\tilde{w}, \tilde{\eta}; \vartheta)
\\[1ex]
& := & \left\{ \begin{array}{ll}
\multicolumn{2}{l}{\ds \int_\Omega \alpha(\tilde{v}) |\nabla \vartheta| \, dx +{\nu} \int_\Omega \beta(\tilde{v}) |\nabla \vartheta|^2 \, dx,}
\\[2ex]
& \mbox{if $ \nu > 0 $ and $ \vartheta \in H^1(\Omega) $,}
\\[1ex]
\multicolumn{2}{l}{\ds \int_\Omega d[\alpha(\tilde{v}) |D \vartheta|], \mbox{ \ if $ \nu = 0 $ and $ \vartheta \in BV(\Omega) $,}}
\\[2ex]
\infty, & \mbox{otherwise,}
\end{array}
\right.
\end{eqnarray*}
and let $ \partial \Phi_\nu(\tilde{v};{}\cdot{}) = \partial \Phi_\nu(\tilde{w}, \tilde{\eta};{}\cdot{}) $ be the $ L^2 $-subdifferential of $ \Phi_\nu(\tilde{v};{}\cdot{}) = \Phi_\nu(\tilde{w}, \tilde{\eta}; {}\cdot{}) $.
\end{em}
\end{notn}

\begin{rem}\label{Rem.FE}
\begin{em}
By virtue of Notation \ref{Note08}, we can uniformly provide proper definitions of the free energies in (\ref{freeEnergy}) and (\ref{freeEnergy0}) by assigning
\begin{equation}\label{FE00}
\begin{array}{rcl}
[\tilde{v}, \tilde{\theta}] = [\tilde{w}, \tilde{\eta}, \tilde{\theta}] & \in & L^2(\Omega)^3 \mapsto \mathscr{F}_\nu(\tilde{v}, \tilde{\theta}) = \mathscr{F}_\nu(\tilde{w}, \tilde{\eta}, \tilde{\theta}) 
\\[1ex]
& := & V_{\rm D}^2(\tilde{w}, \tilde{\eta}) +\Gamma(\tilde{v}) +\mathscr{G}(\tilde{v}) +\Phi_\nu(\tilde{v}; \tilde{\theta}) \mbox{ for all $ \nu \geq 0 $.}
\end{array}
\end{equation}
In this context,
\begin{itemize}
\item[--]$ V_{\rm D}^2 $ is the functional given in (\ref{Dirichlet^d}) for the case $ d = 2 $;
\item[--]$ \Gamma $ is a proper l.s.c.\ and convex function on $ L^2(\Omega)^2 $, defined as
\begin{equation}\label{Gamma(v)}
\tilde{v} = [\tilde{w}, \tilde{\eta}] \in L^2(\Omega)^2 \mapsto \Gamma(\tilde{v}) = \Gamma(\tilde{w}, \tilde{\eta}) := \int_\Omega \gamma(\tilde{w}) \, dx \in (-\infty, \infty];
\end{equation}
\item[--]$ \mathscr{G}(\,\cdot\,) $ is a functional on $ L^2(\Omega)^2 $ given by
\begin{equation*}
\tilde{v} = [\tilde{w}, \tilde{\eta}] \in L^2(\Omega)^2 \mapsto \mathscr{G}(\tilde{v}) = \mathscr{G}(\tilde{w}, \tilde{\eta}) := \int_\Omega g(\tilde{w}, \tilde{\eta}) \, dx \in \R.
\end{equation*} 
\end{itemize}
\end{em}
\end{rem}
\medskip

The main theorems can now be stated. 

\begin{mainTh}[{\boldmath existence for $\mbox{(\hyperlink{(S;u)_nu}{S})}_\nu$ when $ \nu > 0 $}]\label{MainTh01}
Fix the constant $ \nu > 0 $ and assume (\hyperlink{A1}{A1})--(\hyperlink{A4}{A4}). Additionally, assume that
\begin{equation}\label{delta_1}
\delta_1 := \inf \left\{ \begin{array}{l|l}
\alpha_0(\tilde{w}, \tilde{\eta}), \beta(\tilde{w}, \tilde{\eta}) & [\tilde{w}, \tilde{\eta}] \in \R^2 
\end{array} \right\} > 0
\end{equation}
and
\begin{equation}\label{D_1}
[w_0, \eta_0, \theta_0] \in D_1 := \left\{ \begin{array}{l|l}
[\tilde{w}, \tilde{\eta}, \tilde{\theta}] \in H^1(\Omega)^3 & \parbox{4cm}{
$ o_* \leq \tilde{w} \leq \iota_* $ a.e.\ in $ \Omega $, $ 0 \leq \tilde{\eta} \leq 1 $ a.e.\ in $ \Omega $, and $ \tilde{\theta} \in L^\infty(\Omega) $
}
\end{array} \right\}.
\end{equation}
The system $\mbox{(\hyperlink{(S;u)_nu}{S})}_\nu$ then admits at least one solution $ [w, \eta, \theta] $, defined by the following conditions.
\vspace{1ex}
\begin{description}
\item[\textmd{\it$\mbox{(\hypertarget{(S0)_nu}{S0})}_\nu$}]$ [w, \eta, \theta] \in W^{1,2}(0, T; L^2(\Omega)^3) \cap L^\infty(0, T; H^1(\Omega)^3) \cap L^\infty(Q)^3 $, $ o_* \leq w \leq \iota_* $, $ 0 \leq \eta \leq 1 $, and $ |\theta| \leq |\theta_0|_{L^\infty(\Omega)} $, a.e.\ in $ Q $.
\vspace{1ex}
\item[\textmd{\it$\mbox{(\hypertarget{(S1)_nu}{S1})}_\nu$}]$ w $ solves (\ref{1st.eq}) in the following variational sense:
\vspace{-0.5ex}
\begin{equation} \label{1st.VarIneq}
\begin{array}{l}
\ds \hspace{-4ex} \int_\Omega \left( \rule{0pt}{10pt} w_t(t) +g_w(w, \eta)(t) \right)(w(t) -\varphi) \, dx +\int_\Omega \nabla w(t) \cdot \nabla (w(t) -\varphi) \, dx 
\\[2ex]
\ds +\int_\Omega \left( \rule{0pt}{10pt} \alpha_w(w, \eta)(t)|\nabla \theta(t)| +\nu \beta_w(w, \eta)(t) |\nabla \theta(t)|^2 \right) (w(t) -\varphi) \, dx 
\\[2ex]
\ds +\int_\Omega \gamma(w(t)) \, dx \leq \int_\Omega \gamma(\varphi) \, dx 
\\[3ex]
\mbox{for any $ \varphi \in H^1(\Omega) \cap L^\infty(\Omega) $ and a.e.\ $ t \in (0, T) $,}
\end{array}
\end{equation}
with the initial condition $ w(0) = w_0 $ in $ L^2(\Omega) $.
\vspace{1ex}
\item[\textmd{\it$\mbox{(\hypertarget{(S2)_nu}{S2})}_\nu$}]$ \eta $ solves (\ref{2nd.eq}) in the following variational sense:
\vspace{-0.5ex}
\begin{equation} \label{2nd.VarIneq}
\begin{array}{l}
\ds \hspace{-4ex} \int_\Omega \left( \rule{0pt}{10pt} \eta_t(t) +g_\eta(w, \eta)(t) \right) \psi \, dx +\int_\Omega \nabla \eta(t) \cdot \nabla \psi \, dx 
\\[2ex]
\ds +\int_\Omega \left( \rule{0pt}{10pt} \alpha_\eta(w, \eta)(t) |\nabla \theta(t)| +\nu \beta_\eta(w, \eta)(t)|\nabla \theta(t)|^2 \right) \psi \, dx = 0 
\\[3ex]
\mbox{for any $ \psi \in H^1(\Omega) \cap L^\infty(\Omega) $ and a.e.\ $ t \in (0, T) $,}
\end{array}
\end{equation}
with the initial condition $ \eta(0) = \eta_0 $ in $ L^2(\Omega) $.
\vspace{1ex}
\item[\textmd{\it$\mbox{(\hypertarget{(S3)_nu}{S3})}_\nu$}]$ \theta $ solves (\ref{3rd.eq}) in the following variational sense:
\vspace{-0.5ex}
\begin{equation} \label{3rd.VarIneq}
\begin{array}{l}
\ds \hspace{-4ex} \int_\Omega \alpha_0(w, \eta)(t) \, \theta_t(t) \, (\theta(t) -\omega) \, dx +2 \nu \int_\Omega \beta(w, \eta)(t) \nabla \theta(t) \cdot \nabla (\theta(t) -\omega) \, dx 
\\[2ex]
\ds +\int_\Omega \alpha(w, \eta)(t) |\nabla \theta(t)| \, dx \leq \int_\Omega \alpha(w, \eta)(t) |\nabla \omega| \, dx
\\[3ex]
\mbox{for any $ \omega \in H^1(\Omega) $ and a.e.\ $ t \in (0, T) $,}
\end{array}
\end{equation}
with the initial condition $ \theta(0) = \theta_0 $ in $ L^2(\Omega) $.
\end{description}
\end{mainTh}
\medskip

\begin{mainTh}[{\boldmath existence for $\mbox{(\hyperlink{(S;u)_0}{S})}_0$}]\label{MainTh02}
Assume (\hyperlink{A1}{A1})--(\hyperlink{A4}{A4}), and in addition, assume that
\begin{equation}\label{delta_0}
\delta_0 := \inf \left\{ \begin{array}{l|l}
\alpha_0(\tilde{w}, \tilde{\eta}), \alpha(\tilde{w}, \tilde{\eta}) & [\tilde{w}, \tilde{\eta}] \in \R^2 
\end{array} \right\} > 0
\end{equation}
and
\begin{equation}\label{D_0}
[w_0, \eta_0, \theta_0] \in D_0 := \left\{ \begin{array}{l|l}
[\tilde{w}, \tilde{\eta}, \tilde{\theta}] \in X_0 & \parbox{4cm}{
$ o_* \leq \tilde{w} \leq \iota_* $ a.e.\ in $ \Omega $, $ 0 \leq \tilde{\eta} \leq 1 $ a.e.\ in $ \Omega $, and $ \tilde{\theta} \in L^\infty(\Omega) $
}
\end{array} \right\}.
\end{equation}
The system $\mbox{(\hyperlink{(S;u)_0}{S})}_0$ then admits at least one solution $ [w, \eta, \theta] $, defined by the following conditions.
\vspace{1ex}
\begin{description}
\item[\textmd{\it $\mbox{(\hypertarget{(S0)_0}{S0})}_0$}]$ [w, \eta, \theta] \in W^{1,2}(0, T; L^2(\Omega)^3) \cap L^\infty(0, T; X_0) \cap L^\infty(Q)^3 $, 
$ o_* \leq w \leq \iota_* $, $ 0 \leq \eta \leq 1 $, and $ |\theta| \leq |\theta_0|_{L^\infty(\Omega)} $, a.e.\ in $ Q $.
\vspace{-0.5ex}
\item[\textmd{\it $\mbox{(\hypertarget{(S1)_0}{S1})}_0$}]$ w $ solves (\ref{1st.eq0}) in the following variational sense:
\vspace{-0.5ex}
\begin{equation} \label{1st.VarIneq0}
\begin{array}{l}
\ds \hspace{-4ex} \int_\Omega \left( \rule{0pt}{10pt} w_t(t) +g_w(w, \eta)(t) \right)(w(t) -\varphi) \, dx  +\int_\Omega \nabla w(t) \cdot \nabla (w(t) -\varphi) \, dx
\\[2ex]
\ds+\int_\Omega d \left[ \rule{0pt}{10pt} (w(t) -\varphi) \alpha_w(w, \eta)(t)|D \theta(t)| \right]  +\int_\Omega \gamma(w(t)) \, dx \leq \int_\Omega \gamma(\varphi) \, dx 
\\[3ex]
\mbox{for any $ \varphi \in H^1(\Omega) \cap L^\infty(\Omega) $ and a.e.\ $ t \in (0, T) $,}
\end{array}
\end{equation}
with the initial condition $ w(0) = w_0 $ in $ L^2(\Omega) $.
\vspace{1ex}
\item[\textmd{\it $\mbox{(\hypertarget{(S2)_0}{S2})}_0$}]$ \eta $ solves (\ref{2nd.eq0}) in the following variational sense:
\vspace{-0.5ex}
\begin{equation} \label{2nd.VarIneq0}
\begin{array}{l}
\ds \hspace{-10ex} \int_\Omega \left( \rule{0pt}{10pt} \eta_t(t) +g_\eta(w, \eta)(t) \right) \psi \, dx +\int_\Omega \nabla \eta(t) \cdot \nabla \psi \, dx
\\[2ex]
\ds 
+\int_\Omega d \left[ \rule{0pt}{10pt} \psi \alpha_\eta(w, \eta)(t) |D \theta(t)| \right] = 0 
\\[3ex]
\mbox{for any $ \psi \in H^1(\Omega) \cap L^\infty(\Omega) $ and a.e.\ $ t \in (0, T) $,}
\end{array}
\end{equation}
with the initial condition $ \eta(0) = \eta_0 $ in $ L^2(\Omega) $.
\vspace{1ex}
\item[\textmd{\it $\mbox{(\hypertarget{(S3)_0}{S3})}_0$}]$ \theta $ solves (\ref{3rd.eq0}) in the following variational sense:
\vspace{-0.5ex}
\begin{equation} \label{3rd.VarIneq0}
\begin{array}{l}
\ds \hspace{-7ex} \int_\Omega \alpha_0(w, \eta)(t) \, \theta_t(t) \, (\theta(t) -\omega) \, dx 
\\[2ex]
\ds +\int_\Omega d[\alpha(w, \eta)(t) |D \theta(t)|] \leq \int_\Omega d [\alpha(w, \eta)(t) |D \omega|]
\\[3ex]
\hspace{-7ex}\mbox{for any $ \omega \in BV(\Omega) \cap L^2(\Omega) $ and a.e.\ $ t \in (0, T) $,}
\end{array}
\end{equation}
with the initial condition $ \theta(0) = \theta_0 $ in $ L^2(\Omega) $.
\end{description}
\end{mainTh}

\begin{rem}\label{Rem.1.0}
\begin{em}
Hereafter, whenever $ \nu \geq 0 $, we set
\begin{equation*}
v := [w, \eta] \mbox{ in $ C([0, T]; L^2(\Omega)^2) $, \ for the solution $ [w, \eta, \theta] $ to $ \mbox{(\hyperlink{(S;u)_nu}{S})}_\nu $.}
\end{equation*}
Thus, if $ \nu > 0 $, then the variational inequalities (\ref{1st.VarIneq})--(\ref{2nd.VarIneq}) for $ v = [w, \eta] $ can be unified as follows:
\begin{equation}\label{unify01}
\begin{array}{rl}
\multicolumn{2}{l}{\ds \hspace{-2ex} \left( v_t(t), v(t) -\varpi \right)_{L^2(\Omega)^2} +\left( [\nabla g](v(t)), v(t) -\varpi \right)_{L^2(\Omega)^2}}
\\[2ex]
& \ds + \left( \nabla w(t), \nabla (w(t) -\varphi) \right)_{L^2(\Omega)^N} +\left( \nabla \eta(t), \nabla (\eta(t) -\psi) \right)_{L^2(\Omega)^N}
\\[1ex]
& \ds +\int_\Omega \left( |\nabla \theta(t)| [\nabla \alpha](v(t)) +\nu |\nabla \theta(t)|^2 [\nabla \beta](v(t)) \right) \cdot (v(t) -\varpi) \, dx
\\[2ex]
& + \Gamma(v(t)) \leq \Gamma(\varpi) \mbox{ \ for any $ \varpi = [\varphi, \psi] \in [H^1(\Omega) \cap L^\infty(\Omega)]^2 $.}
\end{array}
\end{equation}
Meanwhile, if $ \nu = 0 $, then the corresponding variational inequalities (\ref{1st.VarIneq0})--(\ref{2nd.VarIneq0}) can be unified as
\begin{equation}\label{unify02}
\begin{array}{rl}
\multicolumn{2}{l}{\ds \hspace{-2ex} \left( v_t(t), v(t) -\varpi \right)_{L^2(\Omega)^2} +\left( [\nabla g](v(t)), v(t) -\varpi \right)_{L^2(\Omega)^2}}
\\[2ex]
& \ds + \left( \nabla w(t), \nabla (w(t) -\varphi) \right)_{L^2(\Omega)^N} +\left( \nabla \eta(t), \nabla (\eta(t) -\psi) \right)_{L^2(\Omega)^N}
\\[1ex]
& \ds +\int_\Omega  d \left[ \left( \rule{0pt}{10pt} (v(t) -\varpi) \cdot [\nabla \alpha](v(t)) \right) |D \theta(t)| \right]
\\[2ex]
& + \Gamma(v(t)) \leq \Gamma(\varpi) \mbox{ \ for any $ \varpi = [\varphi, \psi] \in [H^1(\Omega) \cap L^\infty(\Omega)]^2 $.}
\end{array}
\end{equation}
Moreover, for every $ \nu \geq 0 $ the variational inequalities for $ \theta $, i.e., (\ref{3rd.VarIneq}) when $ \nu > 0 $ and (\ref{3rd.VarIneq0}) when $ \nu = 0 $, uniformly reduce to the following form of an evolution equation:
\begin{equation}\label{simple02}
\alpha_0(v(t)) \, \theta_t(t) +\partial \Phi_\nu(v(t); \theta(t)) \ni 0 \mbox{ \ in $ L^2(\Omega) $, \ a.e.\ $ t \in (0, T) $,}
\end{equation}
governed by the subdifferentials $ \partial \Phi_\nu(v(t);{}\cdot{}) $ of unknown-dependent convex functions $ \Phi_\nu(v(t);{}\cdot{}) $ for $ t \in [0, T] $. 
\end{em}
\end{rem}

\begin{rem}\label{Rem.1.1}
\begin{em}
In fact, reductions similar to (\ref{simple02}) are available for the variational inequalities for $ v = [w, \eta] $ in some special cases. For instance, if $ K_\gamma = (-\infty, \infty) $, then by (\hyperlink{A1}{A1}) the functional $ \Gamma $ becomes locally Lipschitz on $ L^2(\Omega)^2 $, and the subdifferential $ \partial \Gamma \subset [L^2(\Omega)^2]^2 $ coincides with the Fr\'{e}chet differential $ \Gamma' $. Therefore, with (\ref{F}), (\ref{elliptic01}), (\hyperlink{Fact6}{Fact\,6}) in Remark \ref{Rem.Note06} and \citep[Chapter 2]{Brezis} in mind, the variational inequalities (\ref{1st.VarIneq})--(\ref{2nd.VarIneq}) and (\ref{1st.VarIneq0})--(\ref{2nd.VarIneq0}) can be reduced to the following evolution equation form:
\begin{equation*}
\begin{array}{ll}
\displaystyle 
v_t(t) +\left( \rule{0pt}{10pt} F v(t) -v(t) \right) & \hspace{-1.5ex} + \Gamma'(v(t)) +[\nabla g](v(t)) +[|D \theta(t)|[\nabla\alpha](v(t))]
\\[1ex]
\multicolumn{2}{c}{+|\nabla (\nu \theta)(t)|[\nabla \beta](v(t)) = 0 \mbox{ \ in \ } [H^s(\Omega)^*]^2 , \mbox{ \ a.e.\ $ t \in (0, T) $,}}
\end{array}
\end{equation*}
where $ s > N/2 $ is a large exponent to realize the embedding $ H^s(\Omega) \subset C(\overline{\Omega}) $, and for any $ \sigma = [\xi, \zeta] \in [H^1(\Omega) \cap L^\infty(\Omega)]^2 $ and any $ \vartheta \in BV(\Omega) \cap L^2(\Omega) $, $ [|D \vartheta| \sigma] $ denotes a (vectorial) Radon measure, defined as
\begin{equation*}
\ds \varpi := [\varphi, \psi] \in C_{0}(\Omega)^2 \mapsto \int_\Omega \varpi \, d [|D \vartheta| \sigma] := \int_\Omega \varphi \, d [\xi |D \vartheta|] +\int_\Omega \psi \, d [\zeta |D \vartheta|] \in \R.
\end{equation*}
However, it must be noted that such reductions may not be valid for general instances of the convex function $ \gamma $. 
\end{em}
\end{rem}

\section{Approximating problems and auxiliary lemmas}
\ \vspace{-3ex}

In this section, the approximating problems for our systems are presented along with relaxed versions. As mentioned in the introduction, these problems are formulated as time discretization systems for $ \mbox{(\hyperlink{(S;u)_nu}{S})}_\nu $ when $ \nu > 0 $. Hence, throughout the description of the approximating problems, we fix $ \nu $ as a positive constant and assume (\ref{delta_1}) as in Main Theorem \ref{MainTh01}. Additionally, for any time step $ 0 < h < 1 $ we denote by $ \mbox{(\hyperlink{(AP)_h^nu}{AP})}_h^\nu $ the approximating problem for $ \mbox{(\hyperlink{(S;u)_nu}{S})}_\nu $, prescribed as follows: 

\begin{description}
\item[\textmd{$ \mbox{(\hypertarget{(AP)_h^nu}{AP})}_h^\nu $}]For any initial value
\begin{equation}\label{apxInit}
[v_0^\nu, \theta_0^\nu] = [w_0^\nu, \eta_0^\nu, \theta_0^\nu] \in D_1 \mbox{ with $ v_0^\nu = [w_0^\nu, \eta_0^\nu] \in H^1(\Omega)^2 $,}
\end{equation}
find a sequence of triplets
\begin{equation*}
\begin{array}{l}
\{ [v_i^\nu, \theta_i^\nu] = [w_i^\nu, \eta_i^\nu, \theta_i^\nu] \, | \, i \in \N \} \subset D_1
\\[1ex]
\hspace{18ex}\mbox{with $ v_i^\nu = [w_i^\nu, \eta_i^\nu] \in H^1(\Omega)^2 $,  $ i \in \N $,}
\end{array}
\end{equation*}
such that 
\begin{equation}\label{apx01-02}
\begin{array}{rl}
\multicolumn{2}{l}{\ds \hspace{-2ex} \frac{1}{h} \left( v_i^\nu -v_{i -1}^\nu, v_i^\nu -\varpi \right)_{L^2(\Omega)^2} +\left( [\nabla g](v_i^\nu), v_i^\nu -\varpi \right)_{L^2(\Omega)^2}}
\\[2ex]
& \ds + \bigl( \nabla w_i^\nu, \nabla (w_i^\nu -\varphi) \bigr)_{L^2(\Omega)^N} +\bigl( \nabla \eta_i^\nu, \nabla (\eta_i^\nu -\psi) \bigr)_{L^2(\Omega)^N}
\\[1ex]
& \ds +\int_\Omega \left( |\nabla \theta_{i -1}^\nu| [\nabla \alpha](v_i^\nu) +\nu |\nabla \theta_{i -1}^\nu|^2 [\nabla \beta](v_i^\nu) \right) \cdot (v_i^\nu -\varpi) \, dx
\\[2ex]
& + \Gamma(v_i^\nu) \leq \Gamma(\varpi), \mbox{ \ for any $ \varpi = [\varphi, \psi] \in [H^1(\Omega) \cap L^\infty(\Omega)]^2 \cap D(\Gamma) $,}
\end{array}
\end{equation}
\begin{equation}\label{apx03}
\begin{array}{c}
\ds \frac{1}{h} \left( \alpha_0(v_i^\nu) (\theta_i^\nu -\theta_{i -1}^\nu), \theta_i^\nu -\omega \right)_{L^2(\Omega)} +\Phi_\nu(v_i^\nu; \theta_i^\nu) \leq \Phi_\nu(v_i^\nu; \omega), 
\\[1.5ex]
\mbox{for any $ \omega \in H^1(\Omega) $,}
\end{array}
\end{equation}
and
\begin{equation}\label{apx03_RC}
|\theta_i^\nu| \leq |\theta_{i -1}^\nu|_{L^\infty(\Omega)}, \mbox{ a.e.\ in $ \Omega $,}
\end{equation}
for $ i = 1, 2, 3, \ldots $\,.
\end{description}
We call the sequence $ \{ [v_i^\nu, \theta_i^\nu] = [w_i^\nu, \eta_i^\nu, \theta_i^\nu] \, | \, i \in \N \} \subset D_1 $ the solution to $ \mbox{(\hyperlink{(AP)_h^nu}{AP})}_h^\nu $, or the approximating solution for short.
\bigskip

Let us fix the time step $ 0 < h < 1 $. Then, our immediate task is to demonstrate the following theorem about the solvability of the problem $ \mbox{(\hyperlink{(AP)_h^nu}{AP})}_h^\nu $. 
\begin{thm}[solvability of the approximating problem]\label{Th.AP}
There exists a small constant $ h_*^\dag \in (0, 1) $ such that for $ 0 < h < h_*^\dag $, the approximating problem $ \mbox{(\hyperlink{(AP)_h^nu}{AP})}_h^\nu $ admits a unique solution $ \{ [v_i^\nu, \theta_i^\nu] = [w_i^\nu, \eta_i^\nu, \theta_i^\nu] \, | \, i \in \N \} \subset D_1 $, starting from any given initial value $ [v_0^\nu, \theta_0^\nu] = [w_0^\nu, \eta_0^\nu, \theta_0^\nu] \in D_1 $. Moreover, for $ 0 < h < h_*^\dag $, the approximating solution $ \{ [v_i^\nu, \theta_i^\nu] \} $ fulfills the following inequality of energy dissipation:
\begin{equation}\label{AP_dissipative}
\begin{array}{rl}
\ds \frac{1}{2h} |v_i^\nu -v_{i -1}^\nu|_{L^2(\Omega)^2}^2 & \ds + ~ \frac{1}{h} |{\textstyle \sqrt{\alpha_0(v_i^\nu)}}(\theta_i^\nu -\theta_{i -1}^\nu)|_{L^2(\Omega)}^2 
\\[2ex]
& \ds + ~ \mathscr{F}_\nu(v_i^\nu, \theta_i^\nu) \leq \mathscr{F}_\nu(v_{i -1}^\nu, \theta_{i -1}^\nu), ~ i = 1, 2, 3, \ldots.
\end{array}
\end{equation}
\end{thm}

The proof of Theorem \ref{Th.AP} will be quite extended, because some regularizations will be needed to relax the $ L^1 $-terms $ \nu |\nabla \theta_{i -1}^\nu|^2 [\nabla \beta](v_i^\nu) $, $ i \in \N $, in (\ref{apx01-02}). In view of this, we introduce one more relaxation index $ 0 < \varepsilon < 1 $, and we fix a large number $ M \in \N $ satisfying $ M > (N +2)/2 $. Additionally, for any $ \tilde{v} = [\tilde{w}, \tilde{\eta}] \in [H^1(\Omega) \cap L^\infty(\Omega)]^2 $, we define a relaxed convex function $ \Psi_\varepsilon^\nu(\tilde{v};\cdot) $ on $ L^2(\Omega) $ by setting
\begin{equation*}
\vartheta \in L^2(\Omega) \mapsto \Psi_\varepsilon^\nu(\tilde{v}; \vartheta) = \Psi_\varepsilon^\nu(\tilde{w}, \tilde{\eta}; \vartheta) := \Phi_\nu(\tilde{w}, \tilde{\eta}; \vartheta) +\frac{\varepsilon}{2} |\vartheta|_{H^M(\Omega)}^2 \in [0, \infty].
\end{equation*}
Also, for any $ 0 < \varepsilon < 1 $ and any $ \tilde{v} \in [H^1(\Omega) \cap L^\infty(\Omega)]^2 $, we denote by $ \partial \Psi_\varepsilon^\nu(\tilde{v};{}\cdot{}) $ the subdifferential of $ \Psi_\varepsilon^\nu(\tilde{v};{}\cdot{}) $ in the topology of $ L^2(\Omega) $. As can easily be verified, $ \Psi_\varepsilon^\nu(\tilde{v};{}\cdot{}) $ is proper l.s.c and convex on $ L^2(\Omega) $, and
\begin{equation}\label{apx04}
D(\Psi_\varepsilon^\nu(\tilde{v};{}\cdot{})) = H^M(\Omega) \subset W^{1, \infty}(\Omega)
\end{equation}
for all $ 0 < \varepsilon < 1 $ and $ \tilde{v} \in L^\infty(\Omega)^2 $. Based on this, we define a relaxed energy functional $ \mathscr{E}_\varepsilon^\nu $ on $ L^2(\Omega)^3 $, by setting
\begin{equation*}
\begin{array}{rcl}
[\tilde{v}, \tilde{\theta}] & = & [\tilde{w}, \tilde{\eta}, \tilde{\theta}] \in L^2(\Omega)^3 \mapsto \mathscr{E}_\varepsilon^\nu(\tilde{v}, \tilde{\theta}) = \mathscr{E}_\varepsilon^\nu(\tilde{w}, \tilde{\eta}, \tilde{\theta}) 
\\[1ex]
& := & V_{\rm D}^2(\tilde{v}) +\Gamma(\tilde{v}) +\mathscr{G}(\tilde{v}) +\Psi_\varepsilon^\nu(\tilde{v}; \tilde{\theta}) \mbox{ \ for any $ 0 < \varepsilon < 1 $.}
\end{array}
\end{equation*}
\medskip

Next, for any $ 0 < \varepsilon < 1 $, we denote by $ \mbox{(\hyperlink{(RX_eps)_h^nu}{RX$_\varepsilon $})}_h^\nu $ the relaxed system for $ \mbox{(\hyperlink{(AP)_h^nu}{AP})}_h^\nu $ prescribed as follows.
\begin{description}
\item[\textmd{$ \mbox{(\hypertarget{(RX_eps)_h^nu}{RX$_\varepsilon $})}_h^\nu $:}]For any initial value 
\begin{equation*}
\begin{array}{c}
[v_{\varepsilon, 0}^\nu, \theta_0^\nu] = [w_{\varepsilon, 0}^\nu, \eta_{\varepsilon, 0}^\nu, \theta_{\varepsilon, 0}^\nu] \in D_M := D_1 \cap [H^1(\Omega)^2 \times H^M(\Omega)] 
\\[1ex]
\mbox{with $v_{\varepsilon, 0}^\nu = [w_{\varepsilon, 0}^\nu, \eta_{\varepsilon, 0}^\nu] \in H^1(\Omega)^2$,}
\end{array}
\end{equation*}
find a sequence of triplets 
\begin{equation*}
\begin{array}{c}
\{ [v_{\varepsilon, i}^\nu, \theta_{\varepsilon, i}^\nu] = [w_{\varepsilon, i}^\nu, \eta_{\varepsilon, i}^\nu, \theta_{\varepsilon, i}^\nu] \, | \, i \in \N \} \subset D_M
\\[1ex]
\mbox{ with $ v_{\varepsilon, i}^\nu = [w_{\varepsilon, i}^\nu, \eta_{\varepsilon, i}^\nu] \in H^1(\Omega)^2 $ for $ i = 1, 2, 3, \ldots $}
\end{array}
\end{equation*}
such that
\begin{equation}\label{rx01-02}
\begin{array}{rl}
\multicolumn{2}{l}{\ds \hspace{-4ex} \frac{1}{h} (v_{\varepsilon, i}^\nu -v_{\varepsilon, i -1}^\nu) -{\mit \Delta}_{\rm N} v_{\varepsilon, i}^\nu +\partial \Gamma(v_{\varepsilon, i}) +[\nabla g](v_{\varepsilon, i}^\nu)}
\\[2ex]
& \ds +|\nabla \theta_{\varepsilon, i -1}^\nu| [\nabla \alpha](v_{\varepsilon, i}^\nu) +\nu |\nabla \theta_{\varepsilon, i -1}^\nu|^2 [\nabla \beta](v_{\varepsilon, i}^\nu) \ni 0 \mbox{ in  $ L^2(\Omega)^2 $}
\end{array}
\end{equation}
and
\begin{equation}\label{rx03}
\begin{array}{c}
\ds \frac{1}{h} \alpha_0(v_{\varepsilon, i}^\nu) (\theta_{\varepsilon, i}^\nu -\theta_{\varepsilon, i -1}^\nu) +\partial \Psi_\varepsilon^\nu(v_{\varepsilon, i}^\nu; \theta_{\varepsilon, i}^\nu) \ni 0 \mbox{ in $ L^2(\Omega) $,}
\end{array}
\end{equation}
for $ i = 1, 2, 3, \ldots\, $.
\end{description}
\begin{rem}\label{Rem.RX01}
\begin{em}
In the relaxed system $ \mbox{(\hyperlink{(RX_eps)_h^nu}{RX$_\varepsilon $})}_h^\nu$, we note that the inclusions (\ref{rx01-02}) and (\ref{rx03}) are relaxed versions of the variational inequalities (\ref{apx01-02}) and (\ref{apx03}), respectively, and these are expressed in the reduced forms by means of $ L^2 $-subdifferentials. As mentioned in Remark \ref{Rem.1.1}, such reductions may not be available for (\ref{rx01-02}). However, in light of (\ref{apx04}), we observe that
\begin{equation*}
\theta_{\varepsilon, i -1}^\nu \in H^M(\Omega) \subset W^{1, \infty}(\Omega) \mbox{ for $ i \in \N $,}
\end{equation*}
and this enables us to suppose that
\begin{equation*}
|\nabla \theta_{\varepsilon, i -1}^\nu| \in L^\infty(\Omega) \mbox{ \ and \ } \nu |\nabla \theta_{\varepsilon, i -1}^\nu|^2 [\nabla \beta](v_i^\nu) \in L^2(\Omega)^N \mbox{ for $ i \in \N $.}
\end{equation*}
Hence, the relaxed system $ \mbox{(\hyperlink{(RX_eps)_h^nu}{RX$_\varepsilon $})}_h^\nu$ will be in the applicable scope of the general theories of $ L^2 $-subdifferentials. 
\end{em}
\end{rem}

Now we fix $ 0 < \varepsilon < 1 $ and devote the remaining part of this section to confirming some auxiliary lemmas concerned with the key properties of the relaxed system $ \mbox{(\hyperlink{(RX_eps)_h^nu}{RX$_\varepsilon $})}_h^\nu$.

\begin{lem}\label{Lem.RX01}
For arbitrary $ \theta_0^\dag \in W^{1, \infty}(\Omega) $ and $ v_0^\dag \in H^1(\Omega)^2 $, consider an auxiliary inclusion
\begin{align}
\ds \frac{v -v_0^\dag}{h} -{\mit \Delta}_{\rm N} v & +\partial \Gamma(v) +[\nabla g](\mathscr{T}_0^1v) 
\nonumber
\\
& +|\nabla \theta_0^\dag| [\nabla \alpha](v) +\nu |\nabla \theta_0^\dag|^2 [\nabla \beta](v) \ni 0 \mbox{ in $ L^2(\Omega)^2 $.}
\label{rx04}
\end{align}
There exists a small constant $ h_0^\dag \in (0, 1) $, depending only on $ |g|_{C^2([0, 1]^2)} $, and for any $ 0 < h < h_0^\dag $, the inclusion (\ref{rx04}) admits a unique solution $ v \in H^1(\Omega)^2 $. 
\end{lem}

\noindent
\textbf{Proof.} 
Let $ S_h^\dag : H^1(\Omega)^2 \rightarrow H^1(\Omega)^2 $ be an operator that maps any $ v^\dag \in H^1(\Omega)^2 $ to a unique minimizer $ S_h^\dag v^\dag \in H^1(\Omega)^2 $ of a proper l.s.c.\ and strictly convex function on $ H^1(\Omega)^2 $, defined as
\begin{equation*}
\begin{array}{rl}
\varpi \in H^1(\Omega)^2 \mapsto \ds \frac{1}{2h} |\varpi -v_0^\dag|_{L^2(\Omega)^2}^2 & \hspace{-1ex} + V_{\rm D}^2(\varpi) +\Gamma(\varpi) +([\nabla g](\mathscr{T}_0^1 v^\dag), \varpi)_{L^2(\Omega)} 
\\
& \ds \hspace{-1ex} +\int_\Omega \left( |\nabla \theta_0^\dag| \alpha(\varpi) +\nu |\nabla \theta_0^\dag|^2 \beta(\varpi) \right) \, dx \in (-\infty, \infty].
\end{array}
\end{equation*}
On this basis, let us take two functions $ v_k^\dag \in H^1(\Omega)^2 \cap D(\Gamma) $, $ k = 1, 2 $, and consider the smallness condition on $ h $ for $ S_h^\dag $ to become contractive. From the definition of $ S_h^\dag $, the functions $ v_k := S_h^\dag v_k^\dag \in H^1(\Omega)^2 $, $ k = 1, 2 $, fulfill
\begin{equation}\label{rx05}
\begin{array}{rl}
\ds \frac{v_k -v_0^\dag}{h} -{\mit \Delta}_{\rm N} v_k \hspace{-1.5ex} \ & +\partial \Gamma(v_k) +[\nabla g](\mathscr{T}_0^1v_k^\dag) 
\\
& +|\nabla \theta_0^\dag| [\nabla \alpha](v_k) +\nu |\nabla \theta_0^\dag|^2 [\nabla \beta](v_k) \ni 0 \mbox{ in $ L^2(\Omega)^2 $, $ k = 1, 2 $,}
\end{array}
\end{equation}
respectively. Here, taking differences between two inclusions in (\ref{rx05}) and multiplying both sides of the result by $ v_1 -v_2 $, we infer from (\hyperlink{A1}{A1})--(\hyperlink{A3}{A3}) that
\begin{equation*}
\frac{1}{h} |v_1 -v_2|_{L^2(\Omega)^2}^2 +|\nabla (v_1 -v_2)|_{L^2(\Omega)^{2 \times N}}^2 \leq |g|_{C^2([0, 1]^2)}|v_1^\dag -v_2^\dag|_{L^2(\Omega)^2} |v_1 -v_2|_{L^2(\Omega)^2}.
\end{equation*}
Subsequently, by using Young's inequality, 
\begin{equation}\label{rx06}
\frac{1}{2h} |v_1 -v_2|_{L^2(\Omega)^2}^2 +|\nabla (v_1 -v_2)|_{L^2(\Omega)^{2 \times N}}^2 \leq \frac{h}{2} |g|_{C^2([0, 1]^2)}^2 |v_1^\dag -v_2^\dag|_{L^2(\Omega)^2}^2.
\end{equation}
So if we assume that
\begin{equation}\label{rx07}
0 < h < h_0^\dag := \frac{1}{2 (1 \vee |g|_{C^2([0, 1]^2)})},
\end{equation}
then it can be seen from (\ref{rx06}) that $ S_h^\dag $ becomes a contraction mapping from $ H^1(\Omega)^2 $ into itself. Therefore, applying Banach's fixed-point theorem, we find a unique fixed point $ v_*^\dag \in H^1(\Omega)^2 $ of $ S_h^\dag $ under (\ref{rx07}). The identity $ v_*^\dag = S_h^\dag v_*^\dag $ implies that $ v_*^\dag$ solves the auxiliary inclusion (\ref{rx04}). \hfill \ \qed

\begin{lem}\label{Lem.RX03}
Let us assume $ 0 < h < h_0^\dag $ with the constant $ h_0^\dag \in (0, 1) $ as in Lemma \ref{Lem.RX01}. Let $ \theta_0^\dag \in W^{1, \infty}(\Omega) $ and $ v_0^\dag = [w_0^\dag, \eta_0^\dag] \in H^1(\Omega)^2 $ be fixed functions, and let $ v = [w, \eta] \in H^1(\Omega)^2 $ be the unique solution to the auxiliary inclusion (\ref{rx04}). Let $ o_*, \iota_* \in D(\partial \gamma) $ be constants as in (\hyperlink{A4}{A4}), and let $ \check{r},  \hat{r} \in \R^2 $ be two constant vectors given by
\begin{equation*}
\check{r} := [o_*, 0] \mbox{ \ and \ } \hat{r} := [\iota_*, 1].
\end{equation*}
If
\begin{equation}\label{rxInit}
\check{r} \leq v_0^\dag \leq \hat{r}, \mbox{ i.e.\ } v_0^\dag \in [o_*, \iota_*] \times [0, 1], \mbox{ a.e.\ in $ \Omega $,}
\end{equation}
then the following ordering property is preserved:
\begin{equation*}
\check{r} \leq v \leq \hat{r}, \mbox{ i.e.\ } v \in [o_*, \iota_*] \times [0, 1], \mbox{ a.e.\ in $ \Omega $.}
\end{equation*}
\end{lem}

\noindent
\textbf{Proof.}
By (\ref{a4-1}) in (\hyperlink{A4}{A4}), we find two elements $ \check{o}_* \in \partial \gamma(o_*) $ and $ \hat{\iota}_* \in \partial \gamma(\iota_*) $ such that
\begin{equation}\label{rx11-1}
\check{o}_* +g_w(o_*, \tilde{\eta}) \leq 0 \mbox{ \ and \ } \hat{\iota}_* +g_w(\iota_*, \tilde{\eta}) \geq 0, \mbox{ for any $ \tilde{\eta} \in [0, 1] $.}
\end{equation}
Setting
\begin{equation*}
\check{r}^* := [\check{o}_*, 0] \mbox{ \ and \ } \hat{r}^* := [\hat{\iota}_*, 0] \mbox{ \ in $ \R^2 $,}
\end{equation*}
it follows from (\ref{a4-1})--(\ref{a4-2}) in (\hyperlink{A4}{A4}) and (\ref{rxInit})--(\ref{rx11-1}) that
\begin{equation}\label{rx12}
[\check{r}, \check{r}^*] \in \partial \Gamma \mbox{ \ and \ } [\hat{r}, \hat{r}^*] \in \partial \Gamma \mbox{ in $ L^2(\Omega)^2 $}
\end{equation}
and
\begin{equation}\label{rx13}
\begin{array}{ll} 
\ds \frac{\check{r} -v_0^\dag}{h} & - ~ {\mit \Delta}_{\rm N} \check{r} +\check{r}_* +\left[ \hspace{-0.5ex} \begin{array}{c} g_w(o_*, \mathscr{T}_0^1 \eta) \\[0.25ex] g_\eta(\mathscr{T}_0^1 w, 0) \end{array} \hspace{-0.5ex} \right]
\\[2ex]
& + ~ |\nabla \theta_0^\dag| \left[ \hspace{-0.5ex} \begin{array}{c} \alpha_w(o_*, \eta) \\ \alpha_\eta(w, 0) \end{array} \hspace{-0.5ex} \right] +\nu |\nabla \theta_0^\dag|^2
\left[ \hspace{-0.5ex} \begin{array}{c} \beta_w(o_*, \eta) \\ \beta_\eta(w, 0) \end{array} \hspace{-0.5ex} \right] \leq \left[ \hspace{-0.5ex} \begin{array}{c} 0 \\[0.25ex] 0 \end{array} \hspace{-0.5ex} \right], \mbox{ \ a.e.\ in $ Q $,}
\end{array}
\end{equation}
\begin{equation}\label{rx14}
\begin{array}{ll} 
\ds \frac{\hat{r} -v_0^\dag}{h} & - ~ {\mit \Delta}_{\rm N} \hat{r} +\hat{r}^* +\left[ \hspace{-0.5ex} \begin{array}{c} g_w(\iota_*, \mathscr{T}_0^1 \eta) \\[0.25ex] g_\eta(\mathscr{T}_0^1 w, 1) \end{array} \hspace{-0.5ex} \right]
\\[2ex]
& + ~ |\nabla \theta_0^\dag| \left[ \hspace{-0.5ex} \begin{array}{c} \alpha_w(\iota_*, \eta) \\[0.25ex] \alpha_\eta(w, 1) \end{array} \hspace{-0.5ex} \right] +\nu |\nabla \theta_0^\dag|^2 \left[ \hspace{-0.5ex} \begin{array}{c} \beta_w(\iota_*, \eta) \\[0.25ex] \beta_\eta(w, 1) \end{array} \hspace{-0.5ex} \right] \geq \left[ \hspace{-0.5ex} \begin{array}{c} 0 \\[0.25ex] 0 \end{array} \hspace{-0.5ex} \right], \mbox{ \ a.e.\ in $ Q $.}
\end{array}
\end{equation}
Take the difference from (\ref{rx13}) to (\ref{rx04}) and multiply both sides of the result by $ [\check{r} -{v}]^+ $. Then, by virtue of (\hyperlink{A1}{A1})--(\hyperlink{A4}{A4}), (\ref{Gamma(v)}), (\ref{rx11-1})--(\ref{rx12}),
\begin{equation}\label{rx14-1}
\ds \frac{1}{h} |[\check{r} -{v}]^+|_{L^2(\Omega)^2}^2 +|\nabla [\check{r} -{v}]^+|_{L^2(\Omega)^{2 \times N}}^2 
\leq \ds |g|_{C^2([0, 1]^2)} |[\check{r} -{v}]^+|_{L^2(\Omega)^2}^2.
\end{equation}
As well as, we also have:
\begin{equation}\label{rx14-2}
\ds \frac{1}{h} |[{v} -\hat{r}]^+|_{L^2(\Omega)^2}^2 +|\nabla [{v} -\hat{r}]^+|_{L^2(\Omega)^{2 \times N}}^2 
\leq \ds |g|_{C^2([0, 1]^2)} |[{v} -\hat{r}]^+|_{L^2(\Omega)^2}^2,
\end{equation}
by taking the difference from (\ref{rx04}) to (\ref{rx14}) and multiplying both sides of the result by $ [{v} -\hat{r}]^+ $. 

On account of (\ref{rx14-1}) (resp. (\ref{rx14-2})), the inequality
\begin{center}
$ \check{r} \leq {v} $ \ a.e.\ in $ \Omega $ (resp. $ {v} \leq \hat{r} $ \ a.e.\ in $ \Omega $)
\end{center}
can be inferred by using the assumption $ h \in (0, h_0^\dag) $. 
\qed
\medskip

\begin{lem}\label{Lem.RX04}
Let $ h_0^\dag \in (0, 1) $ be a constant as in Lemma \ref{Lem.RX01}, and let $ o_*, \iota_* \in \R $ be constants as in (\hyperlink{A4}{A4}). Let $ \theta_0^\dag \in W^{1, \infty}(\Omega) $ and $ v_0^\dag = [w_0^\dag, \eta_0^\dag] \in H^1(\Omega)^2 $ be fixed functions. If $ 0 < h < h_0^\dag $, and the function $ v_0^\dag = [w_0^\dag, \eta_0^\dag] $ satisfies
\begin{equation*}
v_0^\dag \in [o_*, \iota_*] \times [0, 1], \mbox{ i.e.\ } o_* \leq w_0^\dag \leq \iota_* \mbox{ and } 0 \leq \eta_0^\dag \leq 1, \mbox{ a.e.\ in $ \Omega $,}
\end{equation*}
then the (nontruncated) inclusion
\begin{equation*}
\begin{array}{rl}
\ds \frac{v -v_0^\dag}{h} -{\mit \Delta}_{\rm N} v \hspace{-1.5ex} \ & + \ \partial \Gamma(v) +[\nabla g](v) 
\\
& + \ |\nabla \theta_0^\dag| [\nabla \alpha](v) +\nu |\nabla \theta_0^\dag|^2 [\nabla \beta](v) \ni 0 \mbox{ in $ L^2(\Omega)^2 $}
\end{array}
\end{equation*}
admits a unique solution $ v = [w, \eta] \in H^1(\Omega)^2 $. Moreover, 
\begin{equation*}
v \in [o_*, \iota_*] \times [0, 1], \mbox{ i.e.\ } o_* \leq w \leq \iota_* \mbox{ and } 0 \leq \eta \leq 1, \mbox{ a.e.\ in $ \Omega $.}
\end{equation*}
\end{lem}
\medskip

\noindent
\textbf{Proof.}
This lemma is immediately deduced by combining the conclusions of Lemmas \ref{Lem.RX01} and \ref{Lem.RX03}.
\hfill \qed
\medskip

\begin{lem}\label{Lem.RX05}
Let $ v^\dag \in H^1(\Omega)^2 $ and $ \theta_0^\dag \in H^M(\Omega) $ be fixed functions. Then the inclusion
\begin{equation}\label{rx30}
\alpha_0(v^\dag) \frac{\theta -\theta_0^\dag}{h} +\partial \Psi_\varepsilon^\nu(v^\dag; \theta) \ni 0 \mbox{ \ in $ L^2(\Omega) $}
\end{equation}
admits a unique solution $ \theta \in H^M(\Omega) $. 
\end{lem}

\noindent
\textbf{Proof. }
The inclusion (\ref{rx30}) corresponds to the Euler-Lagrange equation for the following proper l.s.c.\ and convex function on $ L^2(\Omega) $:
\begin{equation*}
\theta \in L^2(\Omega) \mapsto \frac{1}{2h} |\textstyle{\sqrt{\alpha_0(v^\dag)}} (\theta -\theta_0^\dag)|_{L^2(\Omega)}^2 +\Psi_\varepsilon^\nu(v^\dag; \theta) \in [0, \infty].
\end{equation*}
Since this function is coercive and strictly convex on $ L^2(\Omega) $, the lemma is a direct consequence of the general theory of convex analysis (cf.\ \citep[Chapter II]{ET}). \hfill \qed

\begin{lem}[solvability of the relaxed system]\label{Lem.RX06}
Assume $ 0 < h < h_1^\dag := h_0^\dag/2 $ with constant $ h_0^\dag \in (0, 1) $ as in Lemma \ref{Lem.RX01}. Then, the relaxed system $ \mbox{(\hyperlink{(RX_eps)_h^nu}{RX$_\varepsilon $})}_h^\nu$ admits a unique solution $ \{ [v_{\varepsilon, i}^h, \theta_{\varepsilon, i}^h] = [w_{\varepsilon, i}^h, \eta_{\varepsilon, i}^h, \theta_{\varepsilon, i}^h] \, | \, i \in \N \} \subset D_M $, starting from any initial value $ [v_{\varepsilon, 0}^h, \theta_{\varepsilon, 0}^h] = $ $ [w_{\varepsilon, 0}^h, \eta_{\varepsilon, 0}^h, \theta_{\varepsilon, 0}^h] $ $ \in D_M $, and moreover,
\begin{equation}\label{rxEgyIneq}
\begin{array}{rl}
\ds \frac{1}{2h} |v_{\varepsilon, i}^h -v_{\varepsilon, i -1}^h|_{L^2(\Omega)^2}^2 & \ds + ~ \frac{1}{h} \left| {\textstyle \sqrt{\alpha_0(v_{\varepsilon, i}^h)}}(\theta_{\varepsilon, i} -\theta_{\varepsilon, i -1}^h) \right|_{L^2(\Omega)}^2
\\[2ex]
& + ~ \mathscr{E}_\varepsilon^\nu(v_{\varepsilon, i}^h, \theta_{\varepsilon, i}^h) \leq \mathscr{E}_\varepsilon^\nu(v_{\varepsilon, i -1}^h, \theta_{\varepsilon, i -1}^h), \mbox{ for } i = 1, 2, 3, \ldots.
\end{array}
\end{equation}
\end{lem}

\noindent
\textbf{Proof.} 
On the basis of Lemmas \ref{Lem.RX04}--\ref{Lem.RX05}, we can obtain a unique solution $ \{ [v_{\varepsilon, i}^h, \theta_{\varepsilon, i}^h] = [w_{\varepsilon, i}^h, \eta_{\varepsilon, i}^h, \theta_{\varepsilon, i}^h] \, | \, i \in \N \} \subset H^1(\Omega)^2 \times H^M(\Omega) $ of the relaxed system $ \mbox{(\hyperlink{(RX_eps)_h^nu}{RX$_\varepsilon $})}_h^\nu$ as follows:

\begin{description}
\item[\textmd{(Step 0)}]Let $ i = 1 $ and fix $ [v_{\varepsilon, 0}^h, \theta_{\varepsilon, 0}^h] = [w_{\varepsilon, 0}^h, \eta_{\varepsilon, 0}^h, \theta_{\varepsilon, 0}^h] \in D_M $.
\item[\textmd{(Step 1)}]Obtain a unique solution $ v_{\varepsilon, i}^h = [w_{\varepsilon, i}^h, \eta_{\varepsilon, i}^h] \in H^1(\Omega)^2 $ to (\ref{rx01-02}) with the range constraint $ v_{\varepsilon, i} \in [o_*, \iota_*] \times [0, 1] $ a.e.\ in $ \Omega $ by applying Lemma \ref{Lem.RX04} as the case when $ \theta_0^\dag = \theta_{\varepsilon, i -1}^h $, $ v_0^\dag = v_{\varepsilon, i -1}^h $, and $ v = v_{\varepsilon, i}^h $.
\item[\textmd{(Step 2)}]Obtain a unique solution $ \theta_{\varepsilon, i}^h \in H^M(\Omega) $ to (\ref{rx30}) by applying Lemma \ref{Lem.RX05} as the case when $ v^\dag = v_{\varepsilon, i}^h $, $ \theta_0^\dag = \theta_{\varepsilon, i -1}^h $, and $ \theta = \theta_{\varepsilon, i}^\nu $.
\item[\textmd{(Step 3)}]Iterate the value of $ i $, i.e., $ i \leftarrow i +1 $, and return to step 1.
\end{description}

Next we verify the inequality (\ref{rxEgyIneq}). Multiply both sides of (\ref{rx01-02}) by $ v_{\varepsilon, i}^h -v_{\varepsilon, i -1}^h $. By using (\hyperlink{A1}{A1}), (\ref{Gamma(v)}), and Young's inequality, we have
\begin{equation}\label{rx40}
\begin{array}{ll}
\multicolumn{2}{l}{\ds\frac{1}{h} |v_{\varepsilon, i}^h -v_{\varepsilon, i -1}^h|_{L^2(\Omega)}^2 +\frac{1}{2} |\nabla v_{\varepsilon, i}^h|_{L^2(\Omega)^{2 \times N}}^2 -\frac{1}{2} |\nabla v_{\varepsilon, i -1}^h|_{L^2(\Omega)^{2 \times N}}^2}
\\[2ex]
\qquad \quad & \ds +\int_\Omega [\nabla g](v_{\varepsilon, i}^h) \cdot (v_{\varepsilon, i}^h -v_{\varepsilon, i -1}^h) \, dx
\\[2ex]
& \ds +\int_\Omega |\nabla \theta_{\varepsilon, i -1}^h| [\nabla \alpha](\theta_{\varepsilon, i}^h) \cdot (v_{\varepsilon, i}^h -v_{\varepsilon, i -1}^h)  \, dx
\\[2ex]
& \ds +\nu \int_\Omega |\nabla \theta_{\varepsilon, i -1}^h|^2 [\nabla \beta](\theta_{\varepsilon, i}^h) \cdot (v_{\varepsilon, i}^h -v_{\varepsilon, i -1}^h)  \, dx
\\[2ex]
& \ds +\Gamma(v_{\varepsilon, i}^h) -\Gamma(v_{\varepsilon, i -1}^h) \leq 0, \mbox{ for } i = 1, 2, 3, \ldots.
\end{array}
\end{equation}
Invoking (\hyperlink{A2}{A2})--(\hyperlink{A3}{A3}), we compute that
\begin{equation}\label{rx41}
\begin{array}{rl}
& {
\ds \int_\Omega [\nabla g](v_{\varepsilon, i}^h) \cdot (v_{\varepsilon, i}^h -v_{\varepsilon, i -1}^h) \, dx \geq \int_\Omega g(v_{\varepsilon, i}^h) \, dx -\int_\Omega g(v_{\varepsilon, i -1}^h) \, dx
}
\\[2ex]
& \ds \qquad \quad +\int_\Omega \left( [\nabla g](v_{\varepsilon, i}^h) -[\nabla g](v_{\varepsilon, i -1}^h)  \right) \cdot (v_{\varepsilon, i}^h -v_{\varepsilon, i -1}^h) \, dx 
\\[2ex]
& \ds \qquad \quad -\frac{1}{2} |g|_{C^2([0, 1]^2)} |v_{\varepsilon, i}^h -v_{\varepsilon, i -1}^h|_{L^2(\Omega)^2}^2
\\[2ex]
\geq & \ds \mathscr{G}(v_{\varepsilon, i}^h) -\mathscr{G}(v_{\varepsilon, i -1}^h) -\frac{3}{2} |g|_{C^2([0, 1]^2)} |v_{\varepsilon, i}^h -v_{\varepsilon, i -1}^h|_{L^2(\Omega)^2}^2 
\end{array}
\end{equation}
and

\begin{equation}\label{rx42}
\begin{array}{rl}
& {\ds \int_\Omega |\nabla \theta_{\varepsilon, i -1}^h| [\nabla \alpha](v_{\varepsilon, i}^h) \cdot (v_{\varepsilon, i}^h -v_{\varepsilon, i -1}^h) \, dx}
\\[2ex]
& \quad \qquad \ds +\nu \int_\Omega |\nabla \theta_{\varepsilon, i -1}^h|^2 [\nabla \beta](v_{\varepsilon, i}^h) \cdot (v_{\varepsilon, i}^h -v_{\varepsilon, i -1}^h) \, dx
\\[2ex]
\geq & \ds \int_\Omega \alpha(v_{\varepsilon, i}^h) |\nabla \theta_{\varepsilon, i -1}^h| \, dx -\int_\Omega \alpha(v_{\varepsilon, i -1}^h) |\nabla \theta_{\varepsilon, i -1}^h| \, dx
\\[2ex]
& \ds \qquad \quad +\nu \int_\Omega \beta(v_{\varepsilon, i}^h) |\nabla \theta_{\varepsilon, i -1}^h|^2 \, dx - \nu \int_\Omega \beta(v_{\varepsilon, i -1}^h) |\nabla \theta_{\varepsilon, i -1}^h|^2 \, dx 
\end{array}
\end{equation}
for $ i = 1, 2, 3, \ldots $. On the basis of (\ref{rx40})--(\ref{rx42}), it is deduced that
\begin{equation}\label{rx43}
\begin{array}{rl}
& \ds \left( 1 -\frac{3}{2} |g|_{C^2([0, 1]^2)} h \right) \frac{1}{h} |v_{\varepsilon, i}^h -v_{\varepsilon, i -1}^h|_{L^2(\Omega)^2}^2 
\\[2ex]
& \ds \quad \qquad +V_{\rm D}^2(v_{\varepsilon, i}^h) +\Gamma(v_{\varepsilon, i}^h) +\mathscr{G}(v_{\varepsilon, i}^h)
\\[2ex]
& \ds \quad \qquad +\int_\Omega \alpha(v_{\varepsilon, i}^h)|\nabla \theta_{\varepsilon, i -1}^h| \, dx +\nu \int_\Omega \beta(v_{\varepsilon, i}^h) |\nabla \theta_{\varepsilon, i -1}^h|^2 \, dx
\\[2ex]
\leq & V_{\rm D}^2(v_{\varepsilon, i -1}^h) +\Gamma(v_{\varepsilon, i -1}^h) +\mathscr{G}(v_{\varepsilon, i -1}^h)
\\[2ex]
& \ds \quad \qquad +\int_\Omega \alpha(v_{\varepsilon, i -1}^h)|\nabla \theta_{\varepsilon, i -1}^h| \, dx +\nu \int_\Omega \beta(v_{\varepsilon, i -1}^h) |\nabla \theta_{\varepsilon, i -1}^h|^2 \, dx, 
\\[2ex]
& \mbox{for } i = 1, 2, 3, \ldots.
\end{array}
\end{equation}
Meanwhile, by multiplying both sides of (\ref{rx03}) by $ \theta_{\varepsilon, i}^h -\theta_{\varepsilon, i -1}^h $, it can seen that
\begin{equation}\label{rx44}
\begin{array}{rl}
& \ds \frac{1}{h} |{\textstyle \sqrt{\alpha_0(v_{\varepsilon, i}^h)}} (\theta_{\varepsilon, i}^h -\theta_{\varepsilon, i -1}^h)|_{L^2(\Omega)}^2 
\\[2ex]
& \ds \quad \qquad +\int_\Omega \alpha(v_{\varepsilon, i}^h)|\nabla \theta_{\varepsilon, i}^h| \, dx +\nu \int_\Omega \beta(v_{\varepsilon, i}^h) |\nabla \theta_{\varepsilon, i}^h|^2 \, dx
\\[2ex]
& \ds \quad \qquad -\int_\Omega \alpha(v_{\varepsilon, i}^h)|\nabla \theta_{\varepsilon, i -1}^h| \, dx -\nu \int_\Omega \beta(v_{\varepsilon, i}^h) |\nabla \theta_{\varepsilon, i -1}^h|^2 \, dx
\\[2ex]
& \ds \quad \qquad +\frac{\varepsilon}{2} |\theta_{\varepsilon, i}^h|_{H^M(\Omega)}^2
\leq \frac{\varepsilon}{2} |\theta_{\varepsilon, i -1}^h|_{H^M(\Omega)}^2, \mbox{ for } i= 1, 2, 3, \ldots.
\end{array}
\end{equation}
Since 
\begin{equation}\label{1/2}
1 -\frac{3}{2} |g|_{C^2([0, 1]^2)} h > \frac{1}{2},\mbox{ \ if } 0 < h < h_1^\dag < \frac{1}{3(1 \vee |g|_{C^2([0, 1]^2)})},
\end{equation}
the required inequality (\ref{rxEgyIneq}) is obtained by taking the sum of (\ref{rx43}) and (\ref{rx44}) and applying (\ref{1/2}). \hfill \qed

\section{Solvability of approximating problems}
\ \vspace{-3ex}

In this section, we fix $ \nu > 0 $ and $ 0 < h < h_1^\dag $ with the constant as in Lemma \ref{Lem.RX06} and prove Theorem \ref{Th.AP} concerning the approximating problem $ \mbox{(\hyperlink{(AP)_h^nu}{AP})}_h^\nu $. The proof of this theorem is divided into two parts, which respectively concerned with ``the existence'' and ``the uniqueness and energy dissipation.''

\subsection*{Existence of approximating solutions}
\ \vspace{-3ex}

First, we prepare some lemmas for the limiting observations of the relaxed systems $ \mbox{(\hyperlink{(RX_eps)_h^nu}{RX$_\varepsilon $})}_h^\nu$ as $ \varepsilon \searrow 0 $.

\begin{lem}\label{Lem.RX07}
Assume $ v^\dag \in [H^1(\Omega) \cap L^\infty(\Omega))]^2 $, $ \{ v_\varepsilon^\dag \, | \, 0 < \varepsilon < 1 \} \subset  [H^1(\Omega) \cap L^\infty(\Omega))]^2 $, and
\begin{equation}\label{rx50}
\left\{ ~ \parbox{10cm}{
$ v_\varepsilon^\dag \to v^\dag $ in the pointwise sense a.e.\ in $ \Omega $ as $ \varepsilon \searrow 0 $,
\\[1ex]
$ \{ v_\varepsilon^\dag \, | \, 0 < \varepsilon < 1 \} $ is bounded in $ L^\infty(\Omega)^2 $.} 
\right.
\end{equation}
Then, for the sequence of convex functions $ \{ \Psi_\varepsilon^\nu(v_\varepsilon^\dag;\cdot\,) \, | \, 0 < \varepsilon < 1 \} $, it holds that $ \Psi_\varepsilon^\nu(v_\varepsilon^\dag;\cdot\,) $ $\to  \Phi_\nu(v^\dag;\cdot\,) $ on $ L^2(\Omega) $ in the sense of Mosco \citep{Mosco} as $ \varepsilon \searrow 0 $, i.e.,
\begin{description}
\item[$ \mbox{\textmd{\em (\hypertarget{(m1)_eps^nu}{m1})}}_\varepsilon^\nu $](lower bound) $ \ds \liminf_{\varepsilon \searrow 0} \Psi_\varepsilon^\nu(v_\varepsilon^\dag; \theta_\varepsilon^\dag) \geq \Phi_\nu(v^\dag; \theta^\dag) $ if $ \theta^\dag \in L^2(\Omega) $, $ \{ \theta_\varepsilon^\dag \, | \, 0 < \varepsilon < 1 \} \subset L^2(\Omega) $, and $ \theta_\varepsilon^\dag \to \theta^\dag $ weakly in $ L^2(\Omega) $ as $ \varepsilon \searrow 0 $; 
\item[$ \mbox{\textmd{\em (\hypertarget{(m2)_eps^nu}{m2})}}_\varepsilon^\nu $](optimality) for any $ \theta^\ddag \in H^1(\Omega) $, there exists a sequence $ \{ \theta_\varepsilon^\ddag \, | \, 0 < \varepsilon < 1 \} \subset H^M(\Omega) $ such that $ \theta_\varepsilon^\ddag \to \theta^\ddag $ in $ L^2(\Omega) $ and $ \ds \Psi_\varepsilon^\nu(v_\varepsilon^\dag; \theta_\varepsilon^\ddag) \to \Phi_\nu(v^\dag; \theta^\ddag) $ as $ \varepsilon \searrow 0 $. 
\end{description}
Additionally, in light of Remark \ref{Rem.MG}, the above Mosco convergence implies $ \Gamma $-convergence on $ L^2(\Omega) $ as $ \varepsilon \searrow 0 $. 
\end{lem}

\noindent
\textbf{Proof.}
To verify condition $ \mbox{(\hyperlink{(m1)_eps^nu}{m1})}_\varepsilon^\nu $, it is enough to consider only the case in which $ \liminf_{\varepsilon \searrow 0} \Psi_\varepsilon^\nu(v_\varepsilon^\dag; \theta_\varepsilon^\dag) < \infty $, because the other case is trivial. On this basis, we suppose that
\begin{equation}\label{rx51}
\left\{ ~ 
\parbox{12.5cm}{
$ 1 > \varepsilon_1^\dag > \cdots > \varepsilon_m^\dag \searrow 0 $ and $ \theta_{\varepsilon_m^\dag}^\dag \to \theta^\dag $ weakly in $ H^1(\Omega) $ \ as $ m \to \infty $,
\\[1ex]
$ \ds \liminf_{\varepsilon \searrow 0} \Psi_\varepsilon^\nu(v_\varepsilon^\dag; \theta_\varepsilon^\dag) = \lim_{m \to \infty} \Psi_{\varepsilon_m^\dag}^\nu(v_{\varepsilon_m^\dag}^\dag; \theta_{\varepsilon_m^\dag}^\dag) $.
} \right.
\end{equation}
Then, as a result of (\hyperlink{A3}{A3}), (\ref{rx50})--(\ref{rx51}), and Lebesgue's dominated convergence theorem, it is inferred that
\begin{equation*}
\left\{ \parbox{5.5cm}{
$ \alpha(v_{\varepsilon_m^\dag}^\dag) \nabla \theta_{\varepsilon_m^\dag}^\dag \to \alpha(v^\dag) \nabla \theta^\dag $
\\[1ex]
$ {\textstyle \sqrt{\beta(v_{\varepsilon_m^\dag}^\dag)}} \nabla \theta_{\varepsilon_m^\dag}^\dag \to {\textstyle \sqrt{\beta(v^\dag)}} \nabla \theta^\dag $
} \right. \mbox{weakly in $ L^2(\Omega)^N $ as $ m \to \infty $.}
\end{equation*}
Therefore, keeping in mind the lower semi-continuity of the norms, it can be seen that
\begin{equation*}
\begin{array}{rl}
\multicolumn{2}{l}{\ds \liminf_{\varepsilon \searrow 0} \Psi_\varepsilon^\nu(v_\varepsilon^\dag; \theta_\varepsilon^\dag) = \lim_{m \to \infty} \Psi_{\varepsilon_m^\dag}^\nu(v_{\varepsilon_m^\dag}^\dag; \theta_{\varepsilon_m^\dag}^\dag)}
\\[2ex]
\quad \qquad \geq & \ds \liminf_{m \to \infty} \, \bigl| \alpha(v_{\varepsilon_m^\dag}^\dag) \nabla \theta_{\varepsilon_m^\dag}^\dag \bigr|_{L^1(\Omega; \R^N)} +\nu \liminf_{m \to \infty} \,  \bigl| {\textstyle \sqrt{\beta(v_{\varepsilon_m^\dag}^\dag)}} \nabla \theta_{\varepsilon_m^\dag}^\dag \bigr|_{L^2(\Omega)^N}^2
\\[2ex]
\quad \qquad \geq & \ds \bigl| \alpha(v^\dag) \nabla \theta^\dag \bigr|_{L^1(\Omega; \R^N)} +\nu \bigl| {\textstyle \sqrt{\beta(v^\dag)}} \nabla \theta^\dag \bigr|_{L^2(\Omega)^N}^2
\\[2ex]
= & \ds \Phi_\nu(v^\dag; \theta^\dag).
\end{array}
\end{equation*}
Thus, we have verified the condition $ \mbox{(\hyperlink{(m1)_eps^nu}{m1})}_\varepsilon^\nu $. 
\medskip

Next we prove $ \mbox{(\hyperlink{(m2)_eps^nu}{m2})}_\varepsilon^\nu $. For any $ \theta^\ddag \in H^1(\Omega) $, let $ \{ \tilde{\theta}_m^\ddag \, | \, m \in \N \} \subset C^\infty(\overline{\Omega}) $ be the standard approximating sequence of $ \theta^\ddag $, such that
\begin{equation*}
\tilde{\theta}_m^\ddag \to \theta^\ddag \mbox{ in $ H^1(\Omega) $ as $ m \to \infty $.}
\end{equation*}
Also, let us take a sequence $ \{ \varepsilon_m^\ddag \, | \, m \in \Z, ~ m \geq 0 \} $ such that
\begin{equation}\label{ouch1}
\left\{ ~ \parbox{11cm}{
$ 1 =: \varepsilon_0^\ddag > \varepsilon_1^\ddag > \cdots >  \varepsilon_m^\ddag \searrow 0 $ \ as $ m \to \infty $,
\\[1ex]
$ \ds \frac{\varepsilon}{2} \bigl| \tilde{\theta}_{m}^\ddag \bigr|_{H^M(\Omega)}^2 < 2^{-m} $ \ for all $ 0 < \varepsilon < \varepsilon_m^\ddag $ and $ m = 1, 2, 3, \ldots \,$.
} \right.
\end{equation}
Considering (\ref{rx50}), (\ref{ouch1}), and Lebesgue's dominated convergence theorem, it can be observed that
\begin{equation*}
\left\{ ~ \parbox{6.4cm}{
$ \alpha(v_{\varepsilon_m^\ddag}^\dag) \to \alpha(v^\dag) $ in $ L^2(\Omega) $
\\[1ex]
$ \beta(v_{\varepsilon_m^\ddag}^\dag) \nabla \tilde{\theta}_m^\ddag \to \beta(v^\dag) \nabla \theta^\ddag $ in $ L^2(\Omega)^N $
} \right. \mbox{ as $ m \to \infty $}
\end{equation*}
and
\begin{eqnarray}
\Phi_\nu(v_{\varepsilon_m^\ddag}^\dag; \tilde{\theta}_m^\ddag) & = & \bigl( \alpha(v_{\varepsilon_m^\ddag}^\dag), |\nabla \tilde{\theta}_m^\ddag| \bigr)_{L^2(\Omega)} +\nu \bigl( \nabla \tilde{\theta}_m^\ddag, \beta(v_{\varepsilon_m^\ddag}^\dag) \nabla \tilde{\theta}_m^\ddag \bigr)_{L^2(\Omega)^N}
\nonumber
\\
& \to & \bigl( \alpha(v^\dag), |\nabla \theta^\ddag| \bigr)_{L^2(\Omega)} +\nu \bigl( \nabla \theta^\ddag, \beta(v^\dag) \nabla \theta^\ddag \bigr)_{L^2(\Omega)^N}
\nonumber
\\
& = & \Phi_\nu(v^\dag; \theta^\ddag), \mbox{ as $ m \to \infty $.}
\label{ouch2}
\end{eqnarray}

Now, based on (\ref{ouch1})--(\ref{ouch2}), the required sequence $ \{ \theta_\varepsilon^\ddag \, | \, 0 < \varepsilon < 1 \} \subset H^M(\Omega) $ is constructed as follows:
\begin{equation*}
\theta_\varepsilon^\ddag := \tilde{\theta}_{m}^\ddag \mbox{ in $ H^M(\Omega) $ \ if $ \varepsilon_{m +1}^\ddag \leq \varepsilon < \varepsilon_m^\ddag $ \ for some $ m \in \Z $ with $ m \geq 0 $.}
\end{equation*}
\qed

\begin{lem}\label{Lem.RX08}
Assume that
\begin{equation}\label{rx52}
\left\{ ~ \parbox{11cm}{
$ \sigma^\dag \in [H^1(\Omega) \cap L^\infty(\Omega)]^2 $,  $ \{ \sigma_m^\dag \, | \, m \in \N \} \subset [H^1(\Omega) \cap L^\infty(\Omega)]^2 $,
\\[1ex]
$ \{ \sigma_m^\dag \, | \, m \in \N \} $ is bounded in $ L^\infty(\Omega)^2 $,
\\[1ex]
$ \sigma_m^\dag \to \sigma^\dag $ in the pointwise sense, a.e.\ in $ \Omega $, as $ m \to \infty $,
} \right.
\end{equation}
and
\begin{equation}\label{rx53}
\left\{ ~ \parbox{11cm}{
$ \omega^\dag \in H^1(\Omega) $,  $ \{ \omega_m^\dag \, | \, m \in \N \} \subset H^1(\Omega) $,
\\[1ex]
$ \omega_m^\dag \to \omega^\dag $ in $ L^2(\Omega) $ and $ \Phi_\nu(\sigma_m^\dag; \omega_m^\dag) \to \Phi_\nu(\sigma^\dag; \omega^\dag) $, as $ m \to \infty $.
} \right.
\end{equation}
Then $ \omega_m^\dag \to \omega^\dag $ in $ H^1(\Omega) $ as $ m \to \infty $.
\end{lem}

\noindent
\textbf{Proof.}
In light of (\hyperlink{A3}{A3}), (\ref{delta_1}) and (\ref{rx53}), we may suppose that
\begin{equation}\label{rx54}
\omega_m^\dag \to \omega^\dag \mbox{ weakly in $ H^1(\Omega) $ as $ m \to \infty $}
\end{equation}
by taking a subsequence if necessary. Here, keeping in mind (\ref{delta_1}), (\ref{rx52}), (\ref{rx54}), and Lebesgue's dominated convergence theorem, we infer
\begin{equation}\label{rx56}
\left\{ ~ \parbox{5.75cm}{
$ \alpha(\sigma_m^\dag) \nabla \omega_m^\dag \to \alpha(\sigma^\dag) \nabla \omega^\dag $ 
\\[1ex]
$ {\textstyle \sqrt{\beta(\sigma_m^\dag)}} \nabla \omega_m^\dag \to {\textstyle \sqrt{\beta(\sigma^\dag)}} \nabla \omega^\dag $
\\[0ex]
$ \frac{1}{\sqrt{\beta(\sigma_m^\dag)}} \nabla \omega_m^\dag \to \frac{1}{\sqrt{\beta(\sigma^\dag)}} \nabla \omega^\dag $
} \right.
\mbox{weakly in $ L^2(\Omega)^N $ as $ m \to \infty $.}
\end{equation}
Additionally, it follows from (\ref{rx53}) that
\begin{eqnarray}
&& \hspace{-7ex} \limsup_{m \to \infty} \left| \rule{-4pt}{11pt} \right.  {\textstyle \sqrt{\beta(\sigma_m^\dag)} \nabla \omega_m^\dag} \left. \rule{-4pt}{11pt} \right|_{L^2(\Omega)^N}^2 
\nonumber
\\
& = & \ds \frac{1}{\nu} \left[ \lim_{m \to \infty} \Phi_\nu(\sigma_m^\dag; \omega_m^\dag) -\liminf_{m \to \infty} \bigl| {\alpha(\sigma_m^\dag)} \nabla \omega_m^\dag \bigr|_{L^1(\Omega; \R^N)} \right]
\label{rx58}
\\
& \leq & \frac{1}{\nu} \left[ \Phi_\nu(\sigma^\dag; \omega^\dag) -\bigl| \alpha(\sigma^\dag) \nabla \omega^\dag \bigr|_{L^1(\Omega; \R^N)} \right] = \bigl| {\textstyle \sqrt{\beta(\sigma^\dag)}} \nabla \omega^\dag \bigr|_{L^2(\Omega)^N}^2.
\nonumber
\end{eqnarray}
By virtue of (\ref{rx56})--(\ref{rx58}) and the uniform convexity of the $ L^2 $-topology, it can be seen that
\begin{equation}\label{rx59-01}
{\textstyle \sqrt{\beta(\sigma_m^\dag)} \nabla \omega_m^\dag} \to {\textstyle \sqrt{\beta(\sigma^\dag)}} \nabla \omega^\dag \mbox{ in $ L^2(\Omega)^N $ \ as $ m \to \infty $,}
\end{equation}
and hence
\begin{eqnarray}
&& \textstyle |\nabla \omega_m^\dag|_{L^2(\Omega)^N}^2 = \left( {\textstyle \sqrt{\beta(\sigma_m^\dag)}} \nabla \omega_m^\dag, \frac{1}{\sqrt{\beta(\sigma_m^\dag)}} \nabla \omega_m^\dag \right)_{L^2(\Omega)^N}
\nonumber
\\
& \to & {\textstyle \left( {\textstyle \sqrt{\beta(\sigma^\dag)}} \nabla \omega^\dag, \frac{1}{\sqrt{\beta(\sigma^\dag)}} \nabla \omega^\dag \right)_{L^2(\Omega)^N}} = |\nabla \omega^\dag|_{L^2(\Omega)^N}^2 \mbox{ as $ m \to \infty $.}
\label{rx59}
\end{eqnarray}
The strong convergence of $ \{ \omega_m^\dag \} $ in $ H^1(\Omega) $ is demonstrated by taking into account (\ref{rx53})--(\ref{rx54}), (\ref{rx59}) and the uniform convexity of the $ L^2 $-topology. 
\hfill \qed

\begin{rem}\label{Rem.RX09}
\begin{em}
Under the same notation as in Lemma \ref{Lem.RX07}, let us assume (\ref{rx50}). Then, for the sequence $ \{ \theta_\varepsilon^\ddag \, | \, 0 < \varepsilon < 1 \} $ as in $ \mbox{(\hyperlink{(m2)_eps^nu}{m2})}_\varepsilon^\nu $, we deduce that
\begin{equation*}
\theta_\varepsilon^\ddag \to \theta^\ddag \mbox{ in $ H^1(\Omega) $ and } \sqrt{\varepsilon} \theta_\varepsilon^\ddag \to 0 \mbox{ in $ H^M(\Omega) $ as $ \varepsilon \searrow 0 $}
\end{equation*}
by applying a similar observation as in Lemma \ref{Lem.RX08}. 
\end{em}
\end{rem}

Finally, we prepare some lemmas for the $ L^\infty $-estimate (\ref{apx03_RC}) as in $ \mbox{(\hyperlink{(AP)_h^nu}{AP})}_h^\nu $.
\begin{lem}[T-monotonicity]\label{Lem.RX09}
Let $ v^\dag \in H^1(\Omega)^2 $ be a fixed function. Then
\begin{equation}\label{T-mono00}
\begin{array}{c}
(\omega_1^* -\omega_2^*, [\omega_1 -\omega_2]^+)_{L^2(\Omega)} \geq 0 
\\[1ex]
\mbox{if $ [\omega_k, \omega_k^*] \in \partial \Phi_\nu(v^\dag;\cdot\,) $ in $ L^2(\Omega)^2 $, $ k = 1, 2 $.}
\end{array}
\end{equation}
\end{lem}

\noindent
\textbf{Proof. }
This lemma can be proved by applying the theory of T-monotonicity (cf.\ \citep{Brezis,KMN}). According to the general theory, we need to start by checking that
\begin{align*}
\Phi_\nu(v^\dag; \tilde{\omega}_1 \wedge \tilde{\omega}_2) &+\Phi_\nu(v^\dag; \tilde{\omega}_1 \vee \tilde{\omega}_2)
\\
 = & \int_\Omega \alpha(v^\dag) |\nabla (\tilde{\omega}_1 \wedge \tilde{\omega}_2)| \, dx +\nu \int_\Omega \beta(v^\dag) |\nabla (\tilde{\omega}_1 \wedge \tilde{\omega}_2)|^2 \, dx
\\
& +\int_\Omega \alpha(v^\dag) |\nabla (\tilde{\omega}_1 \vee \tilde{\omega}_2)| \, dx +\nu \int_\Omega \beta(v^\dag) |\nabla (\tilde{\omega}_1 \vee \tilde{\omega}_2)|^2 \, dx
\\
 = & \sum_{k = 1}^2 \left[ \int_\Omega \alpha(v^\dag) |\nabla \tilde{\omega}_k| \, dx +\nu \int_\Omega \beta(v^\dag) |\nabla \tilde{\omega}_k|^2 \, dx \right]
\\
 = & \Phi_\nu(v^\dag; \tilde{\omega}_1) +\Phi_\nu(v^\dag; \tilde{\omega}_2) \mbox{ for all $ \tilde{\omega}_k \in H^1(\Omega) $, $ k = 1, 2 $.}
\end{align*}
On this basis, the inequality asserted in (\ref{T-mono00}) is verified as follows:
\begin{equation*}
\begin{array}{c}
\ds (\omega_1^* -\omega_2^*, [\omega_1 -\omega_2]^+)_{L^2(\Omega)} 
= (\omega_1^*, \omega_1 -\omega_1 \wedge \omega_2)_{L^2(\Omega)} +(\omega_2^*, \omega_2 -\omega_1 \vee\omega_2)_{L^2(\Omega)}
\\[1ex]
\geq \Phi_\nu(v^\dag; \omega_1) +\Phi_\nu(v^\dag; \omega_2) -\left( \Phi_\nu(v^\dag; \omega_1 \wedge \omega_2) +\Phi_\nu(v^\dag; \omega_1 \vee \omega_2) \right) = 0.
\end{array}
\end{equation*}
\qed

\begin{lem}\label{Lem.RX10}
Let $ v^\dag \in H^1(\Omega) \cap L^\infty(\Omega) $ and $ \check{\theta}_0^\dag, \hat{\theta}_0^\dag \in H^1(\Omega) $ be fixed functions, and let $ [\check{\theta}, \check{\theta}^*], [\hat{\theta}, \hat{\theta}^*] \in L^2(\Omega)^2 $ be pairs of functions such that
\begin{equation}\label{T-mono01}
\left\{ ~ \parbox{12.5cm}{
$ \ds [\check{\theta}, \check{\theta}^*] \in \partial \Phi_\nu(v^\dag;\cdot\,) $ in $ L^2(\Omega)^2 $ and $ \ds  \frac{1}{h} \alpha_0(v^\dag) (\check{\theta} -\check{\theta}_0^\dag) +\check{\theta}^* \leq 0 $ a.e.\ in $ \Omega $,
\\[1ex]
$ [\hat{\theta}, \hat{\theta}^*] \in \partial \Phi_\nu(v^\dag;\cdot\,) $ in $ L^2(\Omega)^2 $ and $ \ds \frac{1}{h} \alpha_0(v^\dag) (\hat{\theta} -\hat{\theta}_0^\dag) +\hat{\theta}^* \geq 0 $ a.e.\ in $ \Omega $,
} \right.
\end{equation}
respectively. Then
\begin{equation*}
|{\textstyle \sqrt{\alpha_0(v^\dag)}}[\check{\theta} -\hat{\theta}]^+|_{L^2(\Omega)}^2 \leq |{\textstyle \sqrt{\alpha_0(v^\dag)}}[\check{\theta}_0^\dag -\hat{\theta}_0^\dag]^+|_{L^2(\Omega)}^2.
\end{equation*}
Moreover, by (\ref{delta_1}), if $ \check{\theta}_0^\dag \leq \hat{\theta}_0^\dag $ a.e.\ in $ \Omega $, then $ \check{\theta} \leq \hat{\theta} $ a.e.\ in $ \Omega $.
\end{lem}

\noindent
\textbf{Proof. }
This lemma is obtained by taking the difference between the inequalities in (\ref{T-mono01}), multiplying both sides of the result by $ [\check{\theta} -\hat{\theta}]^+ $, and applying Lemma \ref{Lem.RX09}. \hfill \qed
\bigskip

\noindent
\textbf{Proof of Theorem \ref{Th.AP} (existence). }
Fix the initial value $ [v_0^\nu, \theta_0^\nu] = [w_0^\nu, \eta_0^\nu, \theta_0^\nu] \in D_1 $ as in (\ref{apxInit}). Then, in light of Lemmas \ref{Lem.RX07}-\ref{Lem.RX08} and Remark \ref{Rem.RX09}, we can take a sequence $ \{ \tilde{\theta}_{\varepsilon, 0}^\nu \, | \, 0 < \varepsilon < 1 \} \subset H^M(\Omega) $ such that
\begin{equation}\label{exist00}
\tilde{\theta}_{\varepsilon, 0}^\nu \to \theta_0^\nu \mbox{ in $ H^1(\Omega) $ and $ \Psi_\varepsilon^\nu(v_0^\nu; \tilde{\theta}_{\varepsilon, 0}^\nu) \to \Phi_\nu(v_0^\nu; \theta_0^\nu) $ as $ \varepsilon \searrow 0 $.}
\end{equation}

Based on this, assume
\begin{equation}\label{h_*^dag}
0 < h < h_*^\dag := h_1^\dag \left( = \frac{1}{4(1 \vee |g|_{C^2([0, 1]^2)})} \right)
\end{equation}
and denote by $ \{ [\tilde{v}_{\varepsilon, i}^\nu, \tilde{\theta}_{\varepsilon, i}^\nu] = [\tilde{w}_{\varepsilon, i}^\nu, \tilde{\eta}_{\varepsilon, i}^\nu, \tilde{\theta}_{\varepsilon, i}^\nu] \, | \, i \in \N \} \subset D_M $ the solution to $ \mbox{(\hyperlink{(RX_eps)_h^nu}{RX$_\varepsilon $})}_h^\nu$, starting from the initial value $ [v_{\varepsilon, 0}^\nu, \theta_{\varepsilon, 0}^\nu] = [v_0^\nu, \tilde{\theta}_{\varepsilon, 0}^\nu] $. Then, from (\ref{exist00}) and Lemma \ref{Lem.RX06}, it can be seen that $ \{ [\tilde{v}_{\varepsilon, i}^\nu, \tilde{\theta}_{\varepsilon, i}^\nu] \, | \, i \in \N \} \subset D_M $ and this sequence is bounded in $ H^1(\Omega)^3 $. By applying Sobolev's embedding theorem, we find a sequence $ \{ \varepsilon_n \, | \, n \in \N \} \subset (0, 1) $ and a sequence of triplets $ \{ [v_{i}^\nu, \theta_{
i}^\nu] = [w_{i}^\nu, \eta_{i}^\nu, \theta_{i}^\nu] \, | \, i \in \N \} \subset H^1(\Omega)^3 $ with $ v_{i}^\nu = [w_{i}^\nu, \eta_{i}^\nu] $, $ i \in \N $ such that \ 
\begin{equation*}
1 > \varepsilon_1 >  \cdots > \varepsilon_n \searrow 0  \mbox{ as $ n \to \infty $,}
\end{equation*}
\begin{equation}\label{exist01}
\left\{ \parbox{11.25cm}{
\begin{tabular}{ll}
\multicolumn{2}{l}{$ \tilde{v}_{n, i}^\nu = [\tilde{w}_{n, i}^\nu, \tilde{\eta}_{n, i}^\nu] := \tilde{v}_{\varepsilon_n, i}^\nu \to v_{i}^\nu $ in $ L^2(\Omega)^2 $, weakly in $ H^1(\Omega)^2 $,}
\\
& weakly-$ * $ in $ L^\infty(\Omega)^2 $, and 
\\ 
& in the pointwise sense a.e.\ in $ \Omega $
\\[0.5ex]
$ \tilde{\theta}_{n, i}^\nu := \tilde{\theta}_{\varepsilon_n, i}^\nu \to \theta_{i}^\nu $ & in $ L^2(\Omega) $, weakly in $ H^1(\Omega) $, and 
\\
& in the pointwise sense a.e.\ in $ \Omega $
\end{tabular}
} \right. \mbox{as $ n \to \infty $,}
\end{equation}
and
\begin{equation}\label{exist02}
v_{i}^\nu = [w_{i}^\nu, \eta_{i}^\nu] \in [o_*, \iota_*] \times [0, 1], \mbox{ a.e.\ in $ \Omega $,}
\end{equation}
for any $ 0 \leq i \in \Z $.

Subsequently, from (\ref{rx03}), (\ref{exist01}), Lemmas \ref{Lem.RX07}--\ref{Lem.RX08} and (\hyperlink{Fact7}{Fact\,7}) in Remark \ref{Rem.MG}, we observe that
\begin{equation}\label{exist03}
\left\{ ~ \parbox{12cm}{
$ \bigl[ \theta_{i}^\nu, -\frac{1}{h}{\textstyle {\alpha_0(v_{i}^\nu)}}(\theta_{i}^\nu -\theta_{i -1}^\nu) \bigr] \in \partial \Phi_\nu(v_{i}^\nu;{}\cdot{}) $ in $ L^2(\Omega)^2 $,
\\[1ex]
$ \Psi_{\varepsilon_n}^{\nu}(\tilde{v}_{n, i}^\nu; \tilde{\theta}_{n, i}^\nu) \to \Phi_\nu(v_{i}^\nu; \theta_{i}^\nu) $ and $ \tilde{\theta}_{n, i}^\nu \to \theta_{i}^\nu $ in $ H^1(\Omega) $ as $ n \to \infty $
} \right.
\end{equation}
for any $ i \in \N $. Furthermore, since
\begin{equation*}
[c, 0] \in \partial \Phi_\nu(v_{i}^\nu;\cdot\,) \mbox{ in $ L^2(\Omega)^2 $ for any constant $ c \in \R $ and any $ 0 \leq i \in \Z $,}
\end{equation*}
it is inductively observed that
\begin{equation}\label{exist05}
\begin{array}{rl}
\theta_{i -1}^\nu \in L^\infty(\Omega) & \mbox{ and \ } \theta_{i}^\nu \leq |\theta_{i -1}^\nu|_{L^\infty(\Omega)} \mbox{ a.e.\ in $ \Omega $}
\\[1ex]
& ( \mbox{resp. } \theta_{i}^\nu \geq -|\theta_{i -1}^\nu|_{L^\infty(\Omega)} \mbox{ a.e.\ in $ \Omega $} ) \mbox{ for any $ i \in \N $,}
\end{array}
\end{equation}
by applying Lemma \ref{Lem.RX10} as the case in which
\begin{equation*}
\left\{ ~ \parbox{12.5cm}{
$ v^\dag = v_i^\nu $, 
\\[0.5ex] 
$ \check{\theta}_0^\dag = \theta_{i -1}^\nu $, $ \hat{\theta}_0^\dag = |\theta_{i -1}^\nu|_{L^\infty(\Omega)} $ (resp. $ \check{\theta}_0^\dag = -|\theta_{i -1}^\nu|_{L^\infty(\Omega)} $, $ \hat{\theta}_0^\dag = \theta_{i -1}^\nu $),
\\[0.5ex]
$ [\check{\theta}, \check{\theta}^*] = \left[ \theta_{i}^\nu, -\frac{1}{h}\alpha_0(v_{i}^\nu) (\theta_{i}^\nu -\theta_{i -1}^\nu) \right] $ (resp. $ [\check{\theta}, \check{\theta}^*] = [ -|\theta_{i -1}^\nu|_{L^\infty(\Omega)}, 0] $),
\\[0.5ex]
$ [\hat{\theta}, \hat{\theta}^*] = [ |\theta_{i -1}^\nu|_{L^\infty(\Omega)}, 0] $ (resp. $ [\hat{\theta}, \hat{\theta}^*] = \left[ \theta_{i}^\nu, -\frac{1}{h}\alpha_0(v_{i}^\nu) (\theta_{i}^\nu -\theta_{i -1}^\nu) \right] $),
} \right. \mbox{for $ i \in \N $.}
\end{equation*}

By invoking (\hyperlink{A1}{A1})--(\hyperlink{A3}{A3}), (\ref{rx01-02}), (\ref{exist01})--(\ref{exist03}), and Lebesgue's dominated convergence theorem and taking further subsequences if necessary, we observe that
\begin{equation}\label{exist04}
\begin{array}{rl}
& \ds \frac{1}{h} (v_{i}^\nu -v_{i -1}^\nu, v_{i}^\nu -\varpi)_{L^2(\Omega)^2} +([\nabla g](v_{i}^\nu), v_{i}^\nu -\varpi)_{L^2(\Omega)^2}
\\[2ex]
& \ds \qquad +\int_\Omega \left( |\nabla \theta_{i -1}^\nu| [\nabla \alpha](v_{i}^\nu) +\nu |\nabla \theta_{i -1}^\nu|^2 [\nabla \beta](v_{i}^\nu) \right) \cdot (v_{i}^\nu -\varpi) \, dx
\\[3ex]
& \ds \qquad +(\nabla v_{i}^\nu, \nabla (v_{i}^\nu -\varpi))_{L^2(\Omega)^{2 \times N}} +\Gamma(v_{i}^\nu)
\\[2ex]
\leq & \ds \liminf_{n \to \infty} \left[ \frac{1}{h} (\tilde{v}_{n, i}^\nu -\tilde{v}_{n, i -1}^\nu, \tilde{v}_{n, i}^\nu -\varpi)_{L^2(\Omega)^2} +([\nabla g](\tilde{v}_{n, i}^\nu), \tilde{v}_{n, i}^\nu -\varpi)_{L^2(\Omega)^2} \right]
\\[2ex]
& \ds \qquad +\lim_{n \to \infty} \int_\Omega \left( |\nabla \tilde{\theta}_{n, i -1}^\nu| [\nabla \alpha](\tilde{v}_{n, i}^\nu) +\nu |\nabla \tilde{\theta}_{n, i -1}^\nu|^2 [\nabla \beta](\tilde{v}_{n, i}^\nu) \right) \cdot (\tilde{v}_{n, i}^\nu -\varpi) \, dx
\\[3ex]
& \ds \qquad + \liminf_{n \to \infty} \left[ (\nabla \tilde{v}_{n, i}^\nu, \nabla (\tilde{v}_{n, i}^\nu -\varpi))_{L^2(\Omega)^{2 \times N}} +\Gamma(\tilde{v}_{n, i}^\nu) \right]
\\[2ex]
\leq & \Gamma(\varpi) \mbox{ \ for any $ \varpi \in [H^1(\Omega) \cap L^\infty(\Omega)]^2 $.}
\end{array}
\end{equation}
With (\ref{exist02})--(\ref{exist04}) in mind, we conclude that the limiting sequence $ \{ [v_{i}^\nu, \theta_{i}^\nu] \, | \, i \in \N \} $ must be the solution to the approximating problem $ \mbox{(\hyperlink{(AP)_h^nu}{AP})}_h^\nu $. \hfill \qed

\subsection*{Uniqueness and energy dissipation of approximating solutions}
\ \vspace{-3ex}

We start with auxiliary lemmas to demonstrate uniqueness. 

\begin{lem}\label{Lem.RX11}
Assume $ 0 < h < h_*^\dag $ with the constant $ h_*^\dag $ as in (\ref{h_*^dag}). Let $ \theta_0^\dag \in H^1(\Omega) $ and $ v_{0, k}^\dag \in [H^1(\Omega) \cap L^\infty(\Omega)]^2 $, $ k = 1, 2 $, be fixed functions, and let $ v_k \in [H^1(\Omega) \cap L^\infty(\Omega)]^2 $, $ k = 1, 2 $, be functions such that
\begin{equation}\label{rx60}
v_k \in [o_*, \iota_*] \times [0, 1] \mbox{ a.e.\ in $ \Omega $}
\end{equation}

and
\begin{equation}\label{rx61}
\begin{array}{rl}
\multicolumn{2}{l}{
\ds \frac{1}{h} (v_k -v_{0, k}^\dag, v_k -\varpi)_{L^2(\Omega)^2} +(\nabla v_k, \nabla (v_k -\varpi))_{L^2(\Omega)^{2 \times N}}}
\\[2ex]
& \ds  \quad +( [\nabla g](v_k), v_k -\varpi)_{L^2(\Omega)^2} +\Gamma(v_k)
\\[1ex]
& \ds  \quad +\int_\Omega \left( |\nabla \theta_0^\dag| [\nabla \alpha](v_k) +\nu |\nabla \theta_0^\dag|^2 [\nabla \beta](v_k)  \right) \cdot (v_k -\varpi) \, dx
\\[2ex]
\leq & \Gamma(\varpi), \mbox{ \ for all $ \varpi \in [H^1(\Omega) \cap L^\infty(\Omega)]^2 $ and $ k = 1, 2 $.}
\end{array}
\end{equation}
Then
\begin{equation}\label{rx62}
|v_1 -v_2|_{L^2(\Omega)^2}^2 \leq 2 |v_{1, 0}^\dag -v_{2, 0}^\dag|_{L^2(\Omega)^2}^2.
\end{equation}
\end{lem}

\noindent
\textbf{Proof.} 
We prepare two inequalities by setting $ k^\perp := (k \mbox{ mod } 2) +1 $ and letting $ \varpi = v_{k^\perp} $ in (\ref{rx61}), for $ k = 1, 2 $. By taking the sum of these, (\ref{rx62}) follows by virtue of (\hyperlink{A2}{A2})--(\hyperlink{A3}{A3}), (\ref{h_*^dag}), (\ref{rx60}), and Young's inequality. \hfill \qed

\begin{lem}\label{Lem.RX12}
Let $ v^\dag \in [H^1(\Omega) \cap L^\infty(\Omega)]^2 $ and $ \theta_{0, k}^\dag \in H^1(\Omega) $, $ k = 1, 2 $, be fixed functions, and let $ \theta_k \in H^1(\Omega) $, $ k = 1, 2 $, be functions such that
\begin{equation*}
\frac{1}{h} \alpha_0(v^\dag) (\theta_k -\theta_{0, k}^\dag) \in \partial \Phi_\nu(v^\dag; \theta_k) \mbox{ in $ L^2(\Omega) $, $ k = 1, 2 $.}
\end{equation*}
Then
\begin{equation*}
|{\textstyle \sqrt{\alpha_0(v^\dag)}}(\theta_1 -\theta_2)|_{L^2(\Omega)}^2 \leq |{\textstyle \sqrt{\alpha_0(v^\dag)}}(\theta_{0, 1} -\theta_{0, 2})|_{L^2(\Omega)}^2.
\end{equation*}
\end{lem}

\noindent
\textbf{Proof. }
This lemma is obtained by applying the standard analytic method for the uniqueness of inclusions governed by subdifferentials (see, e.g., \citep{Brezis,Kenmochi}). \hfill \qed
\bigskip
%

\noindent
\textbf{Proof of Theorem \ref{Th.AP} (uniqueness and energy dissipation).}
Assume $ 0 < h < h_*^\dag $ with the constant as in (\ref{h_*^dag}). Then, on the basis of the foregoing lemmas, the uniqueness for $ \mbox{(\hyperlink{(AP)_h^nu}{AP})}_h^\nu $ is verified through the following steps:
\begin{description}
\item[\textmd{(Step 0)}]Let $ i = 1 $ and fix $ [v_{0}^\nu, \theta_{0}^\nu] = [w_{0}^\nu, \eta_{0}^\nu, \theta_{0}^\nu] \in D_1 $.
\item[\textmd{(Step 1)}]Confirm the uniqueness of the component $ v_{i}^\nu = [w_{i}^\nu, \eta_{i}^\nu] $ of the approximating solution by applying Lemma \ref{Lem.RX11} as the case in which $ \theta_0^\dag = \theta_{i -1}^\nu $ and $ v_{0, 1}^\dag = v_{0, 2}^\dag = v_{i -1}^\nu $.
\item[\textmd{(Step 2)}]Confirm the uniqueness of the component $ \theta_{i}^\nu $ of the approximating solution by keeping in mind (\ref{delta_1}) and applying Lemma \ref{Lem.RX12} as the case in which $ v^\dag = v_{i}^\nu $ and $ \theta_{0, 1}^\dag = \theta_{0, 2}^\dag = \theta_{i -1}^\nu $.
\item[\textmd{(Step\,3)}]Iterate the value of $ i $ , i.e., $ i \leftarrow i +1 $ and return to step 1.
\end{description}

We next verify the inequality (\ref{AP_dissipative}) of energy dissipation. Set $ \varpi = v_{i -1}^\nu $ in (\ref{apx01-02}). Then, by applying a similar derivation method as for (\ref{rx43}) and using (\ref{1/2}), it can be seen that
\begin{equation}\label{apx10}
\begin{array}{rl}
& \ds \frac{1}{2h} |v_{i}^\nu -v_{i -1}^\nu|_{L^2(\Omega)^2}^2 +V_{\rm D}^2(v_{i}^\nu) +\mathscr{G}(v_{i}^\nu) +\Gamma(v_{i}^\nu)
\\[2ex]
& \ds \qquad +\int_\Omega \alpha(v_{i}^\nu)|\nabla \theta_{i -1}^\nu| \, dx +\nu \int_\Omega \beta(v_{i}^\nu) |\nabla \theta_{i -1}^\nu|^2 \, dx
\\[2ex]
\leq & V_{\rm D}^2(v_{i -1}^\nu) +\mathscr{G}(v_{i -1}^\nu) +\Gamma(v_{i -1}^\nu)
\\[2ex]
& \ds \qquad +\int_\Omega \alpha(v_{i -1}^\nu)|\nabla \theta_{i -1}^\nu| \, dx +\nu \int_\Omega \beta(v_{i -1}^\nu) |\nabla \theta_{i -1}^\nu|^2 \, dx, \ i = 1, 2, 3, \ldots.
\end{array}
\end{equation}
Conversely, consider $ \omega = \theta_{i -1}^\nu $ in (\ref{apx03}). Then
\begin{equation}\label{apx11}
\begin{array}{rl}
& \ds \frac{1}{h} |{\textstyle \sqrt{\alpha_0(v_{i}^\nu)}} (\theta_{i}^\nu -\theta_{i -1}^\nu)|_{L^2(\Omega)}^2 
\\[2ex]
& \ds \qquad +\int_\Omega \alpha(v_{i}^\nu)|\nabla \theta_{i}^\nu| \, dx +\nu \int_\Omega \beta(v_{i}^\nu) |\nabla \theta_{i}^\nu|^2 \, dx
\\[2ex]
& \ds \qquad -\int_\Omega \alpha(v_{i}^\nu)|\nabla \theta_{i -1}^\nu| \, dx -\nu \int_\Omega \beta(v_{i}^\nu) |\nabla \theta_{i -1}^\nu|^2 \, dx \leq 0, \mbox{ for $ i= 1, 2, 3, \ldots .$}
\end{array}
\end{equation}
The inequality (\ref{AP_dissipative}) of energy dissipation is obtained by taking the sum of (\ref{apx10}) and (\ref{apx11}). \hfill \qed

\section{Proof of Main Theorem \ref{MainTh01}}
\ \vspace{-3ex}

Throughout this section, we fix the constant $ \nu > 0 $ and assume conditions (\ref{delta_1})--(\ref{D_1}) as in Main Theorem \ref{MainTh01}. We also assume $ 0 < h < h_*^\dag $ with the constant as in (\ref{h_*^dag}), and we denote by $ \{ [v_i^\nu, \theta_i^\nu] = [w_i^\nu, \eta_i^\nu, \theta_i^\nu] \, | \, i \in \N \} \subset D_1 $ the solution to the approximating problem $ \mbox{(\hyperlink{(AP)_h^nu}{AP})}_h^\nu $ under the initial condition
\begin{equation*}
[v_0^\nu, \theta_0^\nu] = [w_0^\nu, \eta_0^\nu, \theta_0^\nu] = [w_0, \eta_0, \theta_0] \in D_1.
\end{equation*}

On this basis, we define three kinds of time interpolation: $ [\overline{v}_h^\nu, \overline{\theta}_h^\nu] = [\overline{w}_h^\nu, \overline{\eta}_h^\nu, \overline{\theta}_h^\nu] \in L_{\rm loc}^{2}([0, \infty); L^2(\Omega)^3) $, $ [\underline{v}_h^\nu, \underline{\theta}_h^\nu] = [\underline{w}_h^\nu, \underline{\eta}_h^\nu, \underline{\theta}_h^\nu] \in  L_{\rm loc}^{2}([0, \infty); L^2(\Omega)^3) $, and $ [\widehat{v}_h^\nu, \widehat{\theta}_h^\nu] = [\widehat{w}_h^\nu, \widehat{\eta}_h^\nu, \widehat{\theta}_h^\nu] $ $\in L_{\rm loc}^{2}([0, \infty); L^2(\Omega)^3) $ with the shorthand $ \overline{v}_h^\nu = [\overline{w}_h^\nu, \overline{\eta}_h^\nu] $, $ \underline{v}_h^\nu = [\underline{w}_h^\nu, \underline{\eta}_h^\nu] $, and $ \widehat{v}_h^\nu = [\widehat{w}_h^\nu, \widehat{\eta}_h^\nu] $, as
\begin{equation}\label{apx200}
\hspace{-3ex} \left\{ \hspace{-2ex} \parbox{14cm}{
\vspace{-2ex}
\begin{itemize}
\item $ [\overline{v}_h^\nu(t), \overline{\theta}_h^\nu(t)] = [\overline{w}_h^\nu(t), \overline{\eta}_h^\nu(t), \overline{\theta}_h^\nu(t)] := [v_i^\nu, \theta_i^\nu] = [w_i^\nu, \eta_i^\nu, \theta_i^\nu] $ in $ L^2(\Omega)^3 $ \linebreak if $ t \in ((i -1)h, ih] $ for some $ i \in \N $, 
\vspace{-1ex}
\item $ [\underline{v}_h^\nu(t), \underline{\theta}_h^\nu(t)] = [\underline{w}_h^\nu(t), \underline{\eta}_h^\nu(t), \underline{\theta}_h^\nu(t)] := [v_{i -1}^\nu, \theta_{i -1}^\nu] = [w_{i -1}^\nu, \eta_{i -1}^\nu, \theta_{i -1}^\nu] $ \linebreak in $ L^2(\Omega)^3 $ if $ t \in [(i -1)h, ih) $ for some $ i \in \N $,
\vspace{-1ex}
\item $ [\widehat{v}_h^\nu(t), \widehat{\theta}_h^\nu(t)] = [\widehat{w}_h^\nu(t), \widehat{\eta}_h^\nu(t), \widehat{\theta}_h^\nu(t)] := [v_i^\nu, \theta_i^\nu] +\left( \frac{t}{h} -i \right) [v_i^\nu -v_{i -1}^\nu, \theta_i^\nu -\theta_{i -1}^\nu] $ in $ L^2(\Omega)^3 $ if $ t \in [(i -1)h, ih) $ for some $ i \in \N $,
\vspace{-2ex}
\end{itemize}
} \right.
\end{equation}
for all $ t \geq 0 $. By using these interpolations, the inequality (\ref{AP_dissipative}) of energy dissipation leads to
\begin{equation}\label{AP_dis01}
\begin{array}{c}
{\ds \frac{1}{2} \int_s^t |(\widehat{v}_h^\nu)_t(\tau)|_{L^2(\Omega)^2}^2 \, d \tau +\int_s^t | {\textstyle \sqrt{\alpha_0(\overline{v}_h^\nu)(\tau)}} (\widehat{\theta}_h^\nu)_t(\tau) |_{L^2(\Omega)}^2 \, d \tau +\mathscr{F}_\nu(\overline{v}_h^\nu(t), \overline{\theta}_h^\nu(t))}
\\[2ex]
\leq \mathscr{F}_\nu(\underline{v}_h^\nu(s), \underline{\theta}_h^\nu(s)) \mbox{ for all $ 0 \leq s \leq t < \infty $ and $ \nu > 0 $}
\end{array}
\end{equation}
and
\begin{equation}\label{AP_dis01_01}
|\mathscr{F}_\nu(\overline{v}_h^\nu(t), \overline{\theta}_h^\nu(t))| \leq F_*^\nu := |\mathscr{F}_\nu(v_0^\nu, \theta_0^\nu)| +c_* \mathscr{L}^N(\Omega) \mbox{ for all $ t > 0 $ and $ \nu > 0 $,}
\end{equation}
where $ c_* > 0 $ is the constant as in (\hyperlink{A2}{A2}). From (\ref{FE00}), (\ref{delta_1}) and Theorem \ref{Th.AP}, it may thus be deduced that
\begin{equation}\label{apx21}
\begin{array}{c}
\left\{ \parbox{11.5cm}{
$ \{ \overline{w}_h^\nu(t, x), \underline{w}_h^\nu(t, x), \widehat{w}_h^\nu(t, x) \, | \, 0 < h < h_*^\dag \} \subset [o_*, \iota_*] $,
\\[1ex]
$ \{ \overline{\eta}_h^\nu(t, x), \underline{\eta}_h^\nu(t, x), \widehat{\eta}_h^\nu(t, x) \, | \, 0 < h < h_*^\dag \} \subset [0, 1] $,
\\[0.5ex]
$ \{ \overline{\theta}_h^\nu(t, x), \underline{\theta}_h^\nu(t, x), \widehat{\theta}_h^\nu(t, x) \, | \, 0 < h < h_*^\dag \} \subset [-|\theta_0^\nu|_{L^\infty(\Omega)}, |\theta_0^\nu|_{L^\infty(\Omega)}] $, 
} \right. 
\\[5ex]
\mbox{ a.e.\ $ x \in \Omega $ and any $ t \in [0, T] $;}
\end{array}
\end{equation}
and
\begin{equation}\label{apx20}
\left\{ \hspace{-3ex} \parbox{9.9cm}{
\vspace{-2ex}
\begin{itemize}
\item $ \{ [\widehat{v}_h^\nu, \widehat{\theta}_h^\nu] = [\widehat{w}_h^\nu, \widehat{\eta}_h^\nu, \widehat{\theta}_h^\nu] \, | \, 0 < h < h_*^\dag \} $ is bounded in $ W^{1, 2}(0, T; L^2(\Omega)^3) \cap L^\infty(0, T; H^1(\Omega)^3) \cap L^\infty(Q)^3 $; 
\vspace{-1ex}
\item $ \{ [\overline{v}_h^\nu, \overline{\theta}_h^\nu] = [\overline{w}_h^\nu, \overline{\eta}_h^\nu, \overline{\theta}_h^\nu] \, | \, 0 < h < h_*^\dag \} $ and $ \{ [\underline{v}_h^\nu, \underline{\theta}_h^\nu] = [\underline{w}_h^\nu, \underline{\eta}_h^\nu, \underline{\theta}_h^\nu] \, | \, 0 < h < h_*^\dag \} $ are bounded in $ L^\infty(0, T; H^1(\Omega)^3) \cap L^\infty(Q)^3 $. 
\vspace{-2ex}
\end{itemize}
} \right.
\end{equation}
Taking into account (\ref{apx21})--(\ref{apx20}) and Aubin-type compactness theory (see \citep{Simon}), we find a sequence
\begin{equation*}
h_*^\dag > h_1^\nu > \cdots > h_n^\nu \searrow 0 \mbox{ as $ n \to \infty $}
\end{equation*}
and a triplet of functions $ [v_\nu, \theta_\nu] = [w_\nu, \eta_\nu, \theta_\nu] \in L^2(0, T; L^2(\Omega)^3) $ such that
\begin{equation}\label{apx25}
\left\{ ~ \parbox{8.5cm}{
$ v_\nu \in W^{1, 2}(0, T; L^2(\Omega)^2) \cap L^\infty(0, T; H^1(\Omega)^2) $,
\\[1ex]
$ \theta_\nu \in W^{1, 2}(0, T; L^2(\Omega)) \cap L^\infty(0, T; H^1(\Omega)) $,
} \right.
\end{equation}
\begin{equation}\label{apx22}
\begin{array}{c}
\left\{ ~ \parbox{8.5cm}{
$ v_\nu(t, x) = [w_\nu(t, x), \eta_\nu(t, x)] \in [o_*, \iota_*] \times [0, 1] $,
\\[1ex]
$ \theta_\nu(t, x) \in [-|\theta_0^\nu|_{L^\infty(\Omega)}, |\theta_0^\nu|_{L^\infty(\Omega)}] $,
} \right. 
\\[3ex]
\mbox{ a.e.\ $ x \in \Omega $ and any $ t \in [0, T] $,}
\end{array}
\end{equation}
and
\begin{equation}\label{apx23}
\begin{array}{c}
\left\{ 
\begin{array}{rl}
\multicolumn{2}{l}{\widehat{v}_n^\nu = [\widehat{w}_n^\nu, \widehat{\eta}_n^\nu] := \widehat{v}_{h_n^\nu}^\nu \to v_\nu \mbox{ \ in $ C(\overline{I}; L^2(\Omega)^2) $,}} 
\\
& \mbox{weakly in $ W^{1, 2}(I; L^2(\Omega)^2) $, weakly-$ * $ in $ L^\infty(I; H^1(\Omega)^2) $,}
\\
& \mbox{and weakly-$ * $ in $ L^\infty(I \times \Omega)^2 $,}
\\[1ex]
\multicolumn{2}{l}{\widehat{\theta}_n^\nu := \widehat{\theta}_{h_n^\nu}^\nu \to \theta_\nu \mbox{ \ in $ C(\overline{I}; L^2(\Omega)) $,}} 
\\
& \mbox{weakly in $ W^{1, 2}(I; L^2(\Omega)) $, weakly-$ * $ in $ L^\infty(I; H^1(\Omega)) $,}
\\
& \mbox{and weakly-$ * $ in $ L^\infty(I \times \Omega) $,}
\end{array}
\right.
\ \\[-2ex]
\ \\
\mbox{as $ n \to \infty $, \  for any open interval $ I \subset (0, T) $.} 
\end{array}
\end{equation}

Additionally, noting that
\begin{equation}\label{apx300}
\left\{ \begin{array}{ll}
\multicolumn{2}{l}{\ds |\overline{v}_h^\nu -\widehat{v}_h^\nu|_{L^\infty(0, T; L^2(\Omega)^2)} \vee |\underline{v}_h^\nu -\widehat{v}_h^\nu|_{L^\infty(0, T; L^2(\Omega)^2)}}
\\[1ex]
& \ds \leq \max_{i \in \Z, \, ih \in [0, T]} \int_{ih}^{(i +1)h} \hspace{-4ex} |(\widehat{v}_h^\nu)_t(t)|_{L^2(\Omega)^2} \, dt \leq \sqrt{h} \, |(\widehat{v}_h^\nu)_t|_{L^2(0, T; L^2(\Omega)^2)},
\\[2ex]
\multicolumn{2}{l}{\ds |{\textstyle \sqrt{\alpha_0(\overline{v}_h^\nu)}}(\overline{\theta}_h^\nu -\widehat{\theta}_h^\nu)|_{L^\infty(0, T; L^2(\Omega))} \vee |{\textstyle \sqrt{\alpha_0(\overline{v}_h^\nu)}}(\underline{\theta}_h^\nu -\widehat{\theta}_h^\nu)|_{L^\infty(0, T; L^2(\Omega))}}
\\[1ex]
& \ds \leq \max_{i \in \Z, \, ih \in [0, T]} \int_{ih}^{(i +1)h} \hspace{-4ex} |{\textstyle \sqrt{\alpha_0(\overline{v}_h^\nu)}}(\widehat{\theta}_h^\nu)_t(t)|_{L^2(\Omega)} \, dt \leq \sqrt{h} \, |{\textstyle \sqrt{\alpha_0(\overline{v}_h^\nu)}}(\widehat{\theta}_h^\nu)_t|_{L^2(0, T; L^2(\Omega))},
\end{array} \right.
\end{equation}
we also have
\begin{equation}\label{apx24}
\begin{array}{c}
\left\{ 
\begin{array}{rl}
\multicolumn{2}{l}{\overline{v}_n^\nu = [\overline{w}_n^\nu, \overline{\eta}_n^\nu] := \overline{v}_{h_n^\nu}^\nu \to v_\nu \mbox{ and }\underline{v}_n^\nu = [\underline{w}_n^\nu, \underline{\eta}_n^\nu] := \underline{v}_{h_n^\nu}^\nu \to v_\nu} 
\\
& \mbox{in $ L^\infty(I; L^2(\Omega)^2) $, weakly-$ * $ in $ L^\infty(I; H^1(\Omega)^2) $,}
\\
& \mbox{and weakly-$ * $ in $ L^\infty(I \times \Omega)^2 $,}
\\[1ex]
\multicolumn{2}{l}{\overline{\theta}_n^\nu  := \overline{\theta}_{h_n^\nu}^\nu \to \theta_\nu \mbox{ and }\underline{\theta}_n^\nu := \underline{\theta}_{h_n^\nu}^\nu \to \theta_\nu \mbox{ \ in $ L^\infty(I; L^2(\Omega)) $,}} 
\\
& \mbox{weakly-$ * $ in $ L^\infty(I; H^1(\Omega)) $ and weakly-$ * $ in $ L^\infty(I \times \Omega) $,}
\end{array}
\right. 
\ \\[-2ex]
\ \\
\mbox{as $ n \to \infty $, \  for any open interval $ I \subset (0, T) $,} 
\end{array}
\end{equation}
and, in particular,
\begin{equation}\label{*end}
\begin{array}{c}
\left\{ \hspace{0ex} \parbox{8.75cm}{
\vspace{-2ex}
\begin{itemize}
\item[]\hspace{-4ex}$ \overline{v}_n^\nu(t) \to v_\nu(t) $, $ \underline{v}_n^\nu(t) \to v_\nu(t) $ and $ \widehat{v}_n^\nu(t) \to v_\nu(t) $ in $ L^2(\Omega)^2 $ and weakly in $ H^1(\Omega) $,
\item[]\hspace{-4ex}$ \overline{\theta}_n^\nu(t) \to \theta_\nu(t) $, $ \underline{\theta}_n^\nu(t) \to \theta_\nu(t) $ and $ \widehat{\theta}_n^\nu(t) \to \theta_\nu(t) $ in $ L^2(\Omega) $ and weakly in $ H^1(\Omega) $,
\vspace{-2ex}
\end{itemize}
} \right. 
\ \\[-2ex]
\ \\
\mbox{as $ n \to \infty $, \  a.e.\ $ t \in (0, T) $.} 
\end{array}
\end{equation}

Based on these considerations, we next demonstrate the following auxiliary lemmas.

\begin{lem}[Mosco convergence]\label{Lem.Mosco01}
Let $ I \subset (0, T) $ be any open interval, and let $ {\Phi}_\nu^I : L^2(I; L^2(\Omega)) \rightarrow [0, \infty] $ and $ {\Phi}_{\nu, n}^I : L^2(I; L^2(\Omega)) \rightarrow [0, \infty] $, $ n \in \N $, be functionals defined as
\begin{equation}\label{apx100}
\zeta \in L^2(I; L^2(\Omega)) \mapsto {\Phi}_\nu^I(\zeta) := \int_I \Phi_\nu(v_\nu(t); \zeta(t)) \, dt,
\end{equation}
and
\begin{equation}\label{apx101}
\zeta \in L^2(I; L^2(\Omega)) \mapsto {\Phi}_{\nu, n}^I(\zeta) := \int_I \Phi_\nu(\overline{v}_n^\nu(t); \zeta(t)) \, dt, ~ n \in \N,
\end{equation}
by using the functions $ v_\nu = [w_\nu, \eta_\nu] \in L^\infty(0, T; H^1(\Omega)^2) \cap L^\infty(Q)^2 $ and $ \overline{v}_n^\nu = [\overline{w}_n^\nu, \overline{\eta}_n^\nu] \in L^\infty(0, T; H^1(\Omega)^2) \cap L^\infty(Q)^2 $, $ n \in \N $, as in (\ref{apx25})--(\ref{*end}). Then, the following two statements hold:
\begin{description}
\item[\textmd{\it (\hypertarget{I-1}{I-1})}]$ {\Phi}_\nu^I $ and $ {\Phi}_{\nu, n} $, $ n \in \N $, are proper l.s.c and convex functions on $ L^2(I; L^2(\Omega)) $ such that $ D({\Phi}_\nu^I) = D({\Phi}_{\nu, n}^I) = L^2(I; H^1(\Omega)) $ for all $ n \in \N $. 
\item[\textmd{\it (\hypertarget{I-2}{I-2})}]$ {\Phi}_{\nu, n}^I \to {\Phi}_\nu^I $ on $ L^2(I; L^2(\Omega)) $, in the sense of Mosco \citep{Mosco}, as $ n \to \infty $. 
\end{description}
\end{lem}

\noindent
\textbf{Proof. }
Since (\hyperlink{I-1}{I-1}) is a straightforward consequence of Notation \ref{Note08}, (\hyperlink{A3}{A3}), (\ref{delta_1}), and (\ref{apx20})--(\ref{apx22}), it is enough to prove only (\hyperlink{I-2}{I-2}). 
\medskip

To verify the \hyperlink{m1}{lower bound} condition, let us take a function $ \zeta^\dag \in L^2(I; L^2(\Omega)) $, a sequence $ \{ \zeta_n^\dag \, | \, n \in \N \} \subset L^2(I; L^2(\Omega)) $, and a subsequence $ \{ \zeta_{n_j}^\dag \, | \, j \in \N \} \subset \{ \zeta_n^\dag \}
$ to posit the following nontrivial situation:
\begin{equation}\label{apx102}
\left\{ ~ \parbox{8cm}{
$ \ds \zeta_n^\dag \to \zeta^\dag $ weakly in $ L^2(I; L^2(\Omega)) $ as $ n \to \infty $, 
\\[1ex]
$ \ds \liminf_{n \to \infty} {\Phi}_{\nu, n}^I(\zeta_n^\dag) = \lim_{j \to \infty} {\Phi}_{\nu, n_j}^I(\zeta_{n_j}^\dag) < \infty $.
} \right.
\end{equation}
By virtue of (\ref{delta_1}) and (\ref{apx102}), the subsequence $ \{ \zeta_{n_j}^\dag \} $ must be bounded in $ L^2(I; H^1(\Omega)) $. So, taking a subsequence if necessary, we may also suppose that
\begin{equation*}
\zeta_{n_j}^\dag \to \zeta^\dag \mbox{ weakly in $ L^2(I; H^1(\Omega)) $ as $ j \to \infty $.}
\end{equation*}
Furthermore, with (\hyperlink{A3}{A3}) and (\ref{apx21})--(\ref{*end}) in mind, we have
\begin{equation*}
\left\{ \begin{array}{l}
\alpha(\overline{v}_{n_j}^\nu) \nabla \zeta_{n_j}^\dag \to \alpha(v_\nu) \nabla \zeta^\dag
\\[2ex]
{\textstyle \sqrt{\beta(\overline{v}_{n_j}^\nu)}} \nabla \zeta_{n_j}^\dag \to {\textstyle \sqrt{\beta(v_\nu)}} \nabla \zeta^\dag
\end{array} \right. \mbox{weakly in $ L^2(I; L^2(\Omega)^N) $ as $ j \to \infty $.}
\end{equation*}
From the above convergence, satisfaction of the \hyperlink{m1}{lower bound} condition is confirmed as follows:
\begin{equation*}
\begin{array}{rl}
\multicolumn{2}{l}{\ds \liminf_{n \to \infty} {\Phi}_{\nu, n}^I(\zeta_n^\dag) = \lim_{j \to \infty} {\Phi}_{\nu, n_j}^I(\zeta_{n_j}^\dag)}
\\[2ex]
\qquad \quad  & \ds = \liminf_{j \to \infty} \left( \left| \alpha(\overline{v}_{n_j}^\nu) \nabla \zeta_{n_j}^\dag \right|_{L^1(I; L^1(\Omega; \R^N))} +\nu \left| {\textstyle \sqrt{\beta(\overline{v}_{n_j}^\nu)}} \nabla \zeta_{n_j}^\dag \right|_{L^2(I; L^2(\Omega)^N)}^2 \right)
\\[2ex]
\qquad \quad  & \ds \geq \left| \alpha(v_\nu) \nabla \zeta^\dag \right|_{L^1(I; L^1(\Omega; \R^N))} +\nu \left| {\textstyle \sqrt{\beta(v_\nu)}} \nabla \zeta^\dag \right|_{L^2(I; L^2(\Omega)^N)}^2 = {\Phi}_\nu^I(\zeta^\dag).
\end{array}
\end{equation*}

Additionally, owing to (\ref{apx21})--(\ref{*end}) and Lebesgue's dominated convergence theorem, it can be inferred that
\begin{equation*}
\begin{array}{rl}
\multicolumn{2}{l}{|{\Phi}_{\nu, n}^I(\zeta^\ddag) -{\Phi}_\nu^I(\zeta^\ddag)|}
\\[1ex]
\qquad & \ds \leq \int_I \int_\Omega |\alpha(\overline{\eta}_n^\nu) -\alpha(\eta_\nu)||\nabla \zeta^\ddag| \, dx \, dt +\nu \int_I \int_\Omega |\beta(\overline{w}_n^\nu) -\beta(w_\nu)||\nabla \zeta^\ddag|^2 \, dx \, dt
\\[2ex]
\qquad & \to 0 \mbox{ as $ n \to \infty $, \ for any $ \zeta^\ddag \in  D({\Phi}_\nu^I) = L^2(I; H^1(\Omega)) $.}
\end{array}
\end{equation*}
This implies the validity of the condition of \hyperlink{m2}{optimality}, for Mosco convergence $ {\Phi}_{\nu, n}^I \to {\Phi}_\nu^I $ on $ L^2(I; L^2(\Omega)) $ as $ n \to \infty $. \hfill \qed

\begin{lem}\label{Lem.L^2(H^1)-conv}
In addition to the assumptions and notations of Lemma \ref{Lem.Mosco01}, assume that $ \zeta^\ddag \in L^2(I; H^1(\Omega)) $, $ \{ \zeta_n^\ddag \, | \, n \in \N \} \subset L^2(I; H^1(\Omega)) $ and
\begin{equation}\label{MT1-00}
{\Phi}_{\nu, n}^I(\zeta_n^\ddag) \to {\Phi}_\nu^I(\zeta^\ddag) \mbox{ as $ n \to \infty $.}
\end{equation}
Then $ \zeta_n^\ddag \to \zeta^\ddag $ in $ L^2(I; H^1(\Omega)) $ as $ n \to \infty $. 
\end{lem}

\noindent
\textbf{Proof. }
The proof of this lemma is a slight modification of that of Lemma \ref{Lem.RX08}. Indeed, from (\hyperlink{A3}{A3}), (\ref{delta_1}), (\ref{apx21})--(\ref{*end}) and (\ref{MT1-00}), we infer that
\begin{equation}\label{*MT1-00}
\left\{ ~ \parbox{5.5cm}{
$ \alpha(\overline{v}_n^\nu) \nabla \zeta_n^\ddag \to \alpha({v}_\nu) \nabla \zeta^\ddag $ 
\\[1.5ex]
$ {\textstyle \sqrt{\beta(\overline{v}_n^\nu)}} \nabla \zeta_n^\ddag \to {\textstyle \sqrt{\beta(v_\nu)}} \nabla \zeta^\ddag $
\\[1ex]
$ \frac{1}{\sqrt{\beta(\overline{v}_n^\nu)}} \nabla \zeta_n^\ddag \to \frac{1}{\sqrt{\beta(v_\nu)}} \nabla \zeta^\ddag $
} \right.
\mbox{weakly in $ L^2(I; L^2(\Omega)^N) $ as $ n \to \infty $.}
\end{equation}
Therefore, as in the derivations of (\ref{rx59-01})--(\ref{rx59}), convergences (\ref{MT1-00})--(\ref{*MT1-00}) show that
\begin{equation}\label{MT1-01}
{\textstyle \sqrt{\beta(\overline{v}_n^\nu)}} \nabla \zeta_n^\ddag \to \sqrt{\beta(v_\nu)} \nabla \zeta^\ddag \mbox{ in $ L^2(I; L^2(\Omega)^N) $ as $ n \to \infty $,}
\end{equation}
and
\begin{equation}\label{MT1-05}
|\nabla \zeta_n^\ddag|_{L^2(I; L^2(\Omega)^N)}^2 \to |\nabla \zeta^\ddag|_{L^2(I; L^2(\Omega)^N)}^2 \mbox{ as $ n \to \infty $. }
\end{equation}
Thus, strong convergence of $ \{ \zeta_n^\ddag \} $ in $ L^2(I; H^1(\Omega)) $ follows from (\ref{MT1-01}), (\ref{MT1-05}), and the uniform convexity of the $ L^2 $-topology.  
\hfill \qed

\bigskip

\noindent
\textbf{Proof of Main Theorem \ref{MainTh01}. }
From (\ref{apx25})--(\ref{apx22}), we see that the limiting triplet $ [v_\nu, \theta_\nu] = [w_\nu, \eta_\nu, \theta_\nu] $ fulfills the condition $ \mbox{(\hyperlink{(S0)_nu}{S0})}_\nu $. Hence, all we have to do is verify the compatibility of $ [v_\nu, \theta_\nu] $ with conditions $ \mbox{(\hyperlink{(S1)_nu}{S1})}_\nu $--$ \mbox{(\hyperlink{(S3)_nu}{S3})}_\nu $.
\medskip

Fix any open interval $ I \subset (0, T) $. Then, from Remark \ref{Rem.Time-dep.}, (\ref{apx01-02})--(\ref{apx03}), and (\ref{apx200}), the functions $ [\overline{v}_n^\nu, \overline{\theta}_n^\nu] $, $ [\underline{v}_n^\nu, \underline{\theta}_n^\nu] $, $ [\widehat{v}_n^\nu, \widehat{\theta}_n^\nu] $, $ n \in \N $, must satisfy
\begin{equation}\label{MT1-10}
\begin{array}{rl}
\multicolumn{2}{l}{\ds \int_I \left( \rule{0pt}{10pt} (\widehat{v}_n^\nu)_t(t), \overline{v}_n^\nu(t) -\varpi \right)_{L^2(\Omega)^2} \, dt  +\int_I \left( \rule{0pt}{10pt} [\nabla g](\overline{v}_n^\nu(t)), \overline{v}_n^\nu -\varpi \right)_{L^2(\Omega)^2} \, dt}
\\[2ex]
\qquad \quad & \ds +\int_I \left( \rule{0pt}{10pt} \nabla \overline{v}_n^\nu(t), \nabla (\overline{v}_n^\nu(t) -\varpi) \right)_{L^2(\Omega)^{2 \times N}} \, dt
\\[2ex]
\qquad \quad & \ds +\int_I \left( \rule{0pt}{10pt} |\nabla \underline{\theta}_n^\nu(t)| [\nabla \alpha](\overline{v}_n^\nu(t)) +\nu |\nabla \underline{\theta}_n^\nu|^2 [\nabla \beta](\overline{v}_n^\nu(t)) \right) \cdot (\overline{v}_n^\nu(t) -\varpi) \, dt
\\[2ex]
& \ds  +\int_I \Gamma(\overline{v}_n^\nu(t)) \, dt \leq \ds \int_I \Gamma(\varpi) \, dt \left( = \Gamma(\varpi) \cdot \mathscr{L}^1(I) \right) 
\\[2ex]
\multicolumn{2}{l}{\mbox{ for any $ \varpi \in [H^1(\Omega) \cap L^\infty(\Omega)]^2 $ and any $ n \in \N $,}}
\end{array}
\end{equation}
and 
\begin{equation}\label{MT1-11}
[\overline{\theta}_n^\nu, -\alpha_0(\overline{v}_n^\nu) (\widehat{\theta}_n^\nu)_t] \in \partial {\Phi}_{\nu, n}^I \mbox{ in $ L^2(I; L^2(\Omega))^2 $ for any $ n \in \N $.}
\end{equation}
Here, from (\ref{apx23})--(\ref{apx101}), (\ref{MT1-11}), (\hyperlink{I-2}{I-2}) of Lemma \ref{Lem.Mosco01}, and (\hyperlink{Fact7}{Fact\,7}) in Remark \ref{Rem.MG}, it follows that
\begin{equation}\label{MT1-12}
[\theta_\nu, -\alpha_0(v_\nu) (\theta_\nu)_t] \in \partial {\Phi}_\nu^I \mbox{ in $ L^2(I; L^2(\Omega))^2 $}
\end{equation}
and
\begin{equation}\label{MT1-13}
{\Phi}_{\nu, n}^I(\overline{\theta}_n^\nu) \to {\Phi}_\nu^I(\theta_\nu) \mbox{ as $ n \to \infty $.}
\end{equation}
Subsequently, in light of (\ref{MT1-12}), (\hyperlink{I-1}{I-1}) of Lemma \ref{Lem.Mosco01}, and (\hyperlink{Fact1}{Fact 1}) in Remark \ref{Rem.Time-dep.}, we can deduce that the triplet $ [v_\nu, \theta_\nu] = [w_\nu, \eta_\nu, \theta_\nu] $ fulfills condition $ (\mbox{\hyperlink{(S3)_nu}{S3}})_\nu $. 
\bigskip

Next, from (\ref{apx23})--(\ref{apx24}), (\ref{MT1-13}), and Lemma \ref{Lem.L^2(H^1)-conv}, it follows that
\begin{equation}\label{MT1-14}
\overline{\theta}_n^\nu \to \theta_\nu \mbox{ in $ L^2(I; H^1(\Omega)) $ as $ n \to \infty $.}
\end{equation}
Meanwhile, from (\ref{delta_1}), (\ref{AP_dis01})--(\ref{AP_dis01_01}), and (\ref{apx20}), it can be seen that
\begin{equation}\label{MT2-07}
\begin{array}{c}
\ds \left| \int_I \int_\Omega |\nabla \overline{\theta}_h^\nu|^2 \, dx \, dt -\int_I \int_\Omega |\nabla \underline{\theta}_h^\nu|^2 \, dx \, dt \right| \leq \frac{2 F_*^\nu}{\delta_1} \cdot \frac{h}{\nu} 
\\[2ex]
\mbox{ \ for all $ 0 < h < h_*^\dag $ and $ \nu > 0 $.}
\end{array}
\end{equation}
Bearing in mind (\ref{apx24}), (\ref{MT1-14}), and (\ref{MT2-07}), it further follows that
\begin{equation*}
\underline{\theta}_n^\nu \to \theta_\nu \mbox{ and  } \widehat{\theta}_n^\nu \to \theta_\nu \mbox{ in $ L^2(I; H^1(\Omega)) $, as $ n \to \infty $.}
\end{equation*}

Now, on account of (\ref{apx23})--(\ref{apx24}), (\ref{MT1-14})--(\ref{MT2-07}), and Lebesgue's dominated convergence theorem, allowing $ n \to \infty $ in (\ref{MT1-10}) yields
\begin{equation*}
\begin{array}{rl}
\multicolumn{2}{l}{\ds \int_I \left( \rule{0pt}{10pt} (v_\nu)_t(t), v_\nu(t) -\varpi \right)_{L^2(\Omega)^2} \, dt  +\int_I \left( \rule{0pt}{10pt} [\nabla g](v_\nu(t)), v_\nu(t) -\varpi \right)_{L^2(\Omega)^2} \, dt}
\\[2ex]
\qquad \quad & \ds +\int_I \left( \rule{0pt}{10pt} \nabla v_\nu(t), \nabla (v_\nu(t) -\varpi) \right)_{L^2(\Omega)^{2 \times N}} \, dt
\\[2ex]
\qquad \quad & \ds +\int_I \left( \rule{0pt}{10pt} |\nabla \theta_\nu(t)| [\nabla \alpha](v_\nu(t)) +\nu |\nabla \theta_\nu(t)|^2 [\nabla \beta](v_\nu(t)) \right) \cdot (v_\nu(t) -\varpi) \, dt
\\[2ex]
& \ds  +\int_I \Gamma(v_\nu(t)) \, dt \leq \ds \int_I \Gamma(\varpi) \, dt \left( = \Gamma(\varpi) \cdot \mathscr{L}^1(I) \right) 
\\[2ex]
\multicolumn{2}{l}{\mbox{ for any $ \varpi \in [H^1(\Omega) \cap L^\infty(\Omega)]^2 $.}}
\end{array}
\end{equation*}
Since the choice of the open interval $ I \subset (0, T) $ is arbitrary, we can verify the remaining conditions $ \mbox{(\hyperlink{(S1)_nu}{S1})}_\nu $--$ \mbox{(\hyperlink{(S2)_nu}{S2})}_\nu $ on the basis of this inequality and the reformulation (\ref{unify01}) in Remark \ref{Rem.1.0}. \hfill \qed

\section{Proof of Main Theorem \ref{MainTh02}}
\ \vspace{-3ex}

Assume (\ref{D_0}) for the initial value $ [w_0, \eta_0, \theta_0] \in D_0 $ of the system $ \mbox{(S)}_0 $. Then, roughly speaking, our second main theorem is proved through some limiting observations for $ \mbox{(\hyperlink{(AP)_h^nu}{AP})}_h^\nu $ as $ h, \nu \searrow 0 $. With this scenario in mind, we rely upon \citep[Section 2]{MS} for the following auxiliary lemmas. 

\begin{lem}\label{MS_Lem.7}
Let $ \delta_* > 0 $ be a constant and let $ I \subset (0, T) $ be any open interval. Assume that
\begin{equation*}
\left\{ \hspace{-5ex} \parbox{15cm}{
\vspace{-2ex}
\begin{itemize}
\item[]$ \varrho \in C(\overline{I}; L^2(\Omega)) \cap L^\infty(I; H^1(\Omega)) \cap L^\infty(I \times \Omega) $, $ \{ \varrho_n \, | \, n \in \N \} \subset L^2(I; L^2(\Omega)) $,
\vspace{-1ex}
\item[]$ \varrho \geq \delta_* $ and $ \varrho_n \geq \delta_* $, a.e.\ in $ I \times \Omega $, for all $ n \in \N $,
\vspace{-1ex}
\item[]$ \varrho_n(t) \to \varrho(t) $ in $ L^2(\Omega) $ and weakly in $ H^1(\Omega) $, as $ n \to \infty $, a.e.\ $ t \in I $,
\vspace{-2ex}
\end{itemize}
} \right.
\end{equation*}
and
\begin{equation*}
\left\{ \hspace{-5ex} \parbox{14cm}{
\vspace{-2ex}
\begin{itemize}
\item[]$ \zeta \in C(\overline{I}; L^2(\Omega)) \cap L^1(I; BV(\Omega)) $, $ \{ \zeta_n \, | \, n \in \N \} \subset L^2(I; H^1(\Omega)) $,
\vspace{-1ex}
\item[]$ \zeta_n(t) \to \zeta(t) $ in $ L^2(\Omega) $ as $ n \to \infty $, \ a.e.\ $ t \in I $.
\vspace{-2ex}
\end{itemize}
} \right.
\end{equation*}
Then the functions
\begin{equation*}
t \in I \mapsto \int_\Omega d[\varrho(t)|D\zeta(t)|] \mbox{ and } 
t \in I \mapsto \int_\Omega \varrho_n(t)|\nabla \zeta_n(t)| \, dx, ~ n \in \N,
\end{equation*}
are integrable. Moreover, if
\begin{equation*}
\int_I \int_\Omega \varrho_n(t) |\nabla \zeta_n(t)| \, dx \, dt \to \int_I \int_\Omega d[\varrho(t) |D \zeta(t)|] \mbox{ as $ n \to \infty $} 
\end{equation*}
and
\begin{equation*}
\left\{ \hspace{-5ex} \parbox{15cm}{
\vspace{-2ex}
\begin{itemize}
\item[]$ \omega \in C(\overline{I}; L^2(\Omega)) \cap L^\infty(I; H^1(\Omega)) \cap L^\infty(I \times \Omega) $ and $ \{ \omega_n \, | \, n \in \N \} \subset L^2(I; L^2(\Omega)) $,
\vspace{-4ex}
\item[]$ \{ \omega_n \, | \, n \in \N \} $ is a bounded sequence in $ L^\infty(I \times \Omega) $,
\vspace{-1ex}
\item[]$ \omega_n(t) \to \omega(t) $ in $ L^2(\Omega) $ and weakly in $ H^1(\Omega) $ as $ n \to \infty $, a.e.\ $ t \in I $,
\vspace{-2ex}
\end{itemize}
} \right.
\end{equation*}
then
\begin{equation*}
\int_I \int_\Omega \omega_n(t) |\nabla \zeta_n(t)| \, dx \, dt \to \int_I \int_\Omega d [\omega(t) |D \zeta(t)|] \mbox{ as $ n \to \infty $.}
\end{equation*}
\end{lem}

\noindent
\textbf{Proof.}
This lemma is a straightforward consequence of \citep[Lemmas 4 and 7]{MS}. \hfill \qed

\begin{lem}[{\boldmath $ \Gamma $-convergence}]\label{Lem.GammaLim01}
Assume $ v^\ddag \in [H^1(\Omega) \cap L^\infty(\Omega)]^2 $, $ \{ v_\nu^\ddag \, | \, \nu > 0 \} \subset [H^1(\Omega) \cap L^\infty(\Omega)]^2 $, and
\begin{equation}\label{GammaLim01}
\left\{ ~ \parbox{10cm}{
$ v_\nu^\ddag \to v^\ddag $ in the pointwise sense, a.e.\ in $ \Omega $, as $ \nu \searrow 0 $,
\\[1ex]
$ \{ v_\nu^\ddag \, | \, \nu > 0 \} $ is bounded in $ L^\infty(\Omega)^2 $.  
} \right.
\end{equation}
Then, for the sequence of convex functions $ \{ \Phi_\nu(v_\nu^\ddag;\cdot\,) \} $, it holds that $ \Phi_\nu(v_\nu^\ddag;\cdot\,) \to  \Phi_0(v^\ddag;\cdot\,) $ on $ L^2(\Omega) $, in the sense of $ {\mit \Gamma} $-convergence \citep{GammaConv}, as $ \nu \searrow 0 $, i.e.:
\begin{description}
\item[$ \mbox{\textmd{\em (\hypertarget{(gamma1)_nu}{$\gamma $1})}}_\nu $](lower bound) $ \ds \liminf_{\nu \searrow 0} \Phi_\nu(v_\nu^\ddag; \theta_\nu^{\ddag}) \geq \Phi_\nu(v^\ddag; \theta^{\ddag}) $ if $ \theta^{\ddag} \in L^2(\Omega) $, $ \{ \theta_\nu^{\ddag} \, | \, \nu > 0 \} \subset L^2(\Omega) $, and $ \theta_\nu^{\ddag} \to \theta^{\ddag} $ in $ L^2(\Omega) $ as $ \nu \searrow 0 $; 
\item[$ \mbox{\textmd{\em (\hypertarget{(gamma2)_nu}{$\gamma $2})}}_\nu $](optimality) for any $ \theta^{\ddag\ddag} \in H^1(\Omega) $, there exists a sequence $ \{ \theta_\nu^{\ddag\ddag} \, | \, \nu > 0 \} \subset H^1(\Omega) $ such that $ \theta_\nu^{\ddag\ddag} \to \theta^{\ddag\ddag} $ in $ L^2(\Omega) $ and $ \ds \Phi_\nu(v_\nu^\ddag; \theta_\nu^{\ddag\ddag}) \to \Phi_0(v^\ddag; \theta^{\ddag\ddag}) $ as $ \nu \searrow 0 $. 
\end{description}
\end{lem}

\noindent
\textbf{Proof.}
We start by confirming that
\begin{equation}\label{GammaLim02}
\left\{ ~ \parbox{11cm}{
$ \ds R_1 := \sup_{\nu > 0} |v_\nu^\ddag|_{L^\infty(\Omega)^2} < \infty $,
\\[1ex]
$ \alpha(v_\nu^\ddag) \to \alpha(v^\ddag) $ in $ L^2(\Omega) $ and weakly in $ H^1(\Omega) $, as $ \nu \searrow 0 $.
} \right.
\end{equation}
This is easily confirmed from assumption (\hyperlink{A3}{A3}) and (\ref{GammaLim01}). 

Next, for any $ \nu > 0 $, any $ \delta > 0 $, and any $ \tilde{v} \in [H^1(\Omega) \cap L^\infty(\Omega)]^2 $, let $ \tilde{\Phi}_{\nu, \delta}(\tilde{v}; \cdot \, ) $ be a proper l.s.c.\ and convex function on $ L^2(\Omega) $, defined as
\begin{equation}\label{PhiTilde}
\vartheta \in L^2(\Omega) \mapsto \tilde{\Phi}_{\nu, \delta}(\tilde{v}; \vartheta) := \Phi_0(\tilde{v}; \vartheta) +\nu \delta V_{\rm D}^1(\vartheta) \in [0, \infty],
\end{equation}
where $ V_{\rm D}^1 $ is the convex function $ V_{\rm D}^d $ given in (\ref{Dirichlet^d}), in the case that $ d = 1 $. Then, setting
\begin{equation}\label{delta_beta(R)}
\delta_*(R) := |\alpha|_{C([-R, R]^2)} +|\beta|_{C([-R, R]^2)} \mbox{ for any $ R > 0 $,}
\end{equation}
it immediately follows from (\hyperlink{A3}{A3}), (\ref{delta_1}), and (\ref{GammaLim02}) that
\begin{equation}\label{GammaLim04}
\tilde{\Phi}_{\nu, \delta_1}(v_\nu^\ddag; z) \leq \Phi_\nu(v_\nu^\ddag; z) \leq \tilde{\Phi}_{\nu, \delta_*({R_1})}(v_\nu^\ddag; z) \mbox{ for any $ \nu > 0 $ and any $ z \in L^2(\Omega) $.}
\end{equation}
In addition, we can apply \citep[Lemma 3]{MS} to see that
\begin{equation}\label{GammaLim03}
\begin{array}{c}
\tilde{\Phi}_{\nu, \delta_1}(v_\nu^\ddag;\cdot\,) \to \Phi_0(v^\ddag;\cdot\,) \mbox{ and } \tilde{\Phi}_{\nu, \delta_*(R_1)}(v_\nu^\ddag; \cdot\,) \to \Phi_0(v^\ddag; \cdot\,) \mbox{ on $ L^2(\Omega) $,}
\\[1ex]
\mbox{in the sense of $ \Gamma $-convergence, as $ \nu \searrow 0 $. }
\end{array}
\end{equation}

On the basis of these facts, $ \Gamma $-convergence for the sequence $ \{ \Phi_\nu(v_\nu^\ddag; \cdot \, ) \, | \, \nu > 0 \} $ can be verified as follows. 
\medskip

First, to verify the \hyperlink{gamma1}{lower bound}, let us take arbitrary $ \theta^{\ddag} \in L^2(\Omega) $ and $ \{ \theta_\nu^\ddag \, | \, \nu > 0 \} \subset L^2(\Omega) $ such that $ \theta_\nu^\ddag \to \theta^\ddag $ in $ L^2(\Omega) $ as $ \nu \searrow 0 $. Then, using (\ref{GammaLim04})--(\ref{GammaLim03}), we deduce the following inequality:
\begin{equation*}
\Phi_0(v^\ddag; \theta^{\ddag}) \leq \liminf_{\nu \searrow 0} \tilde{\Phi}_{\nu, \delta_1}(v_\nu^\ddag; \theta_\nu^{\ddag}) \leq \liminf_{\nu \searrow 0} \Phi_\nu(v_\nu^\ddag; \theta_\nu^{\ddag}).
\end{equation*}
Thus, the \hyperlink{gamma1}{lower bound} for $ \{ \Phi_\nu(v_\nu^\ddag; \cdot \, ) \} $ is verified. 
\medskip

Second, to verify \hyperlink{gamma2}{optimality}, we take any $ \theta^{\ddag\ddag} \in BV(\Omega) \cap L^2(\Omega) $ $ ({=} D(\Phi_0(v^\ddag; \cdot\,))) $ and use the fact (\ref{GammaLim03}) to take a sequence $ \{ \theta_\nu^{\ddag\ddag} \, | \, \nu > 0 \} \subset H^1(\Omega) $ such that
\begin{equation*}
\theta_\nu^{\ddag\ddag} \to \theta^{\ddag\ddag} \mbox{ in $ L^2(\Omega) $ and $ \tilde{\Phi}_{\nu, \delta_*(R_1)}(v_\nu^\ddag; \theta_\nu^{\ddag\ddag}) \to \Phi_0(v^\ddag; \theta^{\ddag\ddag}) $ as $ \nu \searrow 0 $.}
\end{equation*}
Then, taking into account (\ref{GammaLim04}) and the \hyperlink{gamma1}{lower bound} for $ \{ \Phi_\nu(v_\nu^\ddag; \cdot\,) \} $, we observe that
\begin{equation*}
\limsup_{\nu \searrow 0} \Phi_\nu(v_\nu^\ddag; \theta_\nu^{\ddag\ddag}) \leq \lim_{\nu \searrow 0} \tilde{\Phi}_{\nu, \delta_*(R_1)}(v_\nu^\ddag; \theta_\nu^{\ddag\ddag}) = \Phi_0(v^\ddag; \theta_\nu^{\ddag\ddag}) \leq \liminf_{\nu \searrow 0} \Phi_\nu(v_\nu^\ddag; \theta_\nu^{\ddag\ddag}). 
\end{equation*}
This implies \hyperlink{gamma2}{optimality} for $ \{ \Phi_\nu(v_\nu^\ddag; \cdot \, ) \} $. \hfill \qed
\begin{rem}\label{Rem.Gamma-conv01}
\begin{em}
In previous work \citep[Lemma 3]{MS}, it was reported that strong $ L^2 $-convergence as in $ \mbox{(\hyperlink{(gamma1)_nu}{$\gamma $1})}_\nu $ was necessary to obtain convergence as in (\ref{GammaLim03}). Hence, we still have not succeeded in generalizing the result of Lemma \ref{Lem.GammaLim01} by means of Mosco convergence. 
\end{em}
\end{rem}

In light of Lemma \ref{Lem.GammaLim01} and \citep[Remark 2]{MS}, we find a sequence $ \{ \tilde{\theta}_0^\nu \, | \, \nu > 0 \} \subset H^1(\Omega) $ such that
\begin{equation}\label{Gamma01}
\left\{ ~ \parbox{10cm}{
$ |\tilde{\theta}_0^\nu|_{L^\infty(\Omega)} \leq |\theta_0|_{L^\infty(\Omega)} $ for all $ \nu > 0 $,
\\[1ex]
$ \tilde{\theta}_0^\nu \to \theta_0 $ and $ \Phi_\nu(v_0; \tilde{\theta}_0^\nu) \to \Phi_0(v_0; \theta_0) $ as $ \nu \searrow 0 $.
} \right.
\end{equation}

Now, for arbitrary $ 0 < h < h_*^\dag $ and $ \nu > 0 $, let us denote by $ [\overline{v}_h^\nu, \overline{\theta}_h^\nu] = [\overline{w}_h^\nu, \overline{\eta}_h^\nu, \overline{\theta}_h^\nu] \in L_{\rm loc}^2([0, \infty); L^2(\Omega)^3) $, $ [\underline{v}_h^\nu, \underline{\theta}_h^\nu] = [\underline{w}_h^\nu, \underline{\eta}_h^\nu, \underline{\theta}_h^\nu] \in L_{\rm loc}^2([0, \infty); L^2(\Omega)^3) $, and $ [\widehat{v}_h^\nu, \widehat{\theta}_h^\nu] = [\widehat{w}_h^\nu, \widehat{\eta}_h^\nu, \widehat{\theta}_h^\nu] $ $ \in L_{\rm loc}^2([0, \infty); L^2(\Omega)^3) $ the three time interpolations, as in (\ref{apx200}), consisting of the solution $ \{ [v_i^\nu, \theta_i^\nu] = [w_i^\nu, \eta_i^\nu, \theta_i^\nu] \, | \, i \in \N \} \subset D_1 $ to the approximating problem $ \mbox{(\hyperlink{(AP)_h^nu}{AP})}_h^\nu $, under the assumption (\ref{delta_1}) and the following initial condition: 
\begin{equation}\label{Gamma02}
[v_0^\nu, \theta_0^\nu] = [w_0^\nu, \eta_0^\nu, \theta_0^\nu] := [v_0^\nu, \tilde{\theta}_0^\nu] = [w_0, \eta_0, \tilde{\theta}_0^\nu] \in D_1.
\end{equation}
Then, from (\hyperlink{A1}{A1})--(\hyperlink{A3}{A3}), (\ref{AP_dis01})--(\ref{AP_dis01_01}), (\ref{apx20})--(\ref{apx25}), and (\ref{Gamma01})--(\ref{Gamma02}), it is deduced that
\begin{equation}\label{convF_nu}
\begin{array}{rcl}
F_* & := & \ds \sup_{0 < \nu < \nu_*^\dag} |\mathscr{F}_\nu(v_0^\nu, \theta_0^\nu)| +c_* \mathscr{L}^N(\Omega) 
\\[1ex]
& = & \ds \sup_{0 < \nu < \nu_*^\dag} |\mathscr{F}_\nu(w_0, \eta_0, \tilde{\theta}_0^\nu)| +c_* \mathscr{L}^N(\Omega) < \infty \mbox{ \ for some $ \nu_*^\dag > 0 $,}

\end{array}
\end{equation}
\begin{equation}\label{apx21-2}
\begin{array}{c}
\left\{ \parbox{13.75cm}{
$ \{ \overline{w}_h^\nu(t, x), \underline{w}_h^\nu(t, x), \widehat{w}_h^\nu(t, x) \, | \, 0 < h < h_*^\dag, ~ 0 < \nu < \nu_*^\dag \} \subset [o_*, \iota_*] $,
\\[1ex]
$ \{ \overline{\eta}_h^\nu(t, x), \underline{\eta}_h^\nu(t, x), \widehat{\eta}_h^\nu(t, x) \, | \, 0 < h < h_*^\dag, ~ 0 < \nu < \nu_*^\dag \} \subset [0, 1] $,
\\[0.5ex]
$ \{ \overline{\theta}_h^\nu(t, x), \underline{\theta}_h^\nu(t, x), \widehat{\theta}_h^\nu(t, x) \, | \, 0 < h < h_*^\dag, ~ 0 < \nu < \nu_*^\dag \} \subset [-|\theta_0|_{L^\infty(\Omega)}, |\theta_0|_{L^\infty(\Omega)}] $, 
} \right. 
\\[5ex]
\mbox{a.e.\ $ x \in \Omega $ and any $ t \in [0, T] $,}
\end{array}
\end{equation}
\begin{equation}\label{apx21-3}
\begin{array}{c}
\left\{ \begin{array}{ll}
\multicolumn{2}{l}{\ds {\frac{1}{2}} |(\widehat{v}_h^\nu)_t|_{L^2(0, T; L^2(\Omega))^2}^2 +|{\textstyle \sqrt{\alpha_0(\overline{v}_h^\nu)}} (\widehat{\theta}_h^\nu)_t|_{L^2(0, T; L^2(\Omega))}^2}
\\
& \ds \leq \sum_{{i \in \Z},\,{ih \in [0, T]}} \left( \frac{1}{2h} |v_i^\nu -v_{i -1}^\nu|_{L^2(\Omega)^2}^2 +\frac{1}{h} |{\textstyle \sqrt{\alpha_0(v_i^{\nu})}}(\theta_i^\nu -\theta_{i -1}^\nu)|_{L^2(\Omega)}^2 \right) \leq F_*,
\\[3ex]
\multicolumn{2}{l}{\ds \sup_{t \in [0, T]} |\mathscr{F}_\nu(\overline{v}_h^\nu, \overline{\theta}_h^\nu)| \vee \sup_{t \in [0, T]} |\mathscr{F}_\nu(\underline{v}_h^\nu, \underline{\theta}_h^\nu)| \leq \max_{{i \in \Z},\,{ih \in [0, T]}} |\mathscr{F}_\nu(v_i^\nu, \theta_i^\nu)| \leq F_*,}
\end{array} \right.
\\
\mbox{for any $ 0 < h < h_*^\dag $ and any $ 0 < \nu < \nu_*^\dag $,}
\end{array}
\end{equation}
and, therefore,
\begin{equation}\label{apx20-2}
\left\{ \hspace{-2ex} \parbox{12cm}{
\vspace{-2ex}
\begin{itemize}
\item $ \{ [\widehat{v}_h^\nu, \widehat{\theta}_h^\nu] = [\widehat{w}_h^\nu, \widehat{\eta}_h^\nu, \widehat{\theta}_h^\nu] \, | \, 0 < h < h_*^\dag, ~ 0 < \nu < \nu_*^\dag \} $ is bounded in $ W^{1, 2}(0, T; L^2(\Omega)^3) \cap L^\infty(0, T; X_0) \cap L^\infty(Q)^3 $;
\vspace{-1ex}
\item $ \{ [\overline{v}_h^\nu, \overline{\theta}_h^\nu] = [\overline{w}_h^\nu, \overline{\eta}_h^\nu, \overline{\theta}_h^\nu] \, | \, 0 < h < h_*^\dag, ~ 0 < \nu < \nu_*^\dag \} $ and $ \{ [\underline{v}_h^\nu, \underline{\theta}_h^\nu] = [\underline{w}_h^\nu, \underline{\eta}_h^\nu, \underline{\theta}_h^\nu] \, | \, 0 < h < h_*^\dag, ~ 0 < \nu < \nu_*^\dag \} $ are bounded in $ L^\infty(0, T; L^2(\Omega)^3) \cap L^\infty(0, T; X_0) \cap L^\infty(Q)^3 $.
\vspace{-2ex}
\end{itemize}
} \right.
\end{equation}
Taking into account (\ref{Gamma01})--(\ref{apx20-2}) and Aubin-type compactness theory (see \citep{Simon}), we find sequences
\begin{equation*}
h_*^\dag > h_1 > \cdots > h_n \searrow 0 \mbox{ and } \nu_*^\dag > \nu_1 > \cdots > \nu_n \searrow 0 \mbox{ as $ n \to \infty $}
\end{equation*}
and a triplet of functions $ [v, \theta] = [w, \eta, \theta] \in L^2(0, T; L^2(\Omega)^3) $ such that
\begin{equation}\label{apx25-2}
\left\{ ~ \parbox{7.7cm}{
$ v \in W^{1, 2}(0, T; L^2(\Omega)^2) \cap L^\infty(0, T; H^1(\Omega)^2) $,
\\[1ex]
$ \theta \in W^{1, 2}(0, T; L^2(\Omega)) \cap L^\infty(0, T; BV(\Omega)) $,
} \right.
\end{equation}
\begin{equation}\label{apx22-2}
\begin{array}{c}
\left\{ ~ \parbox{7.7cm}{
$ v(t, x) = [w(t, x), \eta(t, x)] \in [o_*, \iota_*] \times [0, 1] $,
\\[1ex]
$ \theta(t, x) \in [-|\theta_0|_{L^\infty(\Omega)}, |\theta_0|_{L^\infty(\Omega)}] $,
} \right. 
\\[3ex]
\mbox{a.e.\ $ x \in \Omega $ and any $ t \in [0, T] $,}
\end{array}
\end{equation}
and
\begin{equation}\label{apx23-2}
\begin{array}{c}
\left\{ 
\begin{array}{rl}
\multicolumn{2}{l}{\widehat{v}_n = [\widehat{w}_n, \widehat{\eta}_n] := \widehat{v}_{h_n}^{\nu_n} \to v \mbox{ \ in $ C(\overline{I}; L^2(\Omega)^2) $,}} 
\\
& \mbox{weakly in $ W^{1, 2}(I; L^2(\Omega)^2) $, weakly-$ * $ in $ L^\infty(I; H^1(\Omega)^2), $}
\\
& \mbox{and weakly-$ * $ in $ L^\infty(I \times \Omega)^2 $,}
\\[1ex]
\multicolumn{2}{l}{\widehat{\theta}_n := \widehat{\theta}_{h_n^\nu}^{\nu_n} \to \theta \mbox{ \ in $ C(\overline{I}; L^2(\Omega)) $,}}
\\
& \mbox{weakly in $ W^{1, 2}(I; L^2(\Omega)) $, and weakly-$ * $ in $ L^\infty(I \times \Omega) $,}
\end{array}
\right. 
\ \\[-2ex]
\ \\
\mbox{as $ n \to \infty $, \ for any open interval $ I \subset (0, T) $.}
\end{array}
\end{equation}
Additionally, from (\ref{delta_0}) and (\ref{apx300}), we deduced that
\begin{equation}\label{apx24-2}
\begin{array}{c}
\left\{ 
\begin{array}{rl}
\multicolumn{2}{l}{\overline{v}_n = [\overline{w}_n, \overline{\eta}_n] := \overline{v}_{h_n}^{\nu_n} \to v \mbox{ and }\underline{v}_n = [\underline{w}_n, \underline{\eta}_n] := \underline{v}_{h_n}^{\nu_n} \to v} 
\\
& \mbox{in $ L^\infty(I; L^2(\Omega)^2) $, \ \ weakly-$ * $ in $ L^\infty(I; H^1(\Omega)^2) $, \ \  and}
\\
& \mbox{weakly-$ * $ in $ L^\infty(I \times \Omega)^2 $,}
\\[1ex]
\multicolumn{2}{l}{\overline{\theta}_n  := \overline{\theta}_{h_n}^{\nu_n} \to \theta \mbox{ and }\underline{\theta}_n := \underline{\theta}_{h_n}^{\nu_n} \to \theta \mbox{ \ in $ L^\infty(I; L^2(\Omega)) $, and}} 
\\
& \mbox{weakly-$ * $ in $ L^\infty(I \times \Omega) $,}
\end{array}
\right. 
\ \\[-2ex]
\ \\
\mbox{as $ n \to \infty $, \ for any open interval $ I \subset (0, T) $,}
\end{array}
\end{equation}
and, in particular,
\begin{equation}\label{apx24-3}
\begin{array}{c}
\left\{ \hspace{0ex} \parbox{8.5cm}{
\vspace{-2ex}
\begin{itemize}
\item[]\hspace{-4ex}$ \overline{v}_n(t) \to v(t) $, $ \underline{v}_n(t) \to v(t) $ and $ \widehat{v}_n(t) \to v(t) $ in $ L^2(\Omega)^2 $ and weakly in $ H^1(\Omega)^2 $,
\item[]\hspace{-4ex}$ \overline{\theta}_n(t) \to \theta(t) $, $ \underline{\theta}_n(t) \to \theta(t) $ and $ \widehat{\theta}_n(t) \to \theta(t) $ in $ L^2(\Omega) $ and weakly-$ * $ in $ BV(\Omega) $,
\vspace{-2ex}
\end{itemize}
} \right. 
\ \\[-2ex]
\ \\
\mbox{as $ n \to \infty $, \ a.e.\ $ t \in (0, T) $.}
\end{array}
\end{equation}

Based on these, we next check the following lemmas. 

\begin{lem}[cf.\ {\citep[Section 2.2]{MS}}]\label{MS_Lem.468}
Let $ v = [w, \eta] $ and $ \overline{v}_n = [\overline{w}_n, \overline{\eta}_n] $, $ n \in \N $, be pairs of functions as in (\ref{apx25-2})--(\ref{apx24-3}). Also, for any open interval $ I \subset (0, T) $, let $ \Phi_0^I : L^2(I; L^2(\Omega)) $ $ \rightarrow [0, \infty] $ and $ \Phi_n^I : L^2(I; L^2(\Omega)) \rightarrow [0, \infty] $, $ n \in \N $, be functionals defined as
\begin{equation*}
\zeta \in L^2(I; L^2(\Omega)) \mapsto \left\{ \begin{array}{ll}
\multicolumn{2}{l}{\ds \Phi_0^I(\zeta) := \int_I \Phi_0(v(t); \zeta(t)) \, dt,}
\\
& \mbox{if $ \Phi_0(v({}\cdot{}); \zeta({}\cdot{})) \in L^1(I) $,}
\\[1ex]
\infty, & \mbox{otherwise,}
\end{array} \right.
\end{equation*}
and
\begin{equation*}
\zeta \in L^2(I; L^2(\Omega)) \mapsto \Phi_n^I(\zeta) := \int_I \Phi_{\nu_n}(\overline{v}_n; \zeta(t)) \, dt, ~ n \in \N,
\end{equation*}
respectively. Then, for any open interval $ I \subset (0, T) $, the following three statements hold: 
\begin{description}
\item[\textmd{\it (\hypertarget{II-1}{II-1})}]The restriction $ \Phi_0^I|_{C(\overline{I}; L^2(\Omega))} $ is a proper l.s.c.\ and convex function on $ C(\overline{I}; L^2(\Omega)) $ such that $ D(\Phi_0^I|_{C(\overline{I}; L^2(\Omega))}) \supset C(\overline{I}; L^2(\Omega)) \cap L^1(I; BV(\Omega)) $, and also the $ \Phi_n^I $, $ n \in \N $, are proper l.s.c.\ and convex functions on $ L^2(I; L^2(\Omega)) $ such that $ D(\Phi_n^I) = L^2(I; H^1(\Omega)) $, $ n \in \N $. 
\item[\textmd{\it (\hypertarget{II-2}{II-2})}]If $ \zeta^{\dag\dag} \in C(\overline{I}; L^2(\Omega)) $, $ \{ \zeta_n^{\dag\dag} \, | \, n \in \N \} \subset L^2(I; H^1(\Omega)) $, and $ \zeta_n^{\dag\dag}(t) \to \zeta^{\dag\dag}(t) $ in $ L^2(\Omega) $, a.e.\ $ t \in I $, then $ \ds \liminf_{n \to \infty} \Phi_n^I(\zeta_n^{\dag\dag}) \geq \liminf_{n \to \infty} \Phi_0^I(\zeta_n^{\dag\dag}) \geq \Phi_0^I(\zeta^{\dag\dag}) $. 
\item[\textmd{\it (\hypertarget{II-3}{II-3})}]For any $ \zeta^{\ddag\ddag} \in C(\overline{I}; L^2(\Omega)) \cap L^1(I; BV(\Omega)) $, there exists a sequence $ \{ \zeta_n^{\ddag\ddag} \, | \, n \in \N \} \subset C^\infty(\overline{I \times \Omega}) $ such that $ \zeta_n^{\ddag\ddag} \to \zeta^{\ddag\ddag} $ in $ L^2(I; L^2(\Omega)) $ and $ \ds \Phi_n^I(\zeta_n^{\ddag\ddag}) \to \Phi_0^I(\zeta^{\ddag\ddag}) $ as $ n \to \infty $. 
\end{description}
\end{lem}

\noindent
\textbf{Proof.}
From (\hyperlink{A3}{A3}), (\ref{delta_1}), (\ref{delta_0}), and (\ref{apx25-2})--(\ref{apx23-2}), we easily have
\begin{equation}\label{Gamma10}
\left\{ \hspace{-2ex} \parbox{10.5cm}{
\vspace{-2ex}
\begin{itemize}
\item $ \{ \alpha(v), \beta(v) \} \subset C(\overline{I}; L^2(\Omega)) \cap L^\infty(I; H^1(\Omega)) \cap L^\infty(Q) $, $ \delta_0 \leq \alpha(v) \leq \delta_*(1) $, and $ \delta_1 \leq \beta(v) \leq \delta_*(1) $ a.e.\ in $ Q $,
\item $ \{ \alpha(\overline{v}_n), \beta(\overline{v}_n) \, | \, n \in \N \} \subset L^\infty(I; H^1(\Omega)) \cap L^\infty(Q) $, $ \delta_0 \leq \alpha(\overline{v}_n) \leq \delta_*(1) $, and $ \delta_1 \leq \beta(\overline{v}_n) \leq \delta_*(1) $ a.e.\ in $ Q $ for $ n \in \N $,
\item $ \alpha(\overline{v}_n(t)) \to \alpha(v(t)) $ in $ L^2(\Omega) $ and weakly in $ H^1(\Omega) $ as $ n \to \infty $,
\vspace{-2ex}
\end{itemize}
} \right.
\end{equation}
where for any $ R > 0 $, $ \delta_*(R) $ is the constant given in (\ref{delta_beta(R)}). From (\ref{Gamma10}), item (\hyperlink{II-1}{II-1}) is a straightforward consequence of \citep[Lemma 4]{MS}. 
\medskip

Next, let us take $ \zeta^{\dag\dag} \in C(\overline{I}; L^2(\Omega)) $, $ \{ \zeta_n^{\dag\dag} \, | \, n \in \N \} \subset L^2(I; H^1(\Omega)) $, $ \zeta^{\ddag\ddag} \in C(\overline{I}; L^2(\Omega)) \cap L^1(I; BV(\Omega)) $, and $ \{ \zeta_n^{\ddag\ddag} \, | \, n \in \N \} \subset C^\infty(\overline{I \times \Omega}) $, fulfilling the assumptions in (\hyperlink{II-2}{II-2}) and (\hyperlink{II-3}{II-3}). By using the functional given in (\ref{PhiTilde}), we set
\begin{equation*}
\begin{array}{c}
\ds \zeta \in L^2(I; L^2(\Omega)) \mapsto \tilde{\Phi}_{n, \delta}^I(\zeta) := \int_I \tilde{\Phi}_{\nu_n, \delta}(\overline{v}_n(t); \zeta(t)) \, dt \in [0, \infty] 
\\[2ex]
\mbox{for any $ n \in \N $ and any $ \delta > 0 $.}
\end{array}
\end{equation*}
Here, from (\ref{Gamma10}), it is immediately verified that
\begin{equation}\label{Gamma12}
\tilde{\Phi}_{n, \delta_1}^I(\zeta) \leq \Phi_n^I(\zeta) \leq \tilde{\Phi}_{n, \delta_*(1)}^I(\zeta) \mbox{ for any $ z \in L^2(I; L^2(\Omega)) $.}
\end{equation}
Furthermore, by virtue of (\hyperlink{A3}{A3}), (\ref{delta_0}), (\ref{GammaLim04}), (\ref{apx21-2}), (\ref{apx22-2}), (\ref{apx24-3}), and \citep[Lemma 8]{MS}, we find a sequence $ \{ \zeta_n^{\ddag\ddag} \, | \, n \in \N \} \subset C^\infty(\overline{I \times \Omega}) $ such that
\begin{equation}\label{Gamma11}
\zeta_n^{\ddag\ddag} \to \zeta^{\ddag\ddag} \mbox{ in
$ L^2(I; L^2(\Omega)) $ and $ \tilde{\Phi}_{n, \delta}(\zeta_n^{\ddag\ddag}) \to \Phi_0(\zeta^{\ddag\ddag}) $ as $ n \to \infty $ for any $ \delta > 0 $.}
\end{equation}

Now, taking into account (\ref{Gamma12})--(\ref{Gamma11}) and applying \citep[Lemma 6]{MS}, we can verify the remaining items (\hyperlink{II-2}{II-2}) and (\hyperlink{II-3}{II-3}), as follows:
\begin{equation*}
\left\{ ~\parbox{12cm}{
$ \ds \liminf_{n \to \infty} \Phi_n^I(\zeta_n^{\dag\dag}) \geq \liminf_{n \to \infty} \tilde{\Phi}_{n, \delta_1}^I(\zeta_n^{\dag\dag}) \geq \Phi_0(\zeta^{\dag\dag}) $,
\\[1ex]
$ \ds \limsup_{n \to \infty} \Phi_n^I(\zeta_n^{\ddag\ddag}) \leq \lim_{n \to \infty} \tilde{\Phi}_{n, \delta_*(1)}^I(\zeta_n^{\ddag\ddag}) = \Phi_0^I(\zeta^{\ddag\ddag}) \leq \liminf_{n \to \infty} \Phi_n^I(\zeta_n^{\ddag\ddag}) $. 
}\right.
\end{equation*}
\qed
\bigskip

\noindent
\textbf{Proof of Main Theorem \ref{MainTh02}.}
From (\ref{apx25-2})--(\ref{apx22-2}), it can be seen that the limiting triplet $ [v, \theta] = [w, \eta, \theta] $ fulfills the condition $ \mbox{(\hyperlink{(S0)_0}{S0})}_0 $. Hence, all we have to do is verify the compatibility of $ [v, \theta] $ with conditions $ \mbox{(\hyperlink{(S1)_0}{S1})}_0 $--$ \mbox{(\hyperlink{(S3)_0}{S3})}_0 $.
\medskip

Fix any open interval $ I \subset (0, T) $. Then, from (\ref{apx01-02})--(\ref{apx03}) and (\ref{apx200}), the functions $ [\overline{v}_n, \overline{\theta}_n] $, $ [\underline{v}_n, \underline{\theta}_n] $, and $ [\widehat{v}_n, \widehat{\theta}_n] $, $ n \in \N $, must fulfill

\begin{equation}\label{MT2-10}
\begin{array}{rl}
\multicolumn{2}{l}{\ds \int_I \left( \rule{0pt}{10pt} (\widehat{v}_n)_t(t), \overline{v}_n(t) -\varpi \right)_{L^2(\Omega)^2} \, dt  +\int_I \left( \rule{0pt}{10pt} [\nabla g](\overline{v}_n(t)), \overline{v}_n(t) -\varpi \right)_{L^2(\Omega)^2} \, dt}
\\[2ex]
\qquad \quad & \ds +\int_I \left( \rule{0pt}{10pt} \nabla \overline{v}_n(t), \nabla (\overline{v}_n(t) -\varpi) \right)_{L^2(\Omega)^{2 \times N}} \, dt
\\[2ex]
\qquad \quad & \ds +\int_I \int_\Omega \left( \rule{0pt}{10pt} |\nabla \underline{\theta}_n(t)| [\nabla \alpha](\overline{v}_n(t)) +\nu |\nabla \underline{\theta}_n|^2 [\nabla \beta](\overline{v}_n(t)) \right) \cdot (\overline{v}_n(t) -\varpi) \, dx \, dt
\\[2ex]
& \ds  +\int_I \Gamma(\overline{v}_n(t)) \, dt \leq \ds \int_I \Gamma(\varpi) \, dt 
\\[2ex]
\multicolumn{2}{l}{\mbox{ for any $ \varpi \in [H^1(\Omega) \cap L^\infty(\Omega)]^2 $ and any $ n \in \N $}}
\end{array}
\end{equation}
and
\begin{equation}\label{MT2-11}
\begin{array}{rl}
\multicolumn{2}{l}{\ds \int_I \bigl( \alpha_0(\overline{v}_n(t)) (\widehat{\theta}_n)_t(t), \overline{\theta}_n(t) -\zeta(t) \bigr)_{L^2(\Omega)} \, dt +\Phi_n^I(\overline{\theta}_n) \leq \Phi_n^I(\zeta)} 
\\[2ex]
\qquad \qquad \qquad  & \mbox{ for any $ \zeta \in L^2(I; H^1(\Omega)) $ and $ n \in \N $.}
\end{array}
\end{equation}

On this basis, we next take the limit of (\ref{MT2-11}) as $ n \to \infty $. Then, invoking (\hyperlink{A3}{A3}), (\ref{apx23-2})--(\ref{apx24-3}), and (\hyperlink{II-2}{II-2}) of Lemma \ref{MS_Lem.468}, we calculate that
\begin{eqnarray*}
&& \hspace{-7ex} \int_I \bigl( \alpha_0(v(t)) \theta_t(t), \theta(t) -\zeta(t) \bigr)_{L^2(\Omega)} \, dt +\Phi_0^I(\theta)
\\
& \leq & \lim_{n \to \infty} \int_I \bigl( \alpha_0(\overline{v}_n) (\widehat{\theta}_n)_t(t), \overline{\theta}_n(t) -\zeta(t) \bigr)_{L^2(\Omega)} \, dt +\liminf_{n \to \infty} \Phi_n^I(\overline{\theta}_n)
\\
& \leq & \lim_{n \to \infty} \Phi_n^I(\zeta) = \Phi_0^I(\zeta) \mbox{ \ for any $ \zeta \in L^2(I; H^1(\Omega)) $.} 
\end{eqnarray*}
Since the choice of the open interval $ I \subset (0, T) $ is arbitrary, this inequality implies that
\begin{equation}\label{MT2-01}
\begin{array}{c}
\ds \bigl( \alpha_0(v(t)) \theta_t(t), \theta(t) -\omega \bigr)_{L^2(\Omega)} +\Phi_0(v(t); \theta(t)) \leq \Phi_0(v(t); \omega) 
\\[1ex]
\mbox{for any $ \omega \in H^1(\Omega) $.}
\end{array}
\end{equation}
Moreover, by virtue of the strict approximation given by (\hyperlink{Fact5}{Fact 5}) in Remark \ref{Rem.Note06}, we can verify that (\ref{MT2-01}) is valid for any $ \omega \in BV(\Omega) \cap L^2(\Omega) $, and the triplet $ [v, \theta] = [w, \eta, \theta] $ is compatible with condition $ \mbox{(\hyperlink{(S3)_0}{S3})}_0 $. 
\bigskip

Next, with (\ref{apx25-2}) in mind, we can apply (\hyperlink{II-3}{II-3}) of Lemma \ref{MS_Lem.468} for the case in which $ \zeta^{\ddag\ddag} = \theta $ and take a sequence $ \{ \tilde{\theta}_n \, | \, n \in \N \} \subset C^\infty(\overline{I \times \Omega}) $ such that
\begin{equation}\label{MT2-02}
\tilde{\theta}_n \to \theta \mbox{ in $ L^2(I; L^2(\Omega)) $ and } \Phi_n^I(\tilde{\theta}_n) \to \Phi_0^I(\theta) \mbox{ as $ n \to \infty $.}
\end{equation}
Setting $ \zeta = \tilde{\theta}_n $ in (\ref{MT2-11}), from (\ref{apx23-2})--(\ref{apx24-3}), (\ref{Gamma10}), (\ref{MT2-02}), and (\hyperlink{II-2}{II-2}) of Lemma \ref{MS_Lem.468}, we have

\begin{equation*}
\begin{array}{rcl}
&& \ds \Phi_0^I(\theta) \leq \ds \liminf_{n \to \infty} \Phi_0^I(\overline{\theta}_n) \leq \liminf_{n \to \infty} \Phi_n^I(\overline{\theta}_n) \leq \limsup_{n \to \infty} \Phi_n^I(\overline{\theta}_n) 
\\[2ex]
& \leq & \ds \lim_{n \to \infty} \left[ \Phi_n^I(\tilde{\theta}_n) -\int_I \bigl( \alpha_0(\overline{v}_n(t)) (\widehat{\theta}_n)_t(t), \overline{\theta}_n(t) -\tilde{\theta}_n(t) \bigr)_{L^2(\Omega)}\, dt \right]
= \Phi_0^I(\theta), 
\end{array}
\end{equation*}

\noindent
and therefore

\begin{equation}\label{MT2-05}
\begin{array}{rcl}
0 & \leq & \ds \liminf_{n \to \infty} \nu_n \int_I \int_\Omega \beta(\overline{v}_n(t)) |\nabla \overline{\theta}_n(t)|^2 \, dx \, dt 
\\[2ex]
& \leq & \ds \limsup_{n \to \infty} \nu_n \int_I \int_\Omega \beta(\overline{v}_n(t)) |\nabla \overline{\theta}_n(t)|^2 \, dx \, dt 
\\[2ex]& \leq & \ds \lim_{n \to \infty} \Phi_n^I(\overline{\theta}_n) -\liminf_{n \to \infty} \Phi_0^I(\overline{\theta}_n) = 0
\end{array}
\end{equation}

\noindent
and

\begin{equation}\label{MT2-04}
\lim_{n \to \infty} \int_I \int_\Omega \alpha(\overline{v}_n(t)) |\nabla \overline{\theta}_n(t)| \, dx \, dt = \int_I \int_\Omega d[\alpha({v}(t)) |D {\theta}(t)|].
\end{equation}
\medskip

\noindent
Again, keeping in mind (\ref{apx23-2})--(\ref{apx24-3}), (\ref{Gamma10}), and (\ref{MT2-04}), we apply Lemma \ref{MS_Lem.7} for the case in which

\begin{equation*}
\parbox{12.75cm}{$ \varrho = \alpha(v) $, $ \{ \varrho_n \} = \{ \alpha(\overline{v}_n) \} $, $ \zeta = \theta $, $ \{ \zeta_n \} = \{ \overline{\theta}_n \} $, $ \omega = 1 $, and $ \{ \omega_n \} = \{ 1 \} $.}
\end{equation*}

\noindent
We then have

\begin{equation*}
\int_I \int_\Omega |\nabla \overline{\theta}_n(t)| \, dx \, dt \to \int_I \int_\Omega |D \theta(t)| \, dt \mbox{ as $ n \to \infty $.}
\end{equation*}
\medskip

Conversely, as a result of (\ref{delta_0}), (\ref{apx200}) and (\ref{apx21-3}),

\begin{equation*}
\left\{ \begin{array}{l}
\ds |\overline{\theta}_n -\underline{\theta}_n|_{L^\infty(0, T; L^2(\Omega))} \leq \sqrt{h_n} |(\widehat{\theta}_n)_t|_{L^2(0, T; L^2(\Omega))} \to 0
\\[1ex]
\ds \left| \int_I \int_\Omega |\nabla \overline{\theta}_n| \, dx \, dt -\int_I \int_\Omega |\nabla \underline{\theta}_n| \, dx \, dt \right| \leq \frac{2 F_*}{\delta_0} h_n \to 0
\end{array} \right.\mbox{as $ n \to \infty $.}
\end{equation*}
\medskip

\noindent
Therefore, taking any $ \varpi \in [H^1(\Omega) \cap L^\infty(\Omega)]^2 $ and applying Lemma \ref{MS_Lem.7} for the case in which

\begin{equation*}
\parbox{15cm}{$ \varrho = 1 $, $ \{ \varrho_n \} = \{ 1 \} $, $ \zeta = \theta $, $ \{ \zeta_n \} = \{ \underline{\theta}_n \} $, $ \omega = \varpi \cdot [\nabla \alpha](v) $, and $ \{ \omega_n \} = \{ \varpi \cdot [\nabla \alpha](\overline{v}_n) \} $,}
\end{equation*}

\noindent
we obtain

\begin{equation}\label{MT2-06}
\begin{array}{c}
\ds \lim_{n \to \infty} \int_I \int_\Omega \varpi \cdot [\nabla \alpha](\overline{v}_n(t))|\nabla \underline{\theta}_n(t)| \, dx \, dt = \int_I \int_\Omega d[\varpi \cdot [\nabla \alpha](v(t))|D \theta(t)|]
\\[2ex]
\mbox{for any $ \varpi \in [H^1(\Omega) \cap L^\infty(\Omega)]^2 $.}
\end{array}
\end{equation}
Furthermore, from (\ref{delta_1}), (\ref{MT2-07}), and (\ref{MT2-05}), we can compute the following:
\begin{eqnarray}
0 & \leq & \ds \liminf_{n \to \infty} \left| \nu_n \int_I \int_\Omega \varpi \cdot [\nabla \beta](\overline{v}_n(t))|\nabla \underline{\theta}_n(t)|^2 \, dx \, dt \right|
\nonumber
\\
& \leq & \ds \limsup_{n \to \infty} \left| \nu_n \int_I \int_\Omega \varpi \cdot [\nabla \beta](\overline{v}_n(t))|\nabla \underline{\theta}_n(t)|^2 \, dx \, dt \right|
\nonumber
\\
& \leq & \ds |\varpi|_{L^\infty(\Omega)^2} |\beta|_{C^1([0, 1]^2)} \limsup_{n \to \infty} \left| \nu_n \int_I \int_\Omega |\nabla \underline{\theta}_n(t)|^2 \, dx \, dt \right|
\label{MT2-08}
\\
& = & \ds |\varpi|_{L^\infty(\Omega)^2} |\beta|_{C^1([0, 1]^2)} \limsup_{n \to \infty} \left| \nu_n \int_I \int_\Omega |\nabla \overline{\theta}_n(t)|^2 \, dx \, dt \right|
\nonumber
\\
& \leq & \ds \frac{|\varpi|_{L^\infty(\Omega)^2} |\beta|_{C^1([0, 1]^2)}}{\delta_1} \lim_{n \to \infty} \left| \nu_n \int_I \int_\Omega \beta(\overline{v}_n(t))|\nabla \overline{\theta}_n(t)|^2 \, dx \, dt \right|
\nonumber
\\
&  = & 0 \mbox{ \ for any $ \varpi \in [H^1(\Omega) \cap L^\infty(\Omega)]^2 $.}
\nonumber
\end{eqnarray}

Thanks to (\ref{Gamma(v)}), (\ref{apx23-2})--(\ref{apx24-3}), and (\ref{MT2-06})--(\ref{MT2-08}), allowing $ n \to \infty $ in (\ref{MT2-10}) yields
\begin{equation}\label{MT2-09}
\begin{array}{rl}
\multicolumn{2}{l}{\ds \int_I \left( \rule{0pt}{10pt} v_t(t), v(t) -\varpi \right)_{L^2(\Omega)^2} \, dt  +\int_I \left( \rule{0pt}{10pt} [\nabla g](v(t)), v(t) -\varpi \right)_{L^2(\Omega)^2} \, dt}
\\[2ex]
\qquad \quad & \ds +\int_I \left( \rule{0pt}{10pt} \nabla v(t), \nabla (v(t) -\varpi) \right)_{L^2(\Omega)^{2 \times N}} \, dt
\\[2ex]
\qquad \quad & \ds +\int_I \int_\Omega d \bigl[ (v(t) -\varpi) \cdot [\nabla \alpha](v(t)) |D \theta(t)| \bigr] \, dt
\\[2ex]
& \ds  +\int_I \Gamma(v(t)) \, dt \leq \ds \int_I \Gamma(\varpi) \, dt \mbox{ \ for any $ \varpi \in [H^1(\Omega) \cap L^\infty(\Omega)]^2 $.}
\end{array}
\end{equation}
Since the open interval $ I \subset (0, T) $ is arbitrary, from (\ref{MT2-09}) and the reformulation (\ref{unify02}) in Remark \ref{Rem.1.0}, the triplet $ [v, \theta] = [w, \eta, \theta] $ fulfills conditions $ \mbox{(\hyperlink{(S1)_0}{S1})}_0 $--$ \mbox{(\hyperlink{(S2)_0}{S2})}_0 $. \hfill \qed 

\begin{rem}\label{Rem.finalComment}
\begin{em}
From the proof of Main Theorem \ref{MainTh02}, it can be said that the convex function $ \beta \in C^1(\R^2) \cap C^2([0, 1]^2) $ can be one of the approximation components for (\hyperlink{(S;u)_0}{S})$_0$.
\end{em}
\end{rem}

\section{Appendix}
\ \vspace{-3ex}

In this Appendix, we make supplementary statements for some preliminary facts as in Section 1, and the solutions to our systems. 
\medskip

First, we show the proof of (\hyperlink{Fact1}{Fact\,1}) in Remark \ref{Rem.Time-dep.}, because the results in the reference \cite[Chapter 2]{Kenmochi} were discussed under quite general settings, and these might not be to-the-point under our simplified setting.
\medskip

\noindent
\textbf{Proof of (Fact\,1).} 
Let us assume that $ [\zeta, \zeta^*] \in \partial \Psi^I $ in $ L^2(I; X)^2 $, i.e.:
\begin{equation}\label{ap01}
\int_I \bigl( \zeta^*(t), \xi(t) -\zeta(t) \bigr)_X \, dt \leq \Psi^I(\xi) -\Psi^I(\zeta), \mbox{ for any $ \xi \in L^2(I; X) $.}
\end{equation}
Also, let us take any open interval $ A \subset I $, and denote by $ \chi_A : \R \to \{0, 1\} $ the characteristic function  of $ A $. On this basis, we take any $ z \in E_0 $, and set
\begin{equation*}
\xi(t) := \zeta(t) +\chi_{A}(t) (z -\zeta(t)) \mbox{ in $ X $, for a.e. $ t \in I $,}
\end{equation*}
as the test function $ \xi \in L^2(I; X) $ in (\ref{ap01}). Then, it is computed that
\begin{equation}\label{ap02}
\begin{array}{rl}
\multicolumn{2}{l}{\ds \int_{A} \bigl( \zeta^*(t), z -\zeta(t) \bigr)_X \, dt \leq \Psi^I(\xi) -\Psi^I(\zeta)}
\\[2ex]
& \ds \leq \int_A \Psi^t(z) \, dt +\int_{I \setminus A} \Psi^t(\zeta(t)) -\int_I \Psi^t(\zeta(t)) \, dt
\\[2ex]
& \ds = \int_A \Psi^t(z) \, dt -\int_A \Psi^t(\zeta(t)) \, dt,
\\[2ex]
\multicolumn{2}{l}{\mbox{for any $ z \in E_0 $ and any open interval $ A \subset I $.}}
\end{array}
\end{equation}

Now, on account of the assumptions for $ E_0 $ $ (= D(\Psi^t)) $, we can see from (\ref{ap02}) that
\begin{equation}\label{ap03}
[\zeta(t), \zeta^*(t)] \in \partial \Psi^t \mbox{ in $ X^2 $, a.e. $ t \in I $.}
\end{equation}

Conversely, if we suppose (\ref{ap03}) for a pair $ [\zeta, \zeta^*] \in D(\Psi^I) \times L^2(I; X) $, then  we immediately derive (\ref{ap01}) as a straightforward consequence of the definition of subdifferential. \hfill \qed
\bigskip

Next, we briefly see the demonstration scenario of (\hyperlink{Fact7}{Fact\,7}) in Remark \ref{Rem.MG}.
\bigskip

\noindent
\textbf{Proof of (\hyperlink{Fact7}{Fact\,7}).} 
First, let us take any $ \tilde{z} \in D(\Psi) $ with a sequence $ \{ \tilde{z}_n \, | \, n \in \N \} \subset X $ such that $ \tilde{z}_n \in D(\Psi_n) $ for any $ n \in \N $, and $ \tilde{z}_n \to \tilde{z} $ in $ X $ and $ \Psi_n(\tilde{z}_n) \to \Psi(\tilde{z}) $ as $ n \to \infty $. Then, with the assumptions in mind, we compute that:
\begin{equation}\label{ap04}
\begin{array}{rl}
(z^*, \tilde{z} -z)_X +\Psi(z) & \leq \ds \lim_{n \to \infty} (z_n^*, \tilde{z}_n -z_n)_X +\liminf_{n \to \infty} \Psi_n(z_n)
\\[2ex]
& \ds \leq \limsup_{n \to \infty} \bigl[ (z_n^*, \tilde{z}_n -z_n)_X +\Psi_n(z_n) \bigr]
\\[2ex]
& \ds \leq \lim_{n \to \infty} \Psi_n(\tilde{z}_n) = \Psi(\tilde{z}), \mbox{ for any $ \tilde{z} \in D(\Psi) $.}
\end{array}
\end{equation}
Hence, $ [z, z^*] \in \partial \Psi $ in $ X^2 $. Moreover, if we consider (\ref{ap04}) anew, by setting $ \tilde{z} = z $, then we can show that
\begin{equation}\label{ap05}
\limsup_{n \to \infty} \Psi_n(z_n) \leq \lim_{n \to \infty} \bigl[ \Psi(\tilde{z}_n) -(z_n^*, \tilde{z}_n -z_n)_X \bigr] = \Psi(z).
\end{equation}
The lower bound condition and (\ref{ap05}) lead to the convergence $ \lim_{n \to \infty} \Psi_n(z_n) = \Psi(z) $, immediately. \hfill \qed
\bigskip

Finally, we leave the following remark as a further observation for the future works.
\begin{rem}[Energy estimate]\label{Rem.EnergyEst}
\begin{em}
For any $ \nu > 0 $, let $ [v_\nu, \theta_\nu] = [w_\nu, \eta_\nu, \theta_\nu] $ be the solution to (\hyperlink{(S;u)_0}{S})$_\nu$ obtained as the limits as in (\ref{apx25})-(\ref{*end}). On this basis, let us consider the inequality (\ref{AP_dis01}) in the cases when $ s = 0 $ and $ h_n^\nu $, for $ n \in \N $, and take the limit as $ n \to \infty $. Then, it is derived from (\ref{apx25})-(\ref{*end}) and Lemma \ref{Lem.Mosco01} that
\begin{equation}\label{ap10}
\begin{array}{c}
\ds \frac{1}{2} \int_0^t |(v_\nu)_t(\tau)|_{L^2(\Omega)^2}^2 \, d \tau +\int_0^t |{\textstyle \sqrt{\alpha_0(v_\nu(\tau))}} (\theta_\nu)_t(\tau)|_{L^2(\Omega)}^2 \, d \tau
\\[2ex]
+\mathscr{F}_\nu(v_\nu(t), \theta_\nu(t)) \leq \mathscr{F}_\nu(v_0^\nu, \theta_0^\nu), \mbox{ for all $ t \in [0, T] $.}
\end{array}
\end{equation}
Also, in the case of $ \nu = 0 $, we can see the following kindred inequality to (\ref{ap10}):
\begin{equation}\label{ap11}
\begin{array}{c}
\ds \frac{1}{2} \int_0^t |v_t(\tau)|_{L^2(\Omega)^2}^2 \, d \tau +\int_0^t |{\textstyle \sqrt{\alpha_0(v(\tau))}} \theta_t(\tau)|_{L^2(\Omega)}^2 \, d \tau
\\[2ex]
+\mathscr{F}_0(v(t), \theta(t)) \leq \mathscr{F}_0(v_0, \theta_0), \mbox{ for all $ t \in [0, T] $,}
\end{array}
\end{equation}
by putting $ s = 0 $, $ \nu = \nu_n $ and $ h = h_n $ with $ n  \in \N $ in (\ref{AP_dis01}), and letting $ n \to \infty $ with (\ref{Gamma01})-(\ref{Gamma02}), (\ref{apx23-2})-(\ref{apx24-2}) and Lemma \ref{Lem.GammaLim01} in mind. Here, $ [v, \theta] = [w, \eta, \theta] $ in (\ref{ap11}) is the solution to (\hyperlink{(S;u)_0}{S})$_0$ obtained as the limits as in (\ref{apx25-2})-(\ref{apx24-3}).

Note that the above inequalities (\ref{ap10})-(\ref{ap11}) provide energy estimates for the observations of time-global solutions. In fact, thanks to these inequalities, we may suppose that
\begin{equation*}
[v_t, \theta_t] \in L^2(0, \infty; L^2(\Omega)^3), ~ [v, \theta] \in L^\infty(0, \infty; X_0) \mbox{ and } \nu \theta \in L^\infty(0, \infty; H^1(\Omega)),
\end{equation*} 
for the time-global solution $ [v, \theta] = [w, \eta, \theta] $ to (\hyperlink{(S;u)_0}{S})$_\nu$, in any case of $ \nu \geq 0 $. 
\end{em}
\end{rem}


\end{document}